\documentclass[11pt]{amsart}
\usepackage{amsmath,amssymb,amsthm}
\usepackage{tikz}
\usepackage{subcaption}

\usepackage[pagewise]{lineno}

\newtheorem{theorem}{Theorem} 
\newtheorem{lemma}{Lemma}     
\newtheorem{definition}{Definition}

\newtheorem{example}{Example}
\newtheorem{remark}{Remark}
\newtheorem{question}{Question}
\newtheorem*{theorem*}{Theorem}

\begin{document}

\title{A survey on the Topology of  Fractal Squares}
\date{\today}
\author[J. Luo]{Jun Luo}
\address{School of Mathematics\\
    Sun Yat-Sen University\\ Guangzhou 510275, China}
\email{luojun3@mail.sysu.edu.cn}
\author[H. Rao]{Hui Rao}
\address{Department of Mathematics\\ Central China Normal University\\ Wuhan 430079, China}
\email{hrao@mail.ccnu.edu.cn}

\begin{abstract} 
We consider a special type of self-similar sets, called  fractal squares, and give a brief review on recent results and unsolved issues with an emphasis on their topological properties.
\end{abstract}

\maketitle

\tableofcontents

\section{Introduction}

A map $\varphi:\mathbb{R}^d\rightarrow\mathbb{R}^d$ is called a contraction if there is $c\in(0,1)$ such that the Euclidean distance $\|\varphi(x)-\varphi(y)\|$ between $\varphi(x)$ and $\varphi(y)$ is less than $c\|x-y\|$. A family $\mathcal{F}=\{\varphi_j: 1\le j\le q\}$ of $q\ge2$  contractions $\varphi_j:\mathbb{R}^d\rightarrow\mathbb{R}^d$ is called an iterated function system (shortly, an IFS).    Due  to Hutchinson \cite{Hutchinson81},  
for any $\mathcal{F}$ there is a unique nonempty compact set $K$ satisfying $K=\bigcup_j\varphi_j(K)$. For the sake of convenience, we set  
$\Phi(X)=\bigcup_j\varphi_j(X)$ for nonempty compact sets $X\subset\mathbb{R}^d$. Moreover, set $\Phi^n(X)=\Phi\Big(\Phi^{n-1}(X)\Big)$ for $n\ge2$.

\begin{definition}\label{def:attractor}
We  call $\Phi$ the Hutchinson map of $\mathcal{F}$ and  $K$ the attractor of $\mathcal{F}$. If further every $\varphi_j$ is a similitude, we call  $K$  a {\em self-similar set}. 
\end{definition}

This paper studies self-similar sets in the plane. In deed, we mostly consider contractions $\varphi_j\ (1\le j\le q)$ that are of the form $f_{d,N}(x)=\frac{x+d}{N}$ for some $N\ge2$ and  $d\in \mathbb{Z}^2$. 
When $N$ is fixed or clear from the context, we just write $f_d$ instead of $f_{d,N}$. We also consider $d\in\mathbb{Z}^2$ as a vector in $\mathbb{R}^2$ so that every translation $x\mapsto x+d$ 
  is well defined. 

\begin{definition}\label{def:FS}
Given $N\ge2$, a {\em fractal square of order $N$} means the unique nonempty compact set satisfying $\displaystyle K=\bigcup_{d\in\mathcal{D}}f_{d}(K)$ for some nonempty $\mathcal{D}\subset\{0,\ldots,N-1\}^2$. 
\end{definition}

Hereafter, we often write $K_\mathcal{D}$ instead of $K$ to emphasize the digit set $\mathcal{D}$. We also write $K(N,\mathcal{D})$, when we want to emphasize both the order $N$ and the digit set $\mathcal{D}$. 
Set $\mathcal{D}_1=\mathcal{D}$ and for $j\ge1$ further set   
\begin{eqnarray}\label{eq:approx}
\mathcal{D}_{j+1}&=&N \mathcal{D}_j+\mathcal{D}=\left\{Nd_1+d_2: d_1\in\mathcal{D}_j, d_2\in\mathcal{D}\right\}
\\
K^{(j)}&=&\bigcup_{d\in \mathcal{D}_j}\frac{[0,1]^2+d}{N^j}.
\end{eqnarray}
Then $\left\{K^{(j)}: j\ge1\right\}$ is a decreasing sequence and the limit is just $K$.
\begin{definition}\label{def:K^j}
We call $K^{(j)}$  the {\em $j$-th approximation} of $K$. 
\end{definition}

\begin{remark} Given  $\mathcal{D}\subset\{0,\ldots,N-1\}^2$,  $K(N,\mathcal{D})=K(N^j,\mathcal{D}_j)$ for all $j\ge2$.  \end{remark}

\begin{example}\label{exmp:sierpinski}
 Given  $\mathcal{D}_S=\left\{0,1,2\right\}^2\setminus\{(1,1)\}$, the fractal square $K\left(3,\mathcal{D}_S\right)$ is known as the {\em Sierpi\'nski carpet}.  Figure \ref{fig:sierpinski} illustrates the approximations $K^{(j)}$ for $1\le j\le3$.
\begin{figure}[ht]    
\begin{center}
\begin{tabular}{cccc}
\includegraphics[width=2.75cm]{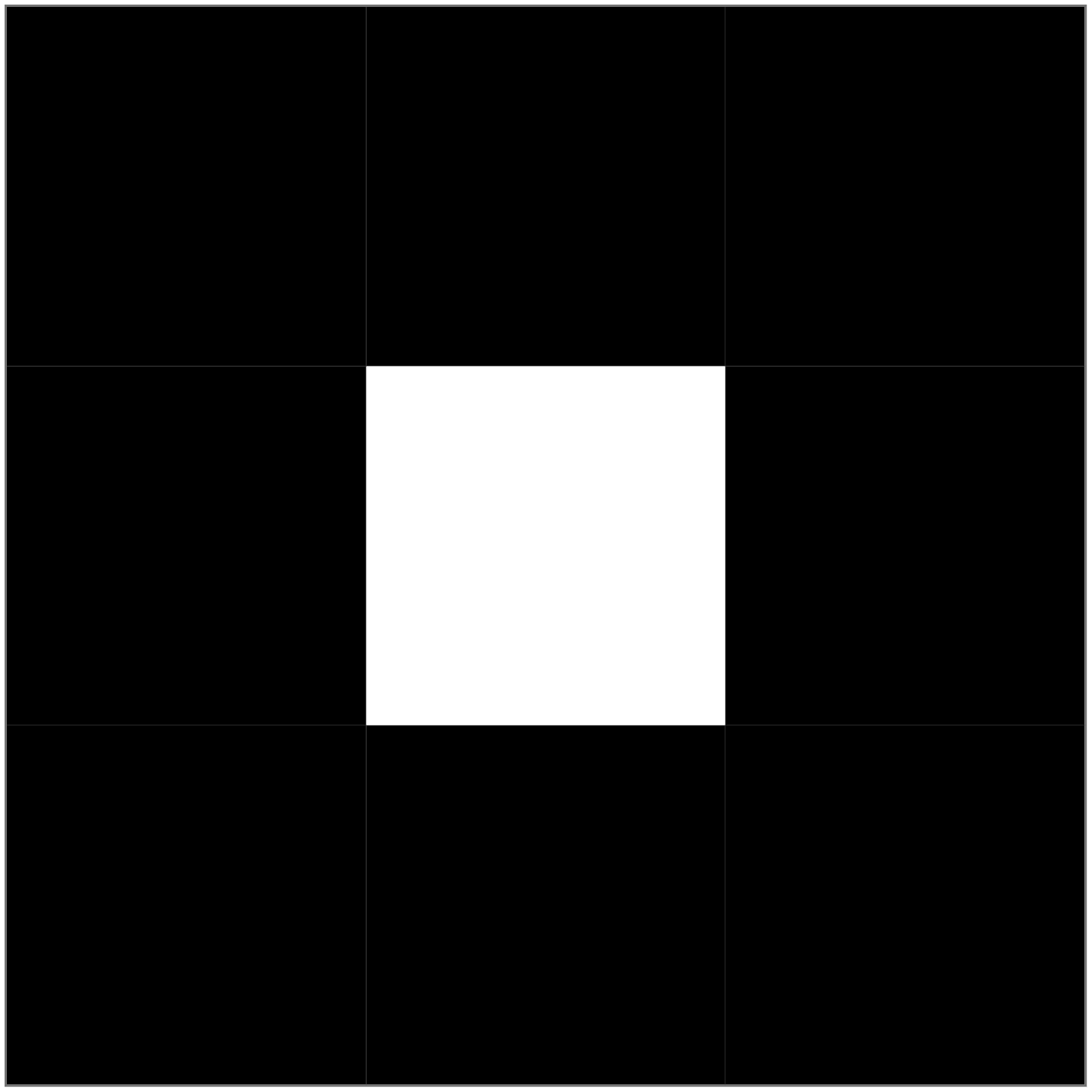} 
&
\includegraphics[width=2.75cm]{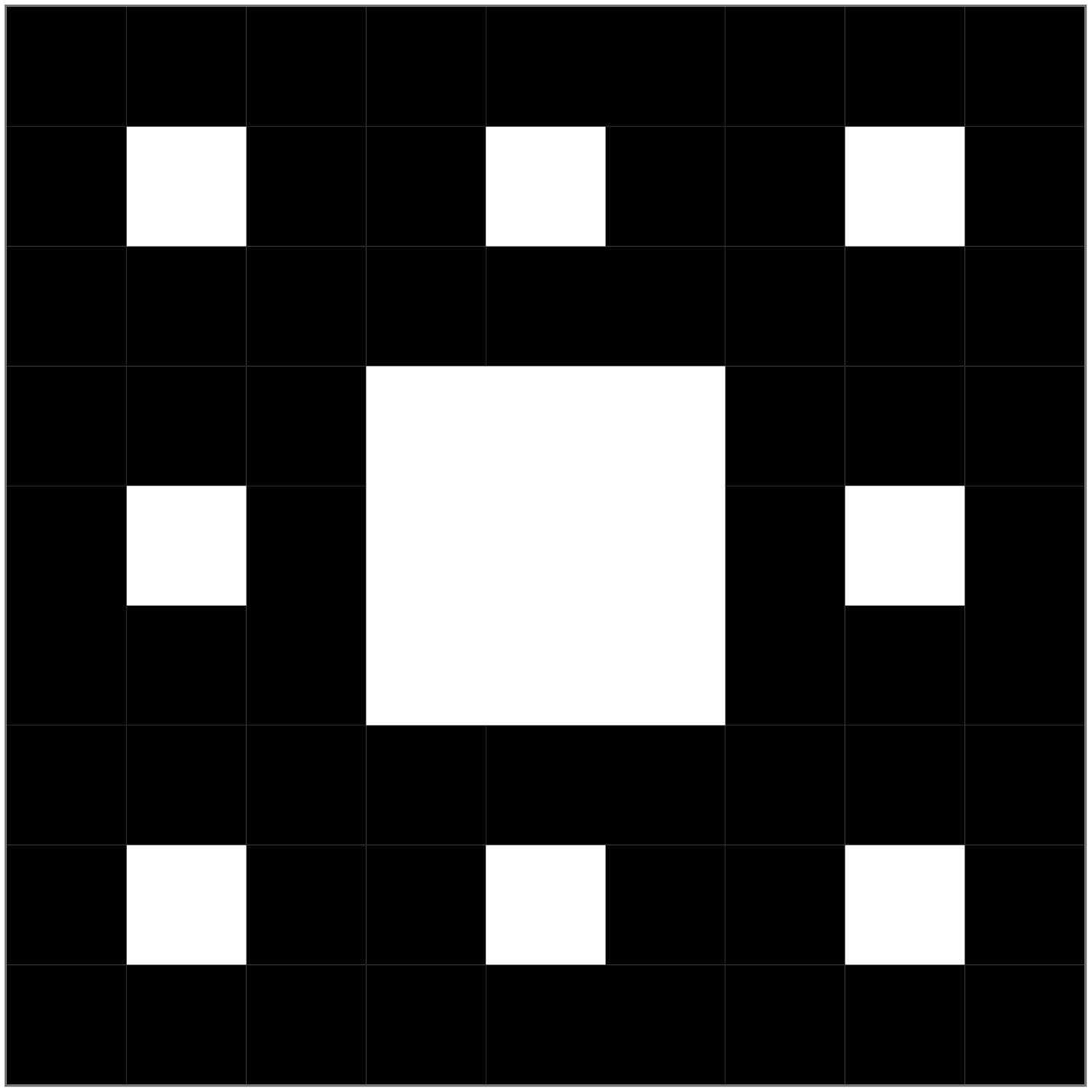} 
&
\includegraphics[width=2.75cm]{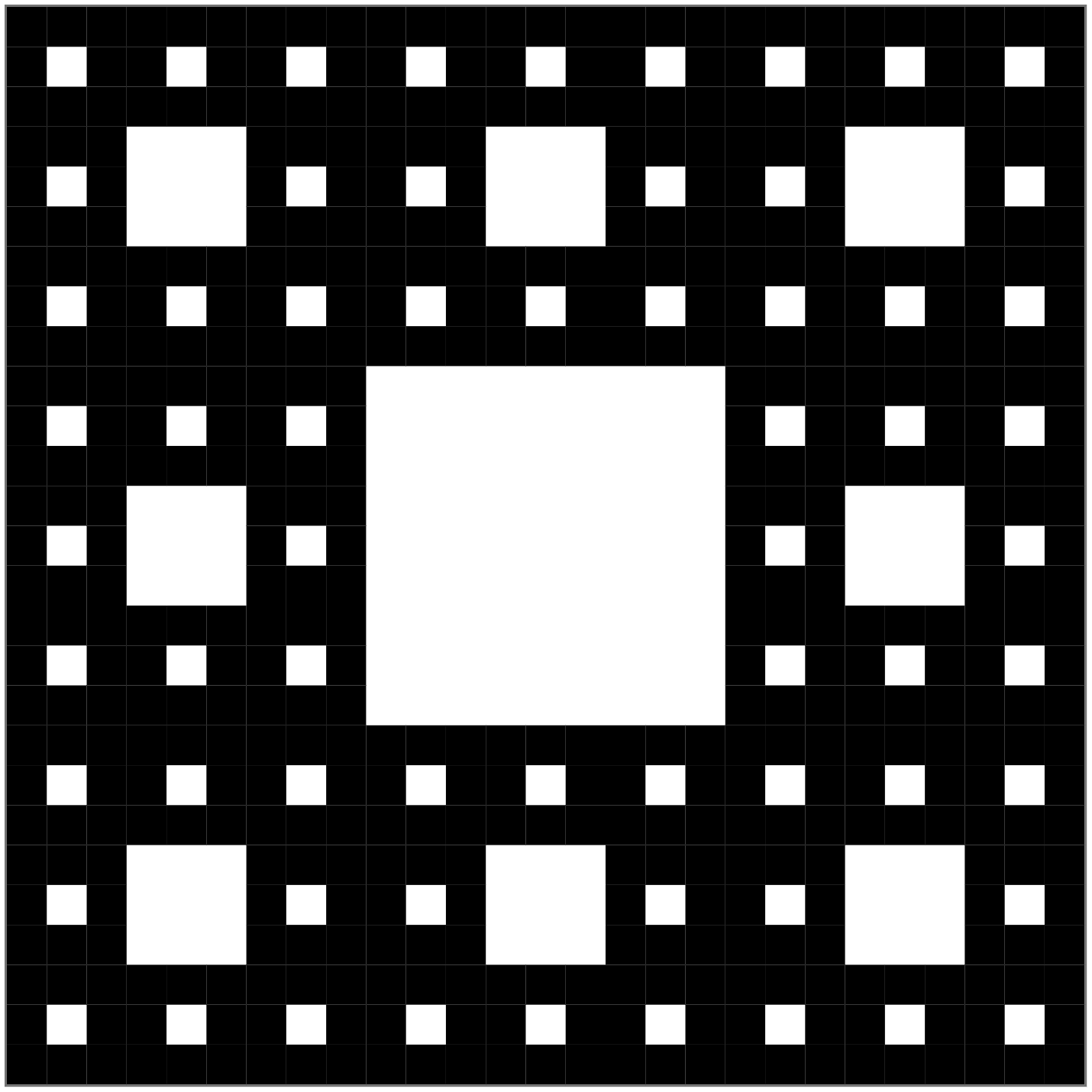} 
\end{tabular}
\end{center}
\vspace{-0.382cm}
\caption{The first three approximations   of $K\left(3,\mathcal{D}_S\right)$ in Example \ref{exmp:sierpinski}. }\label{fig:sierpinski}
\end{figure}
\end{example}

The topology of self-similar sets has been a focus of much attention in the study of fractal geometry. See for instance \cite{Hutchinson81,Falconer90,Akiyama-Thuswaldner04,Hata85}. 
In order to describe the topology of such a set, one may employ a useful map that is to be introduced below.  

Let $\mathcal{F}=\{\varphi_j: 1\le j\le q\}$ be an IFS. Given $w=i_1i_2\cdots i_n\in\{1,\ldots,q\}^n$, we denote by $\varphi_w$ the $n$-fold composition $\varphi_{i_1}\circ\cdots\circ \varphi_{i_n}$. For any infinite sequence $\omega=(i_n)_{n\ge1}\in\{1,\ldots,q\}^\infty$, we further denote by $\omega_n$ the prefix of   $\omega$ of length $n\ge1$. 
Then, by 
\cite[Theorem 3.1(3)]{Hutchinson81}  there is a unique point $x_\omega\in\mathbb{R}^d$, depending only on $\omega$, such that $\displaystyle \lim\limits_{n\rightarrow\infty}\varphi_{\omega_n}(y)=x_\omega$ holds for all $y\in\mathbb{R}^d$. 

\begin{definition}\label{def:symbolic_projection} Let $\mathcal{S}=\mathcal{S}_\mathcal{F}:\{1,\ldots,q\}^\infty\rightarrow K$  be the map  sending $\omega$ to $x_\omega$ and call it the symbolic projection of $\mathcal{F}$ (or of $K$).\end{definition}

\begin{remark}
If $K$ is a fractal square  or a higher dimensional counterpart, the symbolic projection $\mathcal{S}_\mathcal{F}$ provides sufficient information on how the topology of $K$ is related to that of $\{1,\ldots,q\}^\infty$, which becomes a Cantor set  when equipped with an appropriate metric.  Among others, one may refer to \cite{BandtKeller91} concerning such a relation. \end{remark}

In this survey,  
we restrict ourselves to fractal squares. We briefly review known results from recent studies and propose some open questions that arise naturally.  We also give figures that help demonstrate the underlying objects and `Figure X' illustrates `Example X'.
Section \ref{sec:pi_0} concerns the number $\pi_0(K)$ of components of  a fractal square. We will recall known results  and propose open questions concerning how to decide whether $\pi_0(K)$  is equal to one, greater than one but finite, or (uncountably) infinite. 
Section \ref{sec:cut} is about the case $\pi_0(K)=1$ and focuses on global and local  cut points. 
Section \ref{sec:components} is about the case  $\pi_0(K)>1$ and focuses on the shape of non-degenerate components. If non-degenerate components are line segments, we further discuss whether for any $C>0$ at most finitely many of them are of diameter $\ge C$, so that  they form a {\em null sequence}. Those discussions are related to the  classification of fractal squares in terms of the lambda function $\lambda_K$, that is introduced in \cite{FLY-2022} for any compact set $K\subset\mathbb{R}^2$. 
Section \ref{sec:N=3}  reviews all the known results from the literature concerning the classification of fractal squares $K$ with $N=3$. Except for a few unsolved cases, the classification is complete.
Section \ref{sec:pi_1} gives some open questions concerning fundamental groups. Similar questions about fractal cubes, which are three dimensional analogues of fractal squares, become more complicated.

\section{The Number of Components}\label{sec:pi_0}

Fractal squares $K=K(N,\mathcal{D})$ are self-similar sets contained in the unit square $[0,1]^2$. The underlying IFS
\[\mathcal{F}=\left\{f_d(x)=\frac{x+d}{N}: d\in\mathcal{D}\right\}
\]
always possesses  three basic properties.  First, the interior $U={\rm Int}\left([0,1]^2\right)$ is an open set such that $\bigcup_{d\in\mathcal{D}}f_d(U)$ is a disjoint union contained in $U$. Thus the well-known {\em open set condition} \cite[(5.2)]{Hutchinson81} is satisfied. Second, the (double) intersection $f_d(K)\cap f_{d'}(K)$ for any $d\ne d'$ in $\mathcal{D}$ is nonempty {\bf only if} $d-d'$  belongs to
\[
\left\{\pm\left[\begin{array}{c}0\\ 1\end{array}\right], \ \pm\left[\begin{array}{c}1\\ 0\end{array}\right], \
\pm\left[\begin{array}{c}1\\ 1\end{array}\right], \ \pm\left[\begin{array}{c}1\\ -1\end{array}\right]\right\}.\]
Moreover, if  $f_d(K)\cap f_{d'}(K)\ne\emptyset$ then it may consist of finitely many points,  countably many points,  a line segment,  or  a Cantor set together with countably many (or none) isolated points. Third, the {\bf triple} intersection $f_{d_1}(K)\cap f_{d_2}(K)\cap f_{d_3}(K)$ for disticnt $d_1,d_2,d_3\in\mathcal{D}$ is either empty or a singleton. In the latter case, it is easy to infer that for all $u\ne v\in\{d_1,d_2,d_3\}$ either $u-v$ or $v-u$  is a corner of $[0,1]^2$. Therefore,
\[
\left\{\pm(d_1-d_2),\pm(d_1-d_3),\pm(d_2-d_3)\right\}\ \bigcap\
\left\{\left[\begin{array}{c}0\\ 0\end{array}\right], \left[\begin{array}{c}1\\ 0\end{array}\right], \ \left[\begin{array}{c}1\\ 1\end{array}\right], \ \left[\begin{array}{c}0\\ 1\end{array}\right]\right\}
\]
has exactly three elements, each of which  belongs to $\mathcal{D}$. On the one hand,  $K=K(N,\mathcal{D})$ and its components demonstrate very rich properties from a topological point of view. On the other, the connections between those properties and the digit set $\mathcal{D}$ are direct and illustrate a variety of subtleties. A large part of those connections have been well understood, due to recent as well as classical results from the literature.

Let us start from how the digit set $\mathcal{D}$ is related to the number of components in $K=K(N,\mathcal{D})$, to be denoted by $\pi_0(K)$. To do that, we need
a graph that is defined for all IFS's.
\begin{definition}\label{def:Hata_Graph}
Let $\mathcal{F}=\{f_1,\ldots,f_q\}$ be an IFS and $K$ the attractor. The first Hata graph $\mathcal{G}_1(\mathcal{F})$ is the one with vertex set $\mathcal{V}_1=\{1,\ldots,q\}$ in which two vertices $i\ne j$ are incident if and only if $f_i(K)\cap f_j(K)\ne\emptyset$. Similarly, we can define the $k$-th Hata graph  $\mathcal{G}_k(\mathcal{F})$ for $k\ge2$, whose  vertex set $\mathcal{V}_k=\{1,\ldots,q\}^k$. 
\end{definition}

\begin{theorem}[{\bf \cite[Theorem 4.6]{Hata85}}]\label{theo:Hata85}
$\mathcal{G}_1(\mathcal{F})$ is connected if and only if $K$ is. In such a case, $K$ is a locally connected continuum.
\end{theorem}

If $K=K(N,\mathcal{D})$ is a fractal square, we also write $\mathcal{G}_1(\mathcal{D})$ instead of $\mathcal{G}_1(\mathcal{F})$.  Moreover,  one can induce a simple characterization of $\mathcal{G}_1(\mathcal{D})$.

\begin{theorem}[{\cite[Theorem 2.2]{Ruan-Wang17}}]\label{theo:compute_Hata}
$\mathcal{G}_1(\mathcal{D})$ has an edge between $d'\ne d''$ if and only if $K^{(2)}+d'$ intersects $K^{(2)}+d''$. Consequently, $K$ is connected if and only if $K^{(3)}$ is.
\end{theorem}

In order to determine whether $K$ has finitely many components,  one may use the approximations $K^{(n)}(n\ge1)$ instead of $\mathcal{G}_n(\mathcal{D})$. We can rephrase the second part of Theorem \ref{theo:compute_Hata}  as follows.
\begin{theorem}\label{theo:connected_K3}
$\pi_0(K)=1$ if and only if $\pi_0\left(K^{(3)}\right)=1$.
\end{theorem}

Observe that $\pi_0\left(K^{(n+1)}\right)\ge\pi_0\left(K^{(n)}\right)$ holds for all $n\ge1$. In Example \ref{exmp:Hata_Graph}, we give a simple fractal square $K$ satisfying $\pi_0\left(K^{(3)}\right)>\pi_0\left(K^{(2)}\right)=\pi_0\left(K^{(1)}\right)=1$.

\begin{example}\label{exmp:Hata_Graph}
Let $\mathcal{D}_1=\{(0,0),(1,0),(1,1),(1,2),(2,1)\}\subset\{0,1,2\}^2$. The  fractal square $K(3,\mathcal{D}_1)$ was  given in \cite[Example 2.3]{DaiLuoRuanWang_23}. See Figure \ref{fig:Hata_1} for a depiction of the approximations $K^{(j)}\ (j=1,2,3,6)$ and check that $\pi_0\left(K^{(1)}\right)=\pi_0\left(K^{(2)}\right)=1$ while $\pi_0\left(K^{(3)}\right)=3=\pi_0\left(\mathcal{G}_1(\mathcal{D}_1)\right)$. 
\begin{figure}[ht]   
\begin{center}
\begin{tabular}{cccc}
\includegraphics[width=2.75cm]{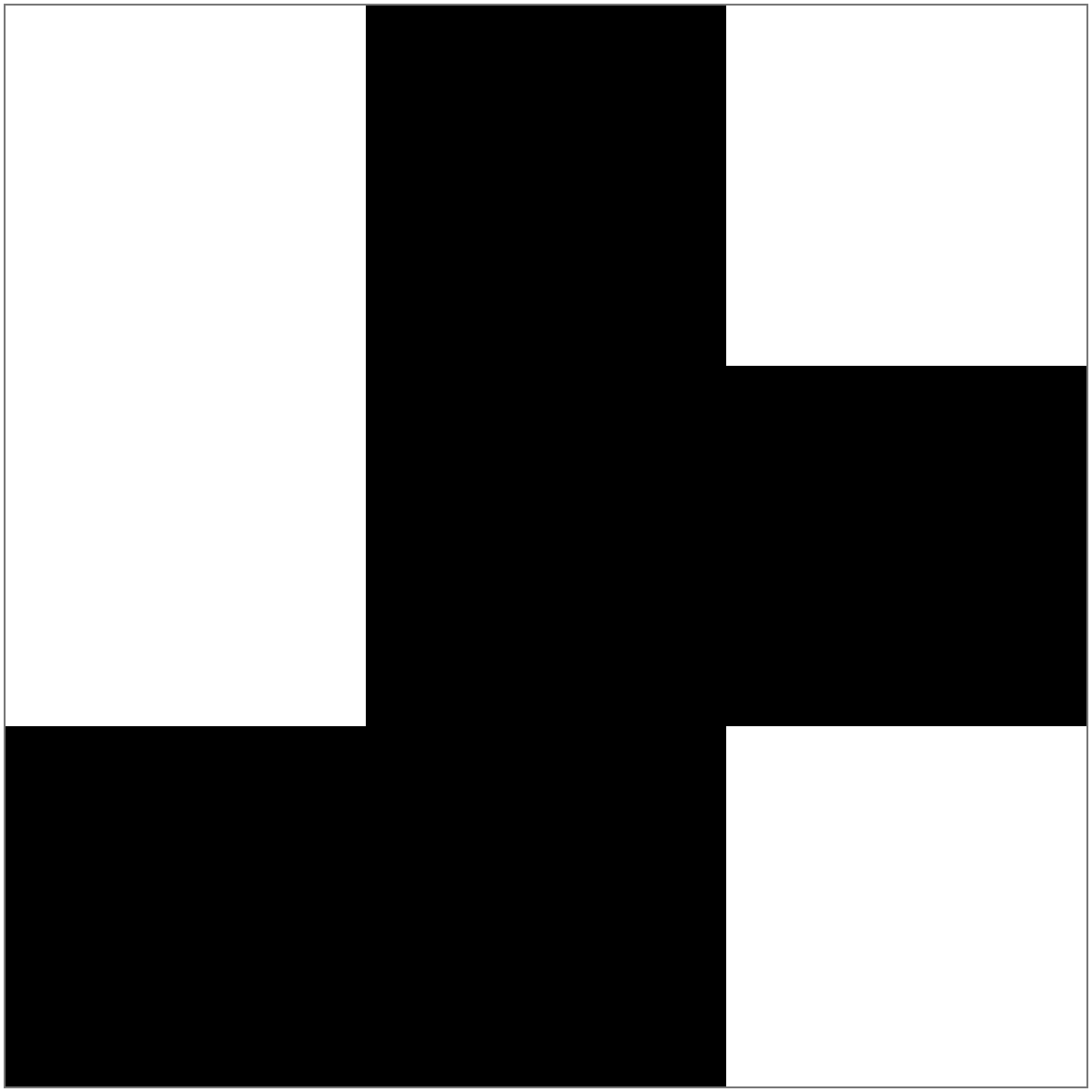}
& \includegraphics[width=2.75cm]{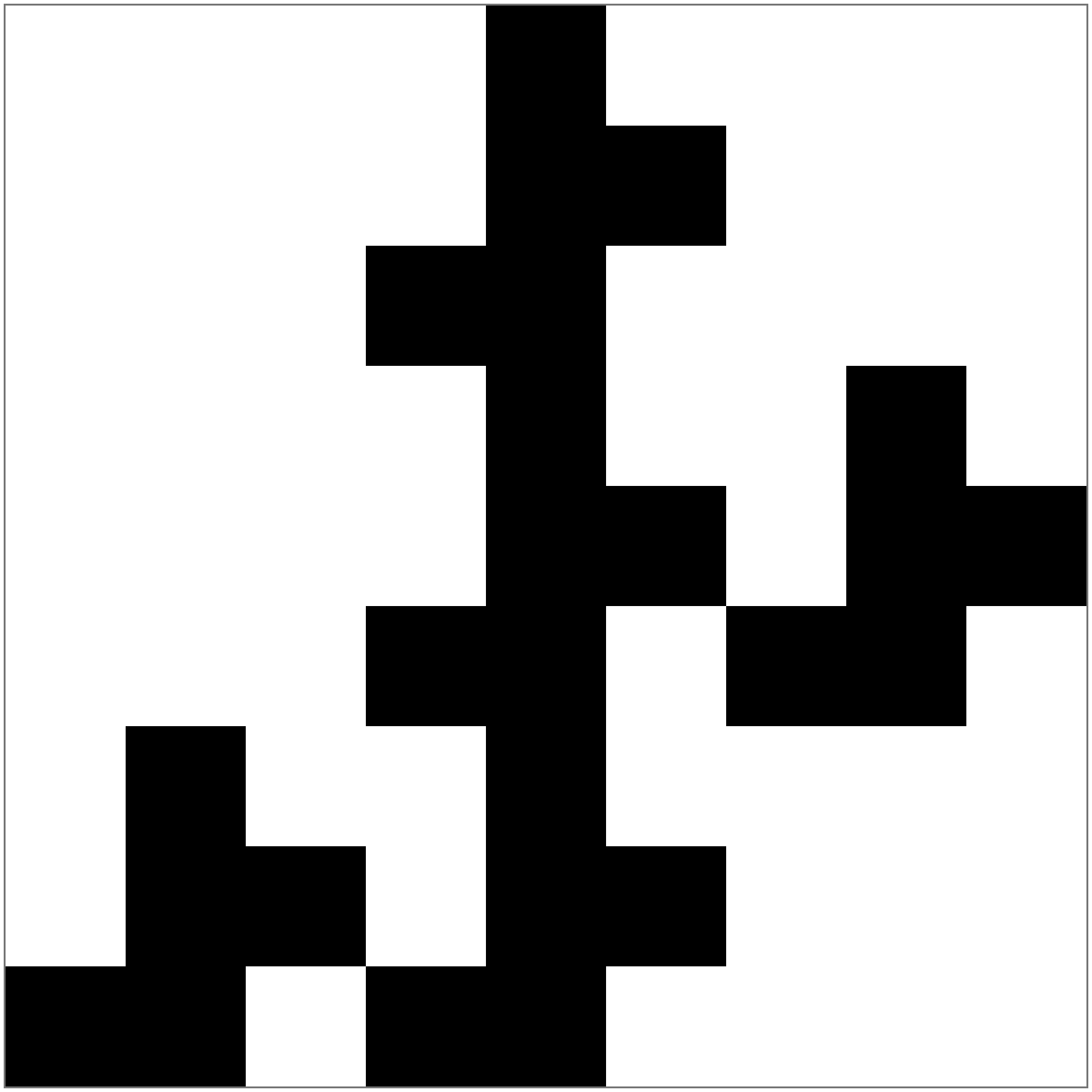}
&\includegraphics[width=2.75cm]{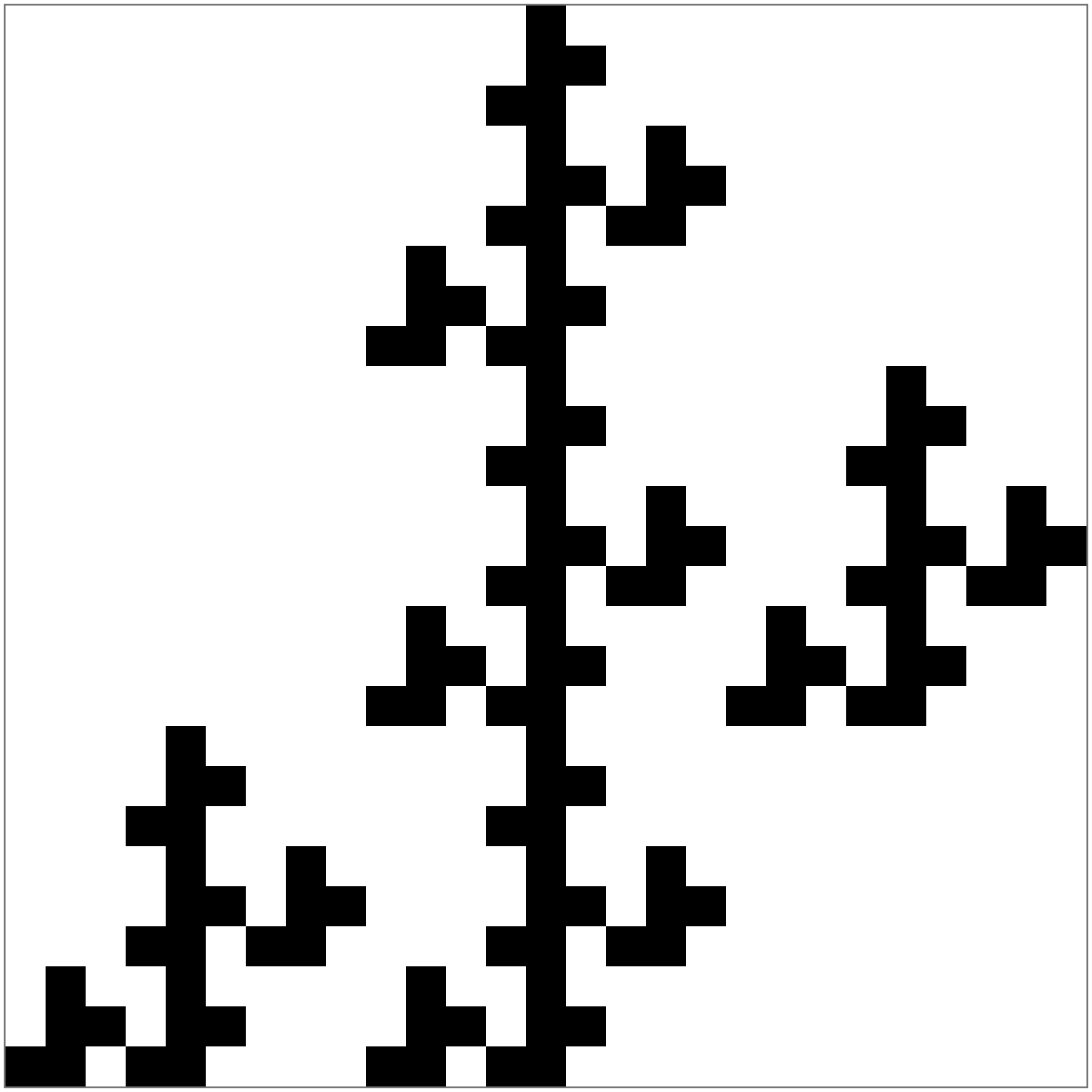}
& \includegraphics[width=2.75cm]{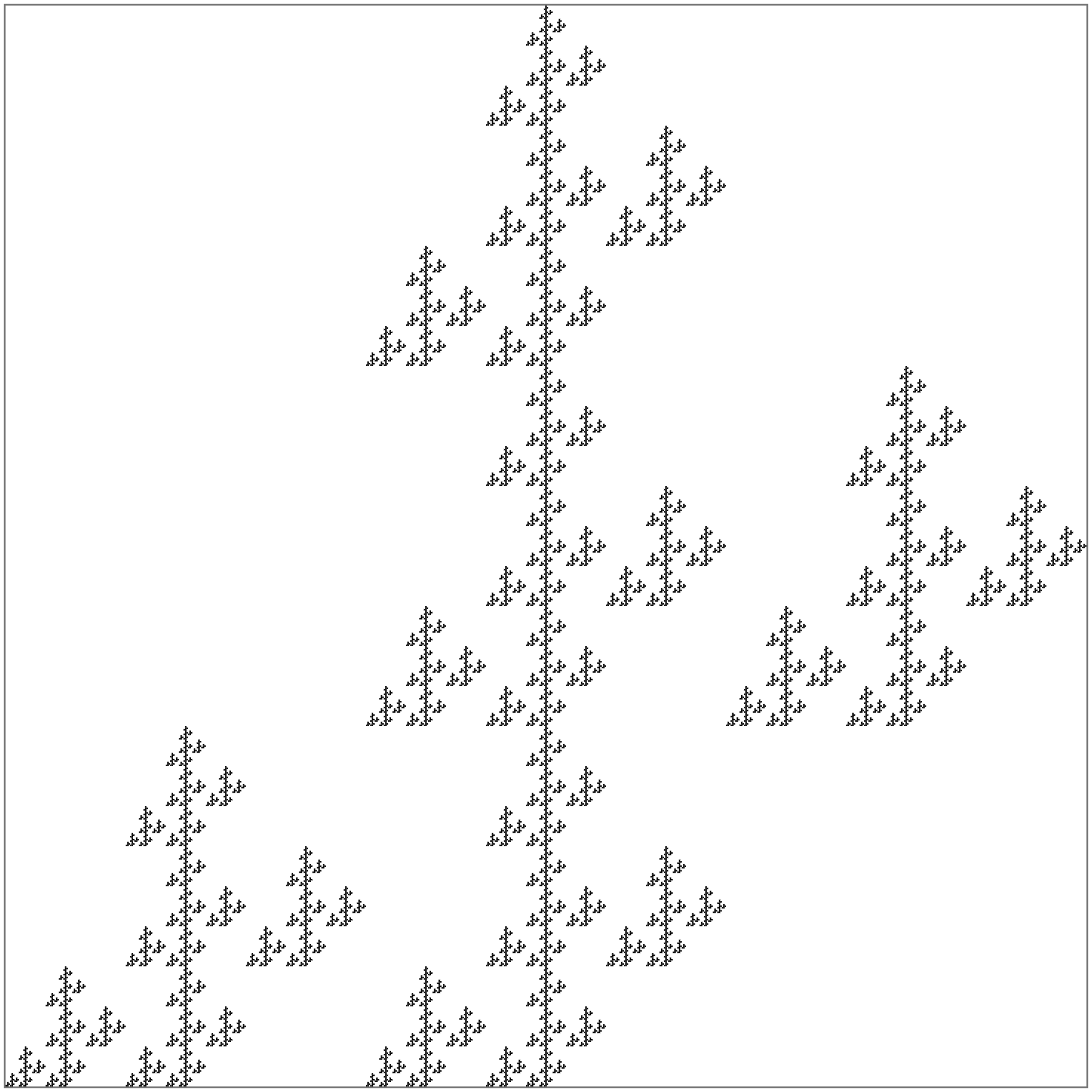}
\end{tabular}
\end{center}
\caption{The approximations $K^{(j)}(j=1,2,3,6)$ of $K(3,\mathcal{D}_1)$ in Example \ref{exmp:Hata_Graph}.}\label{fig:Hata_1}
\end{figure}
\end{example}

Set $q=\pi_0\left(K^{(1)}\right)$ and denote by $\mathcal{P}_1,\ldots,\mathcal{P}_q$ the components of  $K^{(1)}$.  Let $\mathcal{V}^{(1)}_i(1\le i\le q)$ consist of  $d\in \mathcal{D}$ satisfying $\frac{[0,1]^2+d}{N}\subset\mathcal{P}_i$ and  $K_i^{(1)}$ the union of the translates $\frac{K+d}{N}$ with $d\in\mathcal{V}^{(1)}_i$. 

For the sake of convenience, we let $\alpha_i=\{i\}\times[0,1]$ and $\beta_i=[0,1]\times\{i\}$, where $i=0,1$. In the sequel, we give in several lemmas some   observations that are either immediate or may be obtained by using self-similarity and standard arguments from plane topology.

\begin{lemma}\label{lem:GIFS}
The homothetic image $NK_i^{(1)}=\left\{Nx: x\in K_i^{(1)}\right\}(1\le i\le q)$  is a finite union of certain translates of the compact sets $K_j^{(1)}$. Thus $\left\{K_i^{(1)}:1\le i\le q\right\}$ becomes a system of attractors for a graph-directed IFS in the sense of \cite{MauldinWilliams88}.
\end{lemma}

\begin{lemma}\label{lem:sepration_arc}
If no component of $K^{(1)}$ intersects each of the four sides $\alpha_0,\alpha_1,\beta_0,\beta_1$, then one can find either a Jordan arc $\gamma$ in $[0,1]^2\setminus K^{(1)}$ that intersects both $\{0\}\times(0,1)$ and $\{1\}\times(0,1)$ or 
one that intersects both $(0,1)\times\{0\}$ and $(0,1)\times\{1\}$. \end{lemma}

\begin{lemma}\label{lem:points_exist}
If $K^{(1)}$ has a component that is disjoint from one of the four broken lines $\alpha_i\cup\beta_j$ (with $i,j=0,1$) then $K$ has a point component and $\pi_0(K^{(1)})<\pi_0(K^{(2)})$. 
\end{lemma}

\begin{lemma}\label{lem:NS_EW}
If $\pi_0(K^{(1)})=\pi_0(K^{(2)})>1$ then there are two possibilities: (1) every component of $K^{(1)}$ intersects both $(0,1)\times\{0\}$ and $(0,1)\times\{1\}$, (2)  every component of $K^{(1)}$ intersects both $\{0\}\times(0,1)$ and $\{1\}\times(0,1)$. In such cases,  we respectively call 
$\mathcal{D}$ a digit set with a {\bf north-south pattern} and one with an {\bf east-west pattern}.
\end{lemma}

The next result provides further insight concerning  point components.
\begin{theorem}[{\cite[Theorem 1.1]{Huang-Rao21}}]\label{theo:dimension_drop}
If $K$ has a point component then all its non-degenerate components form a subset $K_c$ whose Hausdorff dimension $\dim_HK_c$ is strictly less than $\dim_HK$.
\end{theorem}

The next result  extends  Theorem \ref{theo:connected_K3}. 
\begin{theorem}[{\cite[Theorem 2.6]{CKLY-2025}}]\label{theo:finite_pi0}
If $\pi_0\left(K^{(n)}\right)=\pi_0\left(K^{(n+1)}\right)$ for $n\ge2$ then $\pi_0(K)=\pi_0\left(K^{(n)}\right)$. In particular, if $\pi_0\left(K^{(3)}\right)=1$ (hence  $\pi_0\left(K^{(2)}\right)=1$) then $\pi_0(K)=1$.
\end{theorem}

The assumption $n\ge2$ is necessary. In Example \ref{exmp:Hata_Graph}, we have given a simple fractal square, to show that  $\pi_0\left(K^{(1)}\right)=\pi_0\left(K^{(2)}\right)=1$ while $K$ has uncountably many components.
The  following dichotomy provides a fundamental starting point in exploring the topology of a fractal square.
\begin{theorem}[{\bf \cite{Xiao21}}]\label{theo:Xiao21}
$K$  has either finitely many or uncountably many components.
\end{theorem}

Actually, by combining  Lemmas \ref{lem:points_exist} and \ref{lem:NS_EW} with Theorem \ref{theo:finite_pi0} and  the first part of \cite[Theorem 3.3]{CKLY-2025}, one may  strengthen Theorem \ref{theo:Xiao21} in the following way.
\begin{theorem}\label{theo:Xiao21_strengthened}
Every fractal square $K$ falls into the following categories: (1) $\pi_0\left(K^{(n)}\right)=\pi_0\left(K^{(n+1)}\right)=\pi_0(K)$ for some $n\ge2$, 
(2)  $\pi_0\left(K^{(n)}\right)<\pi_0\left(K^{(n+1)}\right)$ for all $n\ge2$. In the latter case there are just two subcases, either $K$ has  uncountably many point components or $K$ is the product of a Cantor set and $[0,1]$.
\end{theorem}


\begin{example}\label{exmp:two_comp}
Let $K=K(5,\mathcal{D})$ be given  in \cite[Figure 1]{Xiao21}, where the digit set $\mathcal{D}$ consists of three parts: 
$\{(u,v): u=0,1; 0\le v\le4\}$, $\{(u,v): u=1,2; v=3,4\}$ and $\{(3,v): v=0,1\}$.  Figure \ref{fig:Xiao21} illustrates the approximations  $K^{(n)}$ for $n=1,2,4$.  Notice that $\pi_0\left(K^{(1)}\right)=\pi_0\left(K^{(2)}\right)=2$.   
By induction, we can infer that $\pi_0\left(K^{(n)}\right)=2$ for all $n\ge3$. It then follows that  $\pi_0(K)=2$. This type of argument is used in the proof for Theorem \ref{theo:finite_pi0}.
\begin{figure}[ht]    
\begin{center}
\begin{tabular}{ccc}
\includegraphics[width=3.0cm]{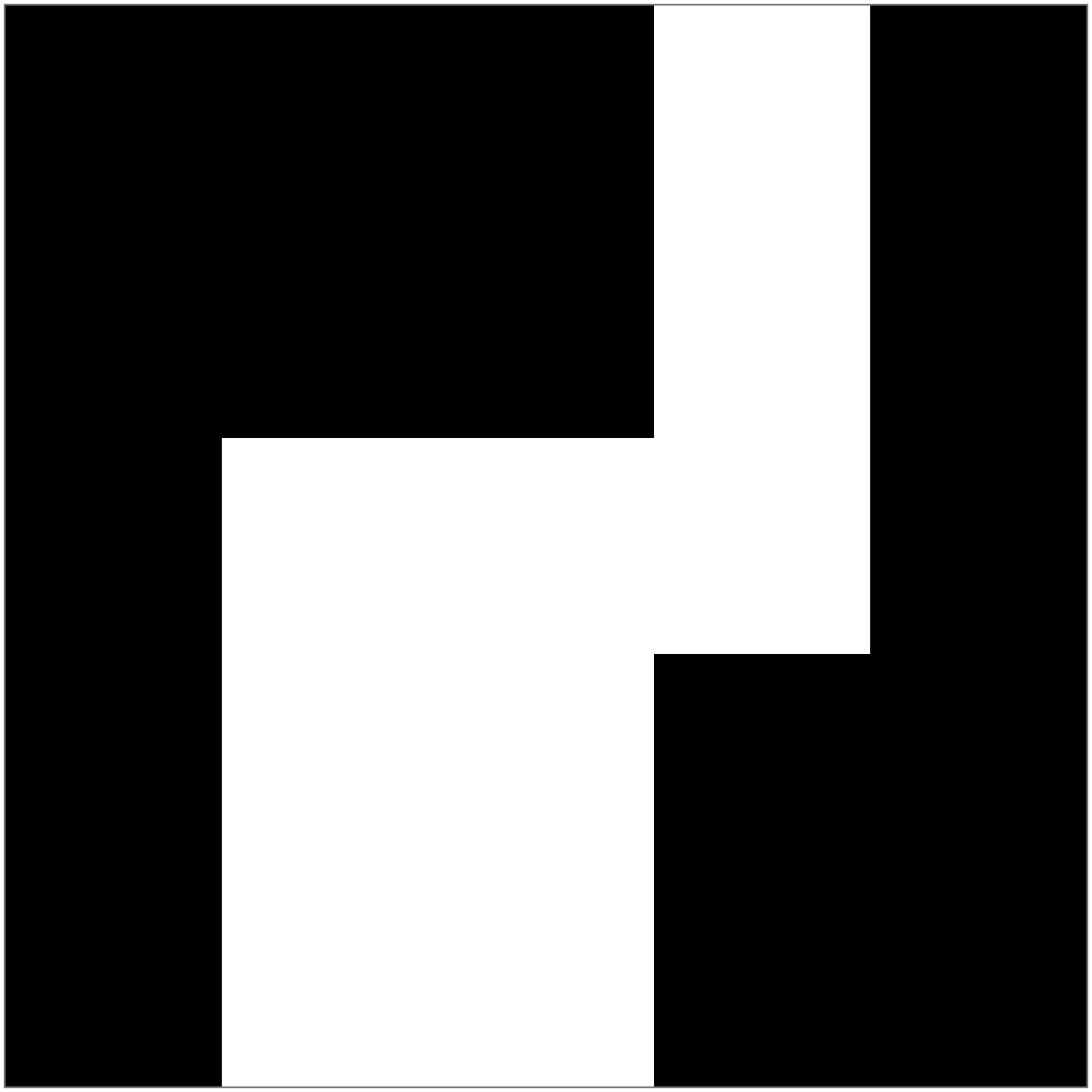}&
\includegraphics[width=3.0cm]{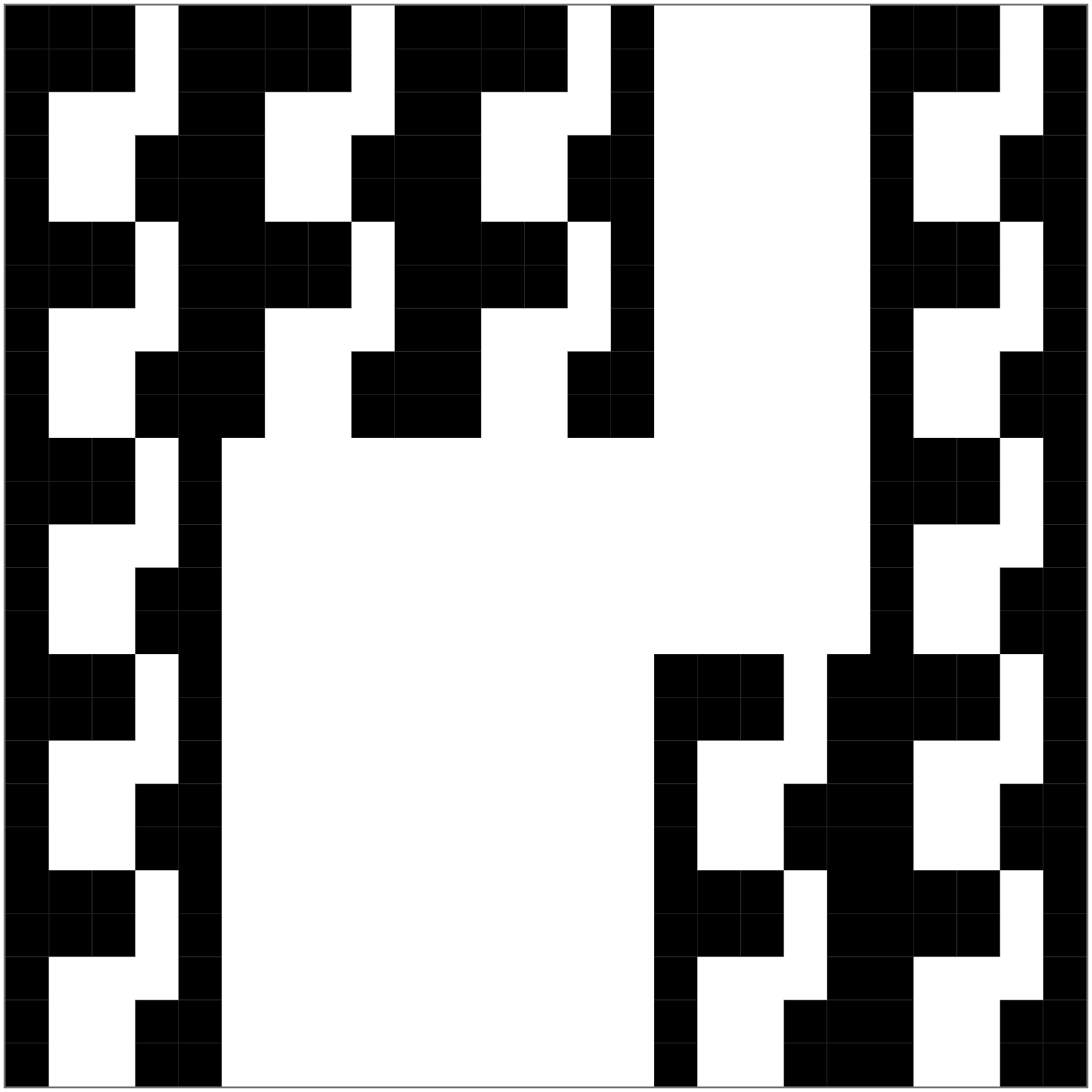}&
\includegraphics[width=3.0cm]{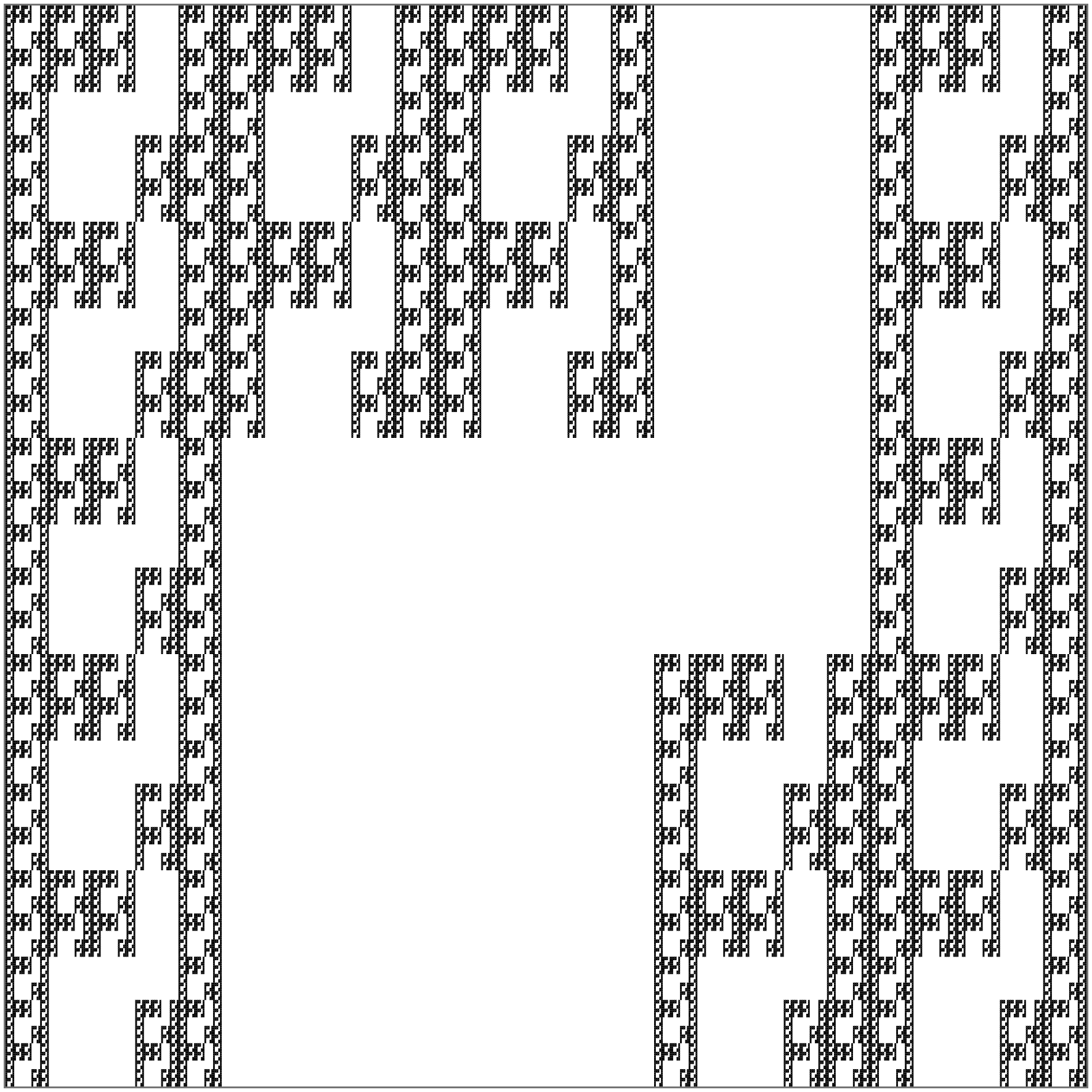}
\end{tabular}
\end{center}
\caption{The approximations $K^{(1)},K^{(2)}, K^{(4)}$ of $K=K(5,\mathcal{D})$ in Example \ref{exmp:two_comp}.}\label{fig:Xiao21}
\end{figure}
\end{example}

To conclude this section, we propose the following open questions.

\begin{question}\label{algo:components}
Given a fractal square $K=K(N,\mathcal{D})$, is it possible to determine $\pi_0(K)$ by checking the initial approximations $K^{(j)}$, say for $j\le4$?
\end{question}   

\begin{question}\label{algo:total_disconn}
Given a fractal square $K=K(N,\mathcal{D})$,  can we decide whether  $K$ is totally disconnected by checking the initial approximations $K^{(j)}$, say for $j\le4$?
\end{question}

\section{Cut Points and Local Cut Points}\label{sec:cut}

Let us start from the following terminology.
\begin{definition}[{\bf \cite[p.41]{Whyburn42}}]\label{def:cut}
A  {\em (global) cut point} of $X$ means a point $x\in X$ whose complement in $X$ is disconnected. A {\em cutting} of $X$ is a subset $C\subset X$ such that $X\setminus C$ is disconnected. And a point $x$ in a continuum $X$ is said to be a {\em local cut point} provided that it has a connected neighborhood $V_x$ with $V_x\setminus\{x\}$ a disconnected set. \end{definition}

\begin{example}\label{exmp:cut}
The Sierpinski's carpet  given in Example \ref{exmp:sierpinski} has no local cut point. The fractal square given in Example \ref{exmp:two_comp} has two components. Let $Q$ be the one to the right.  See Figure \ref{fig:cut}. Then $x_0=\left(\frac{18}{25},\frac15\right)$ is a global cut point of $Q$ and $x_j=\left(\frac{23}{25},\frac{j}{5}\right)$ with $1\le j\le 4$ are local cut points.
\begin{figure}[ht]    
\begin{center}
\includegraphics[width=8.0cm]{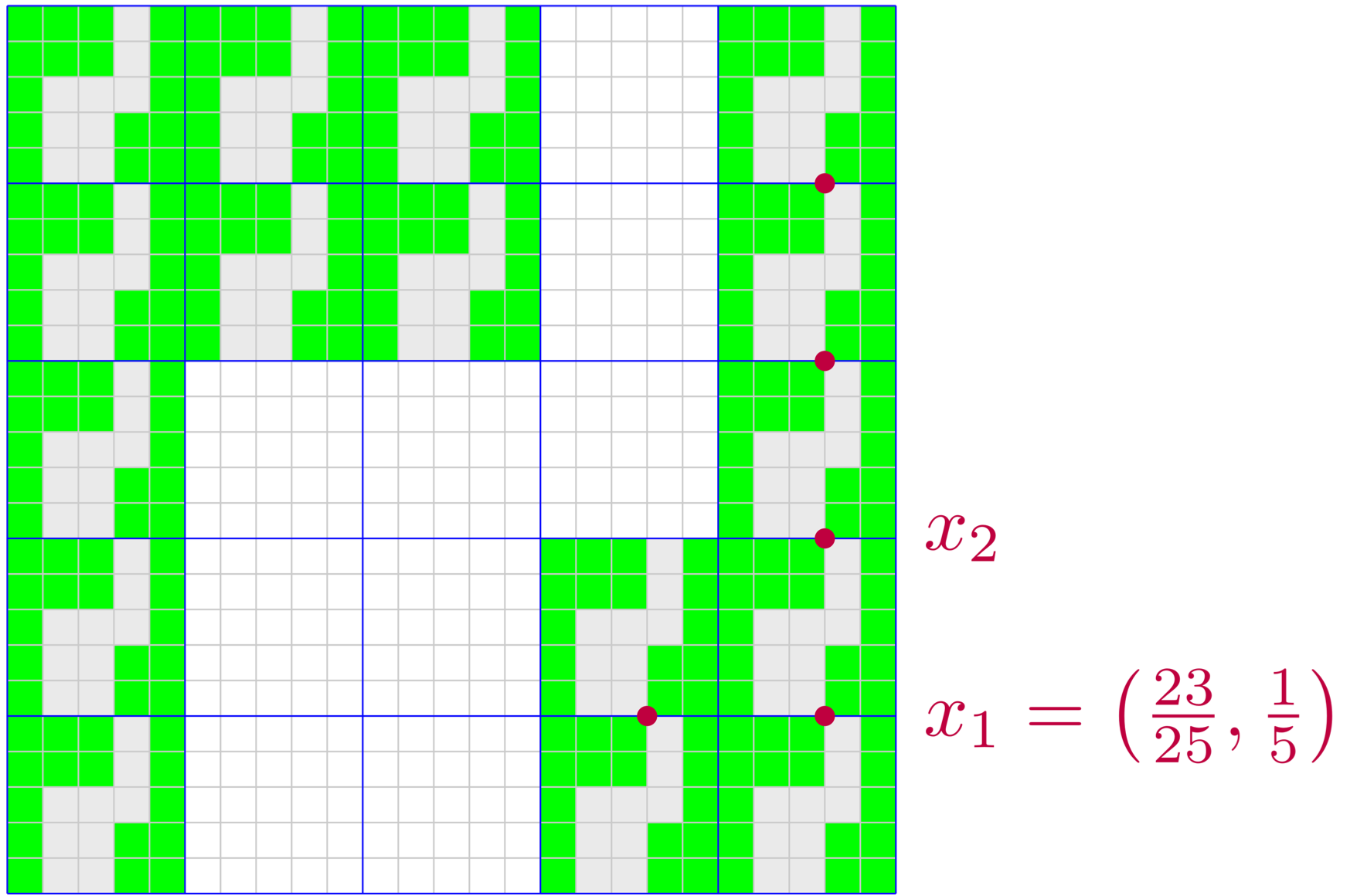}
\end{center}
\caption{A magnified view of the approximation $K^{(2)}$ in Figure \ref{fig:Xiao21} and the points $x_j\ (0\le j\le4)$.}\label{fig:cut}
\end{figure}

\end{example}

With Theorems \ref{theo:Hata85} and \ref{theo:Xiao21}, we know that a fractal square $K(N,\mathcal{D})$ has either  one component (thus is connected), or more than one but finitely many components, or  uncountably many components. In the current section we restrict to those that have finitely many components. 


General self-similar sets with finitely many components are characterized  as follows.
\begin{theorem}[{\cite[Theorem 1]{LRX-2022}}]\label{theo:LC_general}
A self-similar set is locally connected if and only if it has finitely many components.
\end{theorem}

Recall that a Peano continuum is just the image of $[0,1]$ under a continuous map. By the well known Hahn-Mazurkiewicz Theorem, a continuum is locally connected if and only if it is a Peano continuum. Also recall that a compactum means a compact metric space and the notion of Peano continuum is naturally generalized to compacta that may not be connected.

\begin{definition}\label{def:PC}
A {\em Peano compactum}  is a compact metric space with locally connected components such that for any constant $C>0$ at most finitely many components may have a diameter greater than $C$. 
\end{definition}

\begin{remark}
We refer to \cite[Theorems 1-3]{LoridantLuoYang19} concerning how the notion of Peano compactum is related to the study of polynomial Julia sets. In the next section, we will consider a special upper semi-continuous decomposition,  having a Peano compactum as the quotient space, that is well-defined for any compact set $K$ in the plane. With the help of this,  we can define the so-called lambda function $\lambda_K$ and apply it to the study  of fractal squares. 
\end{remark}

Notice that a connected Peano compactum is a Peano continuum. Moreover, by Theorem \ref{theo:Xiao21} a fractal square either is the union of finitely many  Peano continua or has uncountably many components. In the former case, one may wonder whether a component of $K$ is an {\em S-curve} in the sense of Whyburn \cite{Whyburn58}. See Definition \ref{def:S_curve} below. In particular, when $K$ is connected, we want to seek conditions under which $K$ itself is an S-curve. 
\begin{definition}\label{def:S_curve}
An S-curve means a locally connected  continuum  $S$ in the plane  with no interior point    such that the boundary of each complementary domain of $S$ is a simple closed curve and no two of these complementary domain boundaries intersect.
\end{definition}

The following are well known.
\begin{theorem}[{\cite[Theorem 3]{Whyburn58}}]\label{theo:S_curve}
All S-curves are homeomorphic. In particular, every S-curve is homeomorphic with the Sierpinski carpet $K(3,\mathcal{D}_S)$.
\end{theorem}

\begin{theorem}[{\bf \cite[Theorem 4]{Whyburn58}}]\label{theo:Whyburn_criterion}
A Peano continuum in the plane with no interior point is an S-curve if and only if it has no local cut point.
\end{theorem}

Let $\#A$ denote the cardinality of a set. The characterization below is useful.

\begin{theorem}[{\cite[Theorem 1.5]{DaiLuoRuanWang_23}}]\label{theo:global_cut}
A fractal square $K=K(N,\mathcal{D})$ has a cut point if and only if every $\mathcal{G}_k(\mathcal{D})$ with $k\ge2$ contains a vertex $v_k$ such that at least two components of $\mathcal{G}_K(\mathcal{D})\setminus\{v_k\}$ have more than $\#\mathcal{D}^{k-1}$ vertices.
\end{theorem}

Let $\mathcal{S}:\{0,1,\ldots,N-1\}^\infty\rightarrow K$
 be given as in Definition \ref{def:symbolic_projection}. Then for each $x_0\in K$ we shall have 
 $\#\mathcal{S}^{-1}(x_0)\in\{1,2,3,4\}$.  
In what follows, we assume that $K$ is connected and $x_0\in K$ is a local cut point. Sometimes, we also refer to $\#\mathcal{S}^{-1}(x_0)$ as the {\em multiplicity} of $x_0$.

If $\#\mathcal{S}^{-1}(x_0)=1$,  there is a unique $d\in\mathcal{D}$ such that $x_0\in f_d(K)$. Thus $f_d(K)$ is a neighborhood of $x_0$ with respect to the relative topology of $K$. This implies that $x_0$ is a cut point of $f_d(K)$ and $f_d^{-1}(x_0)$ is a cut point of $K$, since $f_d$ restricted to $f_d(K)$ is a homeomorphism onto $K$. 

If $\#\mathcal{S}^{-1}(x_0)\ge2$, one can find $n\ge1$ and distinct digits $d_1,\ldots,d_j\in\mathcal{D}_n$  (with $2\le j\le 4$) such that $\left(\bigcup_i\frac{K+d_i}{N^n}\right)\setminus\{x_0\}$ is disconnected and the following two requirements are both satisfied: 
 \begin{eqnarray}
\{x_0\}&=&\bigcap_{i=1}^j\frac{K+d_i}{N}
\\
x_0&\notin& \bigcup_{d\in\mathcal{D}_n\setminus\{d_1,\ldots,d_j\}}\frac{K+d}{N}
\end{eqnarray}
For related discussions, we refer to \cite[Theorem 1.6]{DaiLuoRuanWang_23}. 

\begin{remark}
Let $K(5,\mathcal{D})$ be given as in Example \ref{exmp:two_comp} and $Q$ its component to the right.  The global cut point $\left(\frac{18}{25},\frac15\right)$ and the four local cut points $\left(\frac{23}{25},\frac{j}{5}\right)$ of  $Q$, as given in Example \ref{exmp:cut}, are each a point with multiplicity two. Note that a connected fractal square has no local cut point with multiplicity four. Up to now, the authors have not found such a result from the literature, although the underlying proof is not so difficult. On the other hand, it remains unclear whether a connected fractal square could have a global or local cut point  whose multiplicity is three. A thorough analysis of the existence of global or local cut points with multiplicity $\ge2$ may be of some help, when one wants to resolve Question \ref{algo:sierpinski} below. \end{remark}

An issue of interest is to characterize fractal squares that are S-curves, and those that have finitely many components each of which is an S-curve. 
Moreover, we wonder whether we can do that by checking the approximations $K^{(n)}$ for small $n$. Here we prefer the approximations than the Hata graphs $\mathcal{G}_j(\mathcal{D})$, since one of the ultimate issues is to compose appropriate programs that may settle the following. 

\begin{question}\label{algo:sierpinski}
How  to determine whether a connected fractal square $K$ is an S-curve by just checking the initial approximations $K^{(j)}$, say for $j\le4$?
\end{question}

\section{The Structure of the Components}\label{sec:components}

The previous two sections  deal with the number $\pi_0(K)$ of components in a fractal square $K=K(N,\mathcal{D})$ and the existence of global or local cut points of $K$ (if $\pi_0(K)=1$) or of its components (if $\pi_0(K)\in(1,\infty)\cap\mathbb{Z}$), respectively. In this section we investigate the geometry of the components of $K$, such as  the local connectedness of non-degenerate components  and the existence of infinitely many components whose diameters are bounded from below by a constant $C>0$.

The components $P$  of $K$ may be divided into three types: (1) a point, (2) a line segment, (3) a continuum that is not a line segment. To carefully analyze the components of $K$, we need to employ the following sets and their complements: (1)  the union $H=K+\mathbb{Z}^2$; (2)  the sequence $H_j=K^{(j)}+\mathbb{Z}^2$ for  $j\ge1$.

Let us recall some fundamental properties concerning $H$ and $H_j$.
\begin{lemma}[{\cite[Equations (2.1)]{Lau-Luo-Rao13}}]\label{lem:containments}
The inclusion  $\displaystyle H_j\supset NH_{j+1}=\left\{Nx: x\in H_{j+1}\right\}$ holds for all $j\ge1$ hence 
$H\supset NH$. Consequently, if $H$ contains a line segment of irrational slope then $H=\mathbb{R}^2$, or equivalently, $K=[0,1]^2$.\end{lemma}

\begin{lemma}[{\cite[Lemma 2.1 and Theorem 2.2]{Lau-Luo-Rao13}}]\label{lem:unbounded_components}
Either all components of $H^c$ are unbounded or none of them is. In the latter case, all components of $H^c$ are of diameter $\displaystyle\le\frac{\sqrt{2}(N^2+1)^2}{N}$.\end{lemma}

\begin{lemma}[{\cite[Lemma 2.3, Corollary 2.4, and Theorem 2.5]{Lau-Luo-Rao13}}]\label{lem:line_in_H}
If $K$ contains a line segment of slope $k_0$ then $H$ contains an infinite line of slope $k_0$. If $K$ contains two line segments of distinct slope then $K$ has a component which is not a line segment. Moreover, $K$ has a component which is not a line segment if and only if every component of $H^c$ is bounded.
\end{lemma}

With  the above results,  Lau et al \cite{Lau-Luo-Rao13} culminate in the following dichotomy.

\begin{theorem}[{\bf \cite[Theorems 1.1 and  2.2]{Lau-Luo-Rao13}}]\label{theo:LauLuoRao}
Given a fractal square $K=K(N,\mathcal{D})$, there are two possible cases: (1)  all components of $H^c$ are unbounded and every component of $K$ is either a point or a line segment;  (2)  all components of $H^c$ are of diameter 
$\displaystyle\le\frac{\sqrt{2}(N^2+1)^2}{N}$. In the former case, all the line segments in $K$, if there are any, are of the same {\bf rational} slope.\end{theorem}

In case one, there are three subcases: 
\begin{enumerate} 
\item[(1.a)]  all components are points hence $K$ is totally disconnected; 
\item[(1.b)]  no component is a point; 
\item[(1.c)]  there are both points and  line segments, among the components. 
\end{enumerate}
In case two, there are two subcases: 
\begin{enumerate} 
\item[(2.a)]  there are finitely many components, each of which is a Peano continuum; 
\item[(2.b)]  there are uncountably many point components and countably many non-degenerate ones. 
\end{enumerate}

The totally disconnected fractal squares, for subcase (1.a), may be characterized by the existence of special sub-continua  of $[0,1]^2$ that are disjoint from the interior of $ K^{(n)}$, for some $n\ge1$.
\begin{definition}[{\cite[p.36, Definition 27]{Roinestad10}}]\label{def:complete_path}
Given a fractal square $K=K(N,\mathcal{D})$ of order $N\ge2$, a {\em complete path $P$} at level $n\ge1$ means the union of squares $\frac{[0,1]^2+d}{N^n}$ for $d$ belonging to a subset $\mathcal{D}_P$ of $\left\{0,1,\ldots,N^n-1\right\}^2\setminus\mathcal{D}_n$ such that  the next two requirements are both satisfied:
\begin{itemize} 
\item[(a)] There exist $0\le i,j\le N-1$ with $\mathcal{D}_P\supset\{(i,0),(i,N^n-1),(0,j),(N^n-1,j)\}$
\item[(b)] The interior of $P$ is connected and $P$ itself is  disjoint from the interior of $K^{(n)}$.
\end{itemize}
\end{definition}
The above notion is more general than the one given in \cite{Roinestad10}.  Using Theorem \ref{theo:LauLuoRao} and standard arguments, one may easily infer the following.
\begin{theorem}[{see \cite{Roinestad10} or \cite[Proposition 2.3]{Luo-Liu16}}]\label{theo:total_disconnected}
A fractal square $K(N,\mathcal{D})$ is totally disconnected if and only if there is a complete path, at level $n$ for some $n\ge1$.
\end{theorem}

\begin{remark}\label{rem:bound_for_level}
An issue of interest is to estimate the level $n$ from above. That is to say, one may wonder whether there is a universal bound, say $n_0$, such that if there is a complete path then there is one at level $n\le n_0$. For fractal squares of order $N=3$, the five $K(3,\mathcal{D})$ given in \cite[Fig. A.1]{Luo-Liu16} are totally disconnected. Clearly, for none of them one can find a complete path at level one.
\end{remark}

Due to recent studies from the literature, we can further clarify all subcases in the following way.
By \cite[Theorem 1.1]{Huang-Rao21}, see also Theorem \ref{theo:dimension_drop}, 
if $K$ has a point component then the non-degenerate components form a subset whose Hausdorff dimension is strictly less than that of $K$.
This points out a significant property shared by fractal squares that have point components. 

We may characterize subcase (1.b) as follows.

\begin{theorem}[{\cite[Theorem 3]{CKLY-2025}}]\label{theo:no_point}
If $K=K(N,\mathcal{D})$ has uncountably many components and none of them is a point, then $K$ is the product of a Cantor set with $[0,1]$. Moreover, this happens if and only if  $\#\mathcal{D}=Nq$ for some $q\ge2$ and  $K^{(1)}$ consists of either $q$ rectangles of the form $\left[\frac{i}{N},\frac{i+1}{N}\right]\times[0,1]$ or  $q$ ones of the form $[0,1]\times\left[\frac{i}{N},\frac{i+1}{N}\right]$.
\end{theorem}

In subcase (1.c) there are infinitely many line segments each of which is a component of $K$. Moreover, we have the following.
\begin{theorem}[{\cite[Corollary 2.6]{Lau-Luo-Rao13}}]\label{theo:single_slope}
If all components of $H^c$ are unbounded and if $K$ has both point components and non-degenerate ones then there is a constant rational $k_0$ such that every non-degenerate component is a line segment  of slope $k_0$.\end{theorem}

In  the subcases (1.a), (1.b), (1.c) and (2.a), it is transparent that all components of $K$ are locally connected. In subcase (2.b), it is also shown that every component of $K$ is locally connected. 
\begin{theorem}[{\bf\cite[Theorem 2]{LRX-2022}}]\label{theo:LRX}
Every component of $K$ is a Peano continuum.
\end{theorem}

Based on the core decomposition of planar compact sets (obtained in \cite{LoridantLuoYang19}), Feng et al introduce in \cite{FLY-2022} the notion of lambda function for all compact sets in the plane. In particular, one may employ the lambda function in studying the topology of fractal squares. To illustrate that, let us prepare some terminology.

For any compact set $K\subset\mathbb{R}^2$, let $\mathfrak{M}_K$ consist of all the upper semi-continuous (shortly, usc) decompositions of $K$ that satisfy two requirements: (1) every element is a sub-continuum of $K$, (2) the resulting quotient space is a Peano compactum. By \cite[Theorem 7]{LoridantLuoYang19}, there exists an usc decomposition $\mathcal{D}_K^{PC}\in\mathfrak{M}_K$ that refines every other $\mathcal{D}\in\mathfrak{M}_K$. Here we note that  a decomposition $\mathcal{D}$ of a compact metric space $X$ corresponds to an equivalence relation $\sim$ on $X$ whose classes are the elements of $\mathcal{D}$.  Moreover, the equivalence relation $\sim$ is closed, as a subset of $X^2$, if and only if $\mathcal{D}$ is an usc decomposition. See \cite[pp.98-99, Theorems 11 and 12]{Kelley55}.
\begin{definition}\label{def:CD}
Call $\mathcal{D}_K^{PC}$ the core decomposition of $K$. An element of $\mathcal{D}_K^{PC}$ is referred to as an (order-one) atom of $K$. If $\delta\in \mathcal{D}_K^{PC}$ then every atom of $\delta$ is called an order-two atom of $K$. Similarly, every atom of an order-$n$ atom is called an order-$(n+1)$ atom of $K$.
\end{definition}

Now, from the nets of ``atoms within atoms'' one may define the lambda function as in \cite{FLY-2022}. Given a compact set $K\subset\mathbb{R}^2$, we define its {\em lambda function} $\lambda_K:\mathbb{R}^2\rightarrow\mathbb{N}\cup\{\infty\}$ as follows.

First, set $\lambda_K(x)=0$ for all $x\notin K$ and all $x$ such that $\{x\}$ is an order-one atom of $K$. Second, let $\lambda_K(x)=m-1$ for all $x\in K$ if there is a minimal integer $m\ge2$ such that $\{x\}$ is an order-$m$ atom of $K$. Finally, set $\lambda_K(x)=\infty$ for $x\in K$ if such an integer $m$ does not exist.

\begin{definition}\label{def:lambda_range}
Call  $\lambda_K(K)=\{\lambda_K(x): x\in K\}$ the {\em lambda range} of $K$.
\end{definition}

The lambda function is useful in the study of plane topology. For instance, a compact set $K\subset\mathbb{R}^2$ is a Peano compactum if and only if its lambda function vanishes everywhere. This is equivalent to  $\lambda_K(K)=\{0\}$. In particular, a planar continuum $K$ is a Peano continuum if and only if  $\lambda_K(K)=\{0\}$. The classical Torhorst Theorem \cite[p.106]{Whyburn42} states that, for any Peano continuum $M\subset\mathbb{R}^2$ and any component $U$ of $\mathbb{R}^2\setminus M$, the boundary $\partial U$ is also a Peano continuum. Below we quote a quantified version of Torhorst Theorem from a recent study.
\begin{theorem}[{\bf \cite[Theorem 2]{FLY-2022}}]\label{theo:Torhorst}
For any compact  $K\subset\mathbb{R}^2$, any component $U$ of \  $\mathbb{R}^2\!\setminus\!K$ and any compact $L\subset\partial U$,  $\lambda_L(x)\le\lambda_{\partial U}(x)\le\lambda_K(x)$ holds for all $x\in\mathbb{R}^2$.
\end{theorem}

Now we apply the lambda function to the study of fractal squares.
\begin{theorem}[{\bf \cite[Theorems 2 and 3]{CKLY-2025}}]\label{theo:lambda_range}
Given a fractal square $K$,  we have $\lambda_K(K)\subset\{0,1\}$. Moreover, $\lambda_K(K)=\{1\}$ if and only if $K$ is the product of a Cantor set and $[0,1]$.
\end{theorem}

It follows that for any fractal square $K$ there are three possible cases. First,  $\lambda_K(K)=\{1\}$ and $K$ is the product of a linear Cantor set and $[0,1]$. Second,  $K$ is a Peano compactum hence $\lambda_K(K)=\{0\}$. 
Third, $K$ is not a Peano compactum and $\lambda_K(K)=\{0,1\}$. In such a case, all the non-degenerate components of $K$ form an infinite family of line segments with the same slope; moreover, there is a constant $C>0$ such that infinitely many of those line segments are of diameter $\ge C$.
Fractal squares $K$ with $\lambda_K(K)=\{1\}$ is characterized in the second part of Theorem  \ref{theo:lambda_range}. Fractal squares $K$ with $\lambda_K(K)=\{1\}$ are analyzed in \cite{CKLY-2025}. In Example \ref{exmp:lambda} we recall two simple fractal squares with $\lambda_K(K)=\{0, 1\}$, which have been discussed in \cite[Example 1.5]{CKLY-2025}. 
\begin{example}\label{exmp:lambda}
We give two fractal squares $K=K(5,\mathcal{D})$ 
with $\lambda_K(K)=\{0,1\}$, whose digit sets are respectively  $\mathcal{D}=\mathcal{D}_A$ or $\mathcal{D}_B$. 
In  Figure \ref{fig:non_PC} we illustrate the first and the fourth approximations, for $K(5,\mathcal{D}_A)$ and $K(5,\mathcal{D}_B)$, that agree with the underlying attractor to the best of our eyesight.
We can verify that for $K=K(5,\mathcal{D}_A)$ the level set  $\lambda_K^{-1}(1)$ contains infinitely many line segments and its Hausdorff dimension $\dim_H\lambda_K^{-1}(1)$ is one. On the other hand, for $K=K(5,\mathcal{D}_B)$ the level set  $\lambda_K^{-1}(1)$ contains uncountably many line segments and satisfies $1<\dim_H\lambda_K^{-1}(1)<\dim_HK$.
\begin{figure}[ht]    
\begin{center}
\begin{tabular}{cccc}
\includegraphics[width=2.75cm]{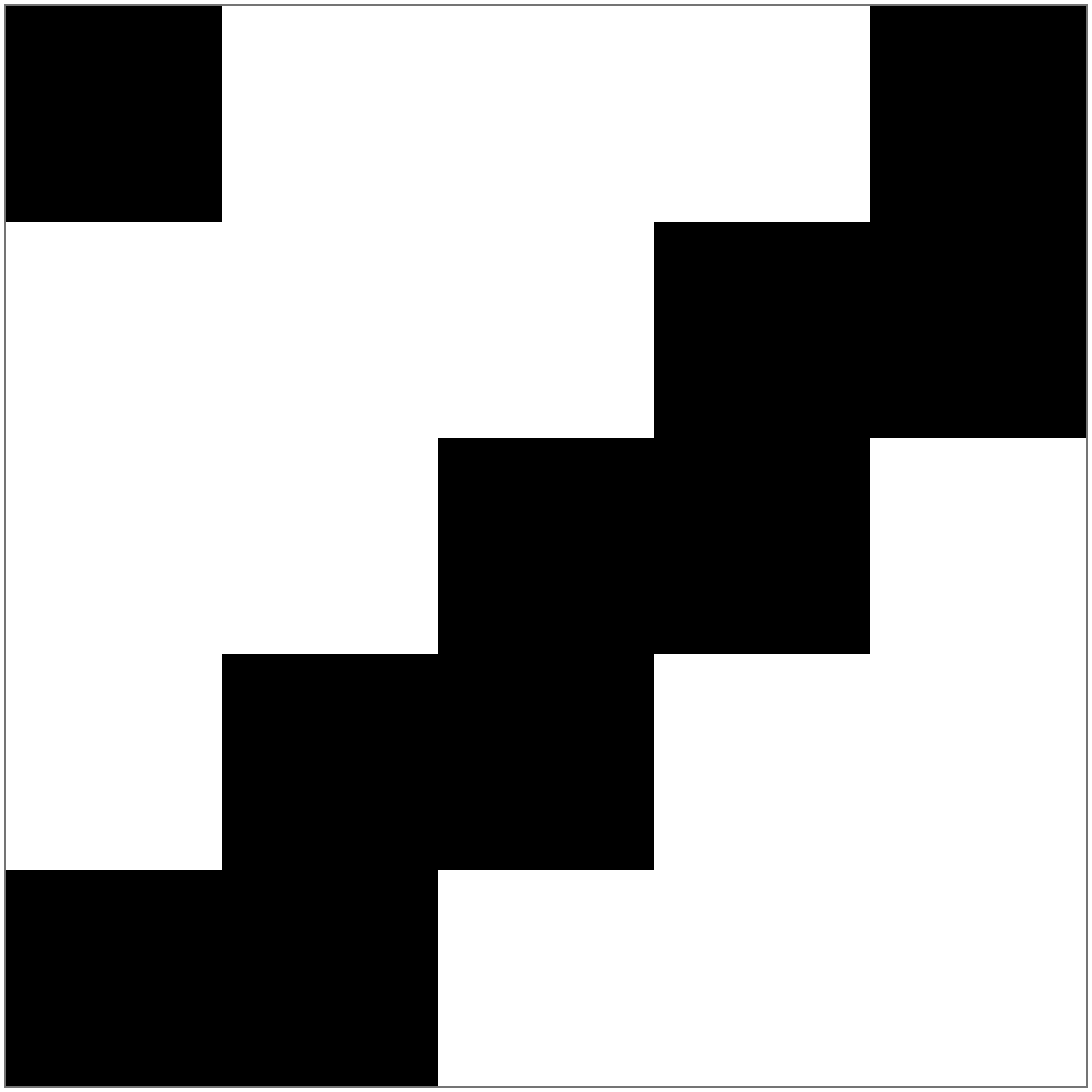} & 
\includegraphics[width=2.75cm]{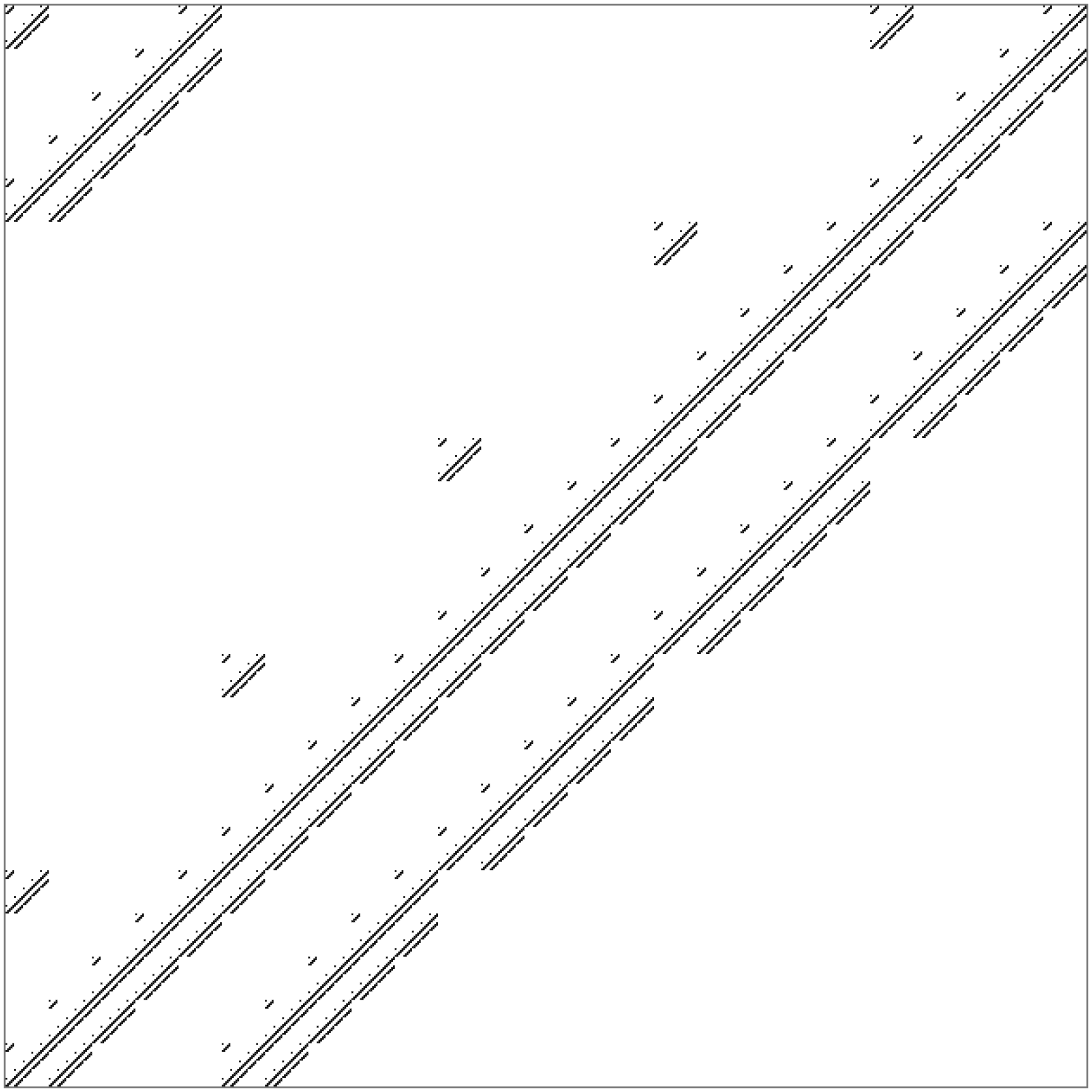}&
\includegraphics[width=2.75cm]{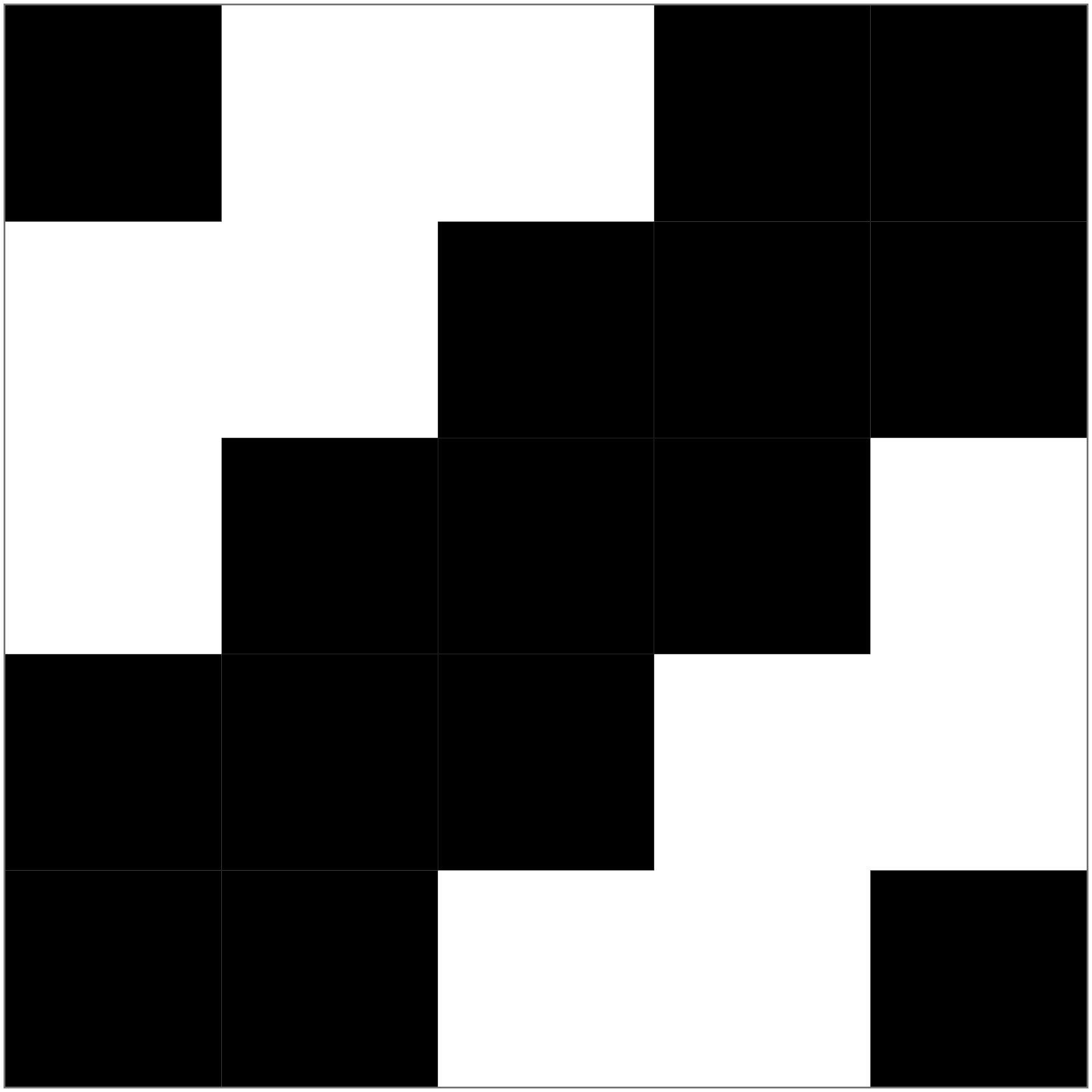} & 
\includegraphics[width=2.75cm]{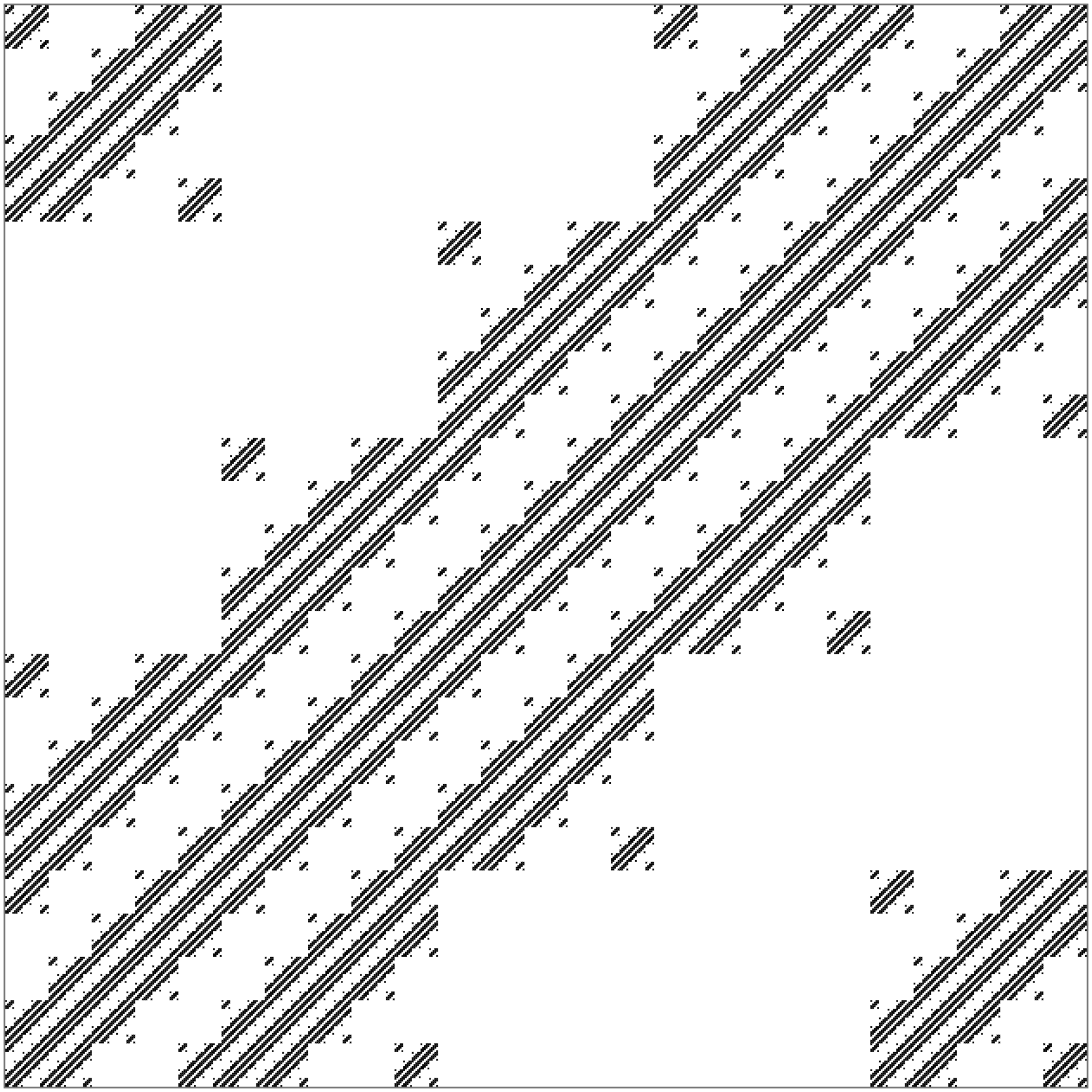}
\end{tabular}
\end{center}
\vspace{-0.25cm}
\caption{The first and the fourth approximations for $K(5,\mathcal{D}_A)$ and $K(5,\mathcal{D}_B)$ in Example \ref{exmp:lambda}.}\label{fig:non_PC}
\end{figure}
\end{example}

The following Question \ref{que:lambda}  is of natural interest.

\begin{question}\label{que:lambda}
Under what conditions,  on the approximations $K^{(j)}$ with small $j$, do we have $\lambda_K(K)=\{0,1\}$? Moreover,  when do we have
$\dim_H\lambda_K^{-1}(1)=1$ or $>1$?
\end{question}

To conclude this section, we want to mention a very basic question in response to  Theorem \ref{theo:single_slope}. For related details, one may refer to \cite[Corollary 2.6]{Lau-Luo-Rao13}. 
 
 Given $N\ge3$, let $\mathcal{L}_N$ consist  of $0$, $1$ and all numbers $\tau\in(0,1)$ such that there is a fractal square $K=K(N,\mathcal{D})$ whose non-degenerate components are line segments of the same slope $\tau$. Such a slope $\tau$ is necessarily rational.  Moreover, let $\mathcal{L}_N^\#$ consist  of  $0$, $1$, and all numbers $\tau\in(0,1)$ such that there is a fractal square $K=K(N,\mathcal{D})$ with $\mathcal{D}\ne\{0,\ldots,N-1\}^2$ that contains a line segment of  slope $\tau$. 
 By the symmetry of $K$, the set $\mathcal{L}_N^\#$ essentially gives all the possible slopes of line segments contained in $K$, while $\mathcal{L}_N$ gives the same information for those $K$ whose non-degenerate components are line segments. 
 
 Write  every rational $\frac{r}{s}\in\mathcal{L}_N^\#$ in reduced form, so that the numerator $r$ and the denominator $s$ are coprime. Then the denominator $s$ is strictly smaller than $N$. Indeed, if $K$ contains a line segment with slope $\frac{r}{s}<1$ for some $s\ge N$, then $H=K+\mathbb{Z}^2$ contains an infinite line of the same slope, which implies that $\mathcal{D}=\{0,\ldots,N-1\}^2$. 
 
 We propose the following.

\begin{question}\label{algo:slope}
How to determine  $\mathcal{L}_N$ and $\mathcal{L}_N^\#$ for  $N\ge3$? Particularly, given $\mathcal{D}\subset\{0,\ldots,N-1\}^2$, how to determine all $\tau\in\mathcal{L}_N^\#$ such that $K(N,\mathcal{D})$ contains a line segment of slope $\tau$? Moreover, is it always true that $\mathcal{L}_N\subsetneq\mathcal{L}_N^\#$?
\end{question}

\section{Topological Classification for $N=3$}\label{sec:N=3}

In this section we focus on the fractal squares $K(3,\mathcal{D})$ for all $\mathcal{D}\subset\{0,1,2\}^2$ with $2\le \#\mathcal{D}\le 8$. Based on recent studies, we have a nearly complete classification of those fractal squares in terms of topological equivalence.

For the sake of convenience, let $\chi(q)$ be the number of topological equivalence classes in the family of all $K(3,\mathcal{D})$ with 
$\#\mathcal{D}=q\in\{2,\ldots, 8\}$, to be denoted by $\mathcal{K}_{3,q}$. 
Knowledge of those numbers $\chi(q)$ may be summarized in the table below.

\begin{center}
\begin{tabular}{l|c|c|c|c|c|c|c}
$q$&\hspace{0.25cm}        $2$  \hspace{0.25cm}&\hspace{0.25cm} 
$3$ \hspace{0.25cm}&\hspace{0.25cm} $4$ 
\hspace{0.25cm}&\hspace{0.25cm} $5$ \hspace{0.25cm}&\hspace{0.25cm} $6$  \hspace{0.25cm}&\hspace{0.25cm}
$7$ \hspace{0.25cm}&\hspace{0.25cm} $8$ \\ \hline
$\chi(q)$&\hspace{0.25cm} $1$ \hspace{0.25cm}&\hspace{0.25cm} 
$2$ \hspace{0.25cm}&\hspace{0.25cm} $2$
\hspace{0.25cm}&\hspace{0.25cm} $\le8$ \hspace{0.25cm}&\hspace{0.25cm} $13$ \hspace{0.25cm}&\hspace{0.25cm} $8$
\hspace{0.25cm}&\hspace{0.25cm} $2$
\end{tabular}
\end{center}

The topology of $K\in\mathcal{K}_{3,q}$ is not complicated for small $q$. It is trivial that every $K\in\mathcal{K}_{3,2}$  is a Cantor set hence we have $\chi(2)=1$. 
It is also transparent that every $K\in\mathcal{K}_{3,3}$  is either a Cantor set or a line segment, whose slope may be $0$, $\infty$, or $\pm1$. 
Moreover,  every $K\in\mathcal{K}_{3,4}$  is either a Cantor set or homeomorphic with $K_0=K(3,\mathcal{D}_0)$ with $\mathcal{D}_0=\{(i,i): 0\le i\le 2\}\cup\{(2,0)\}$. See the left most two pictures in Figure \ref{fig:q=4}. Notice that $K_0$  is a Peano compactum containing infinitely many line segments of slope $1$. 

\begin{example}\label{exmp:N=4}
Let $K=K(3,\mathcal{D})$ be a fractal square of order three, where the digit set $\mathcal{D}\subset\{0,1,2\}^2$ contains $4$ points and is given as in Figure \ref{fig:q=4}. Then $K$ is  homeomorphic with $K_0$.
   
\begin{figure}[ht]    
\begin{center}
\begin{tabular}{cccc}
\includegraphics[width=2.75cm]{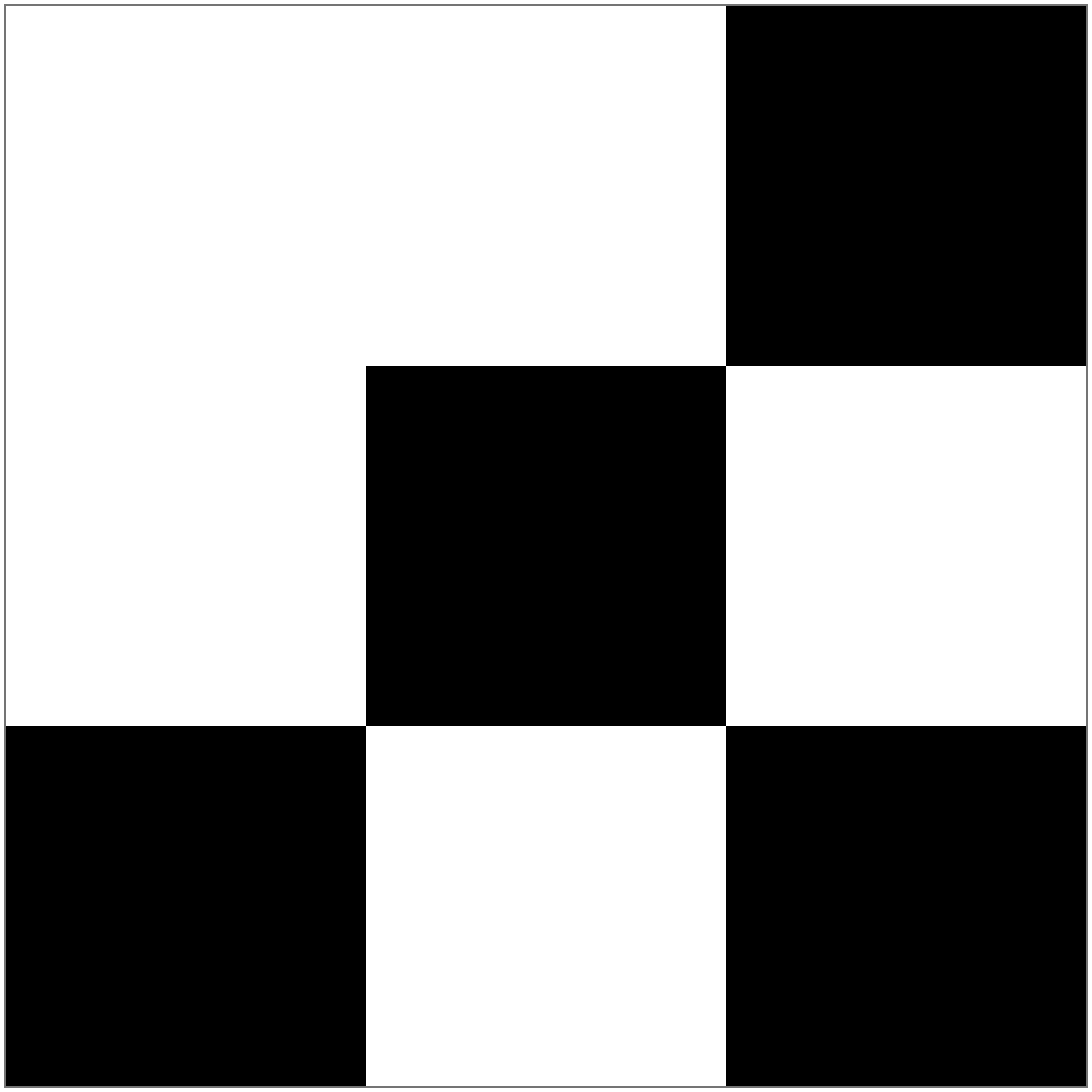} &
\includegraphics[width=2.75cm]{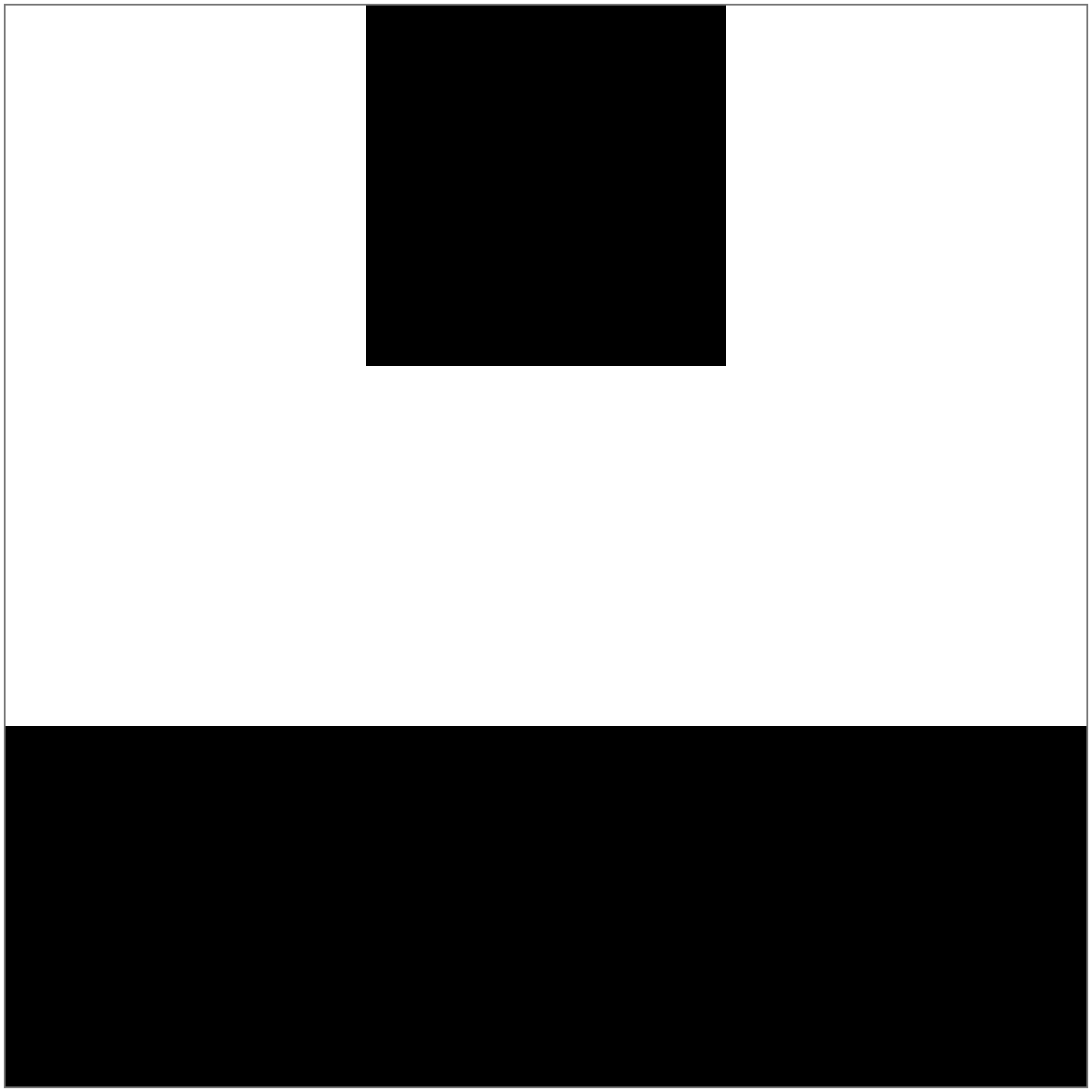}&  
\includegraphics[width=2.75cm]{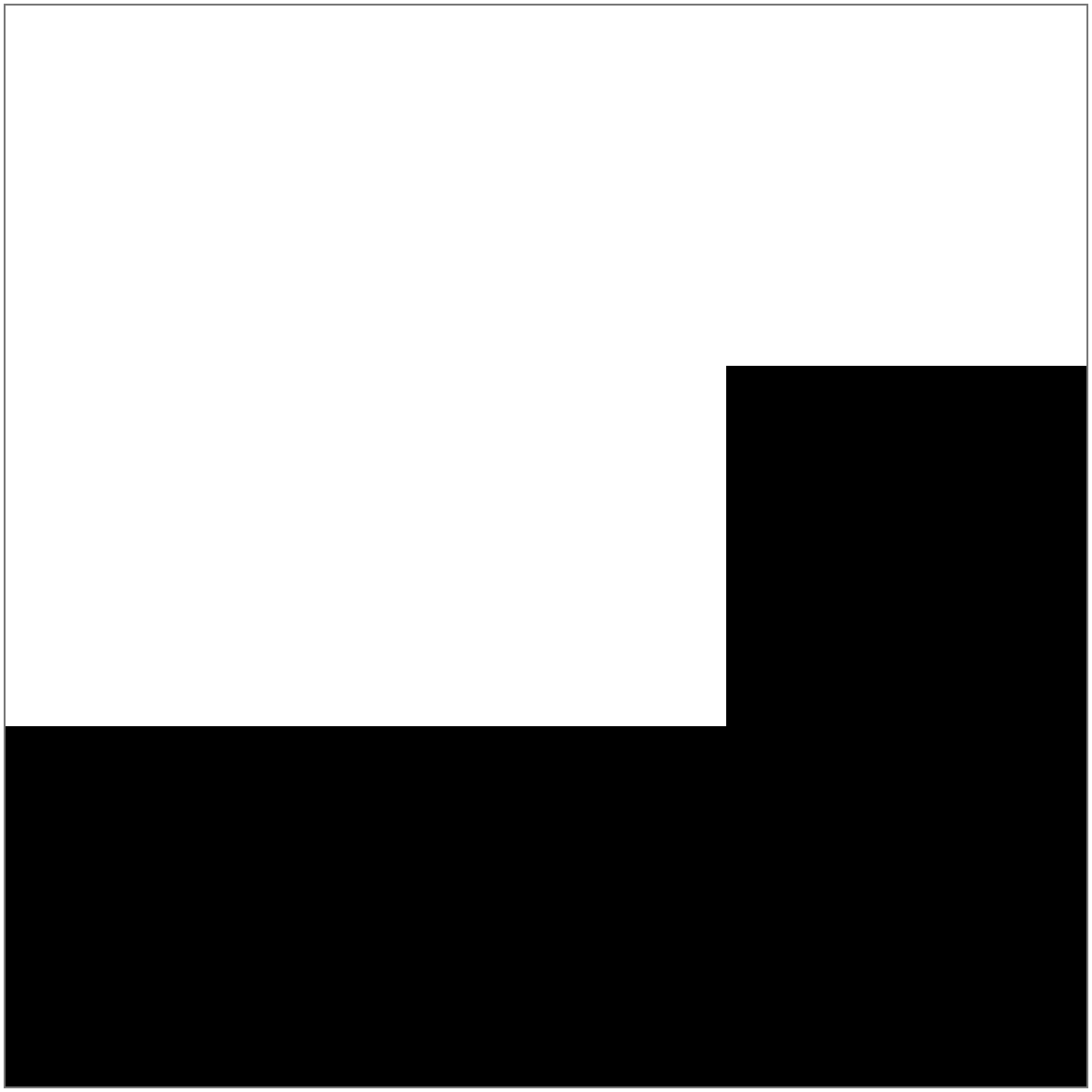}&
\includegraphics[width=2.75cm]{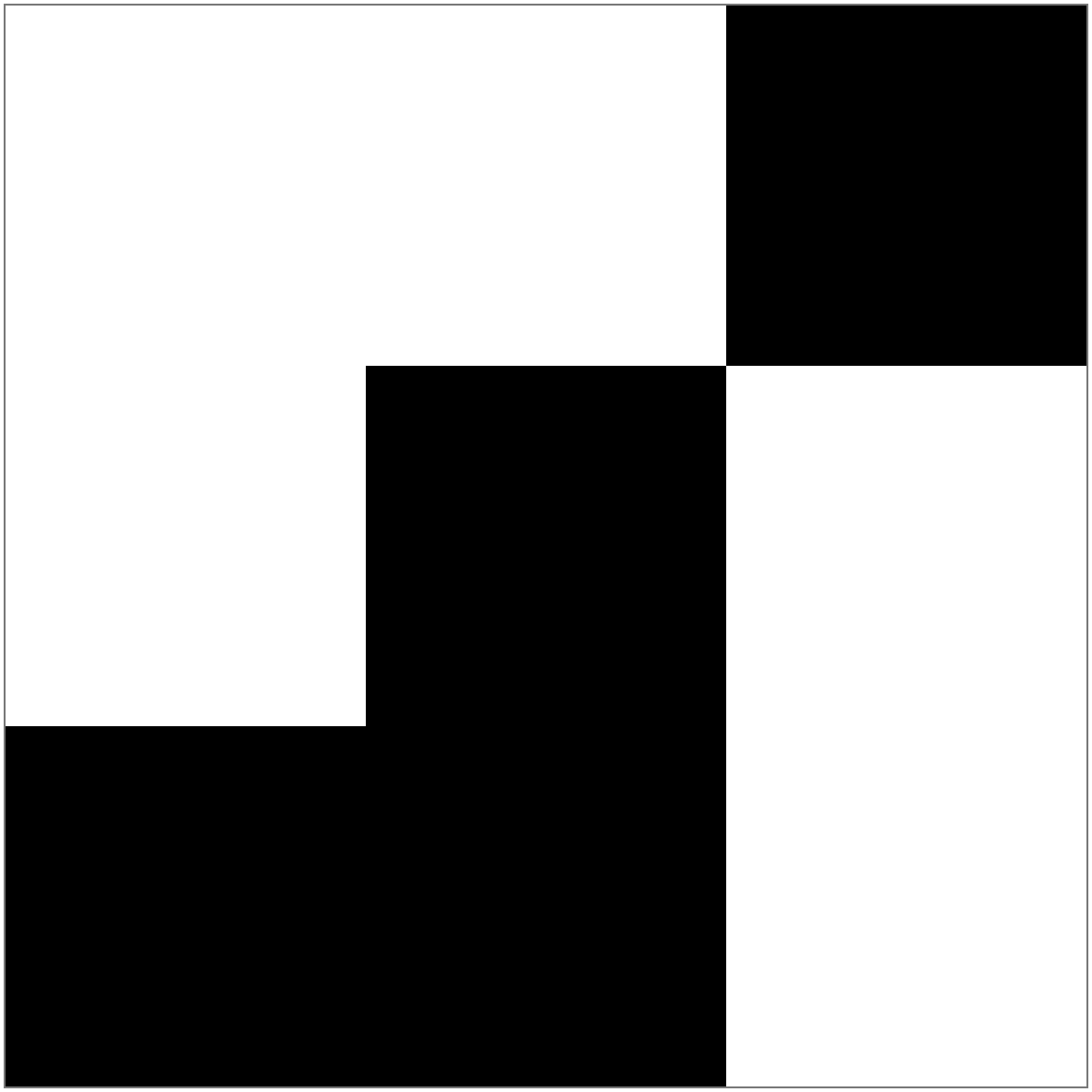}
 \\
 \includegraphics[width=2.75cm]{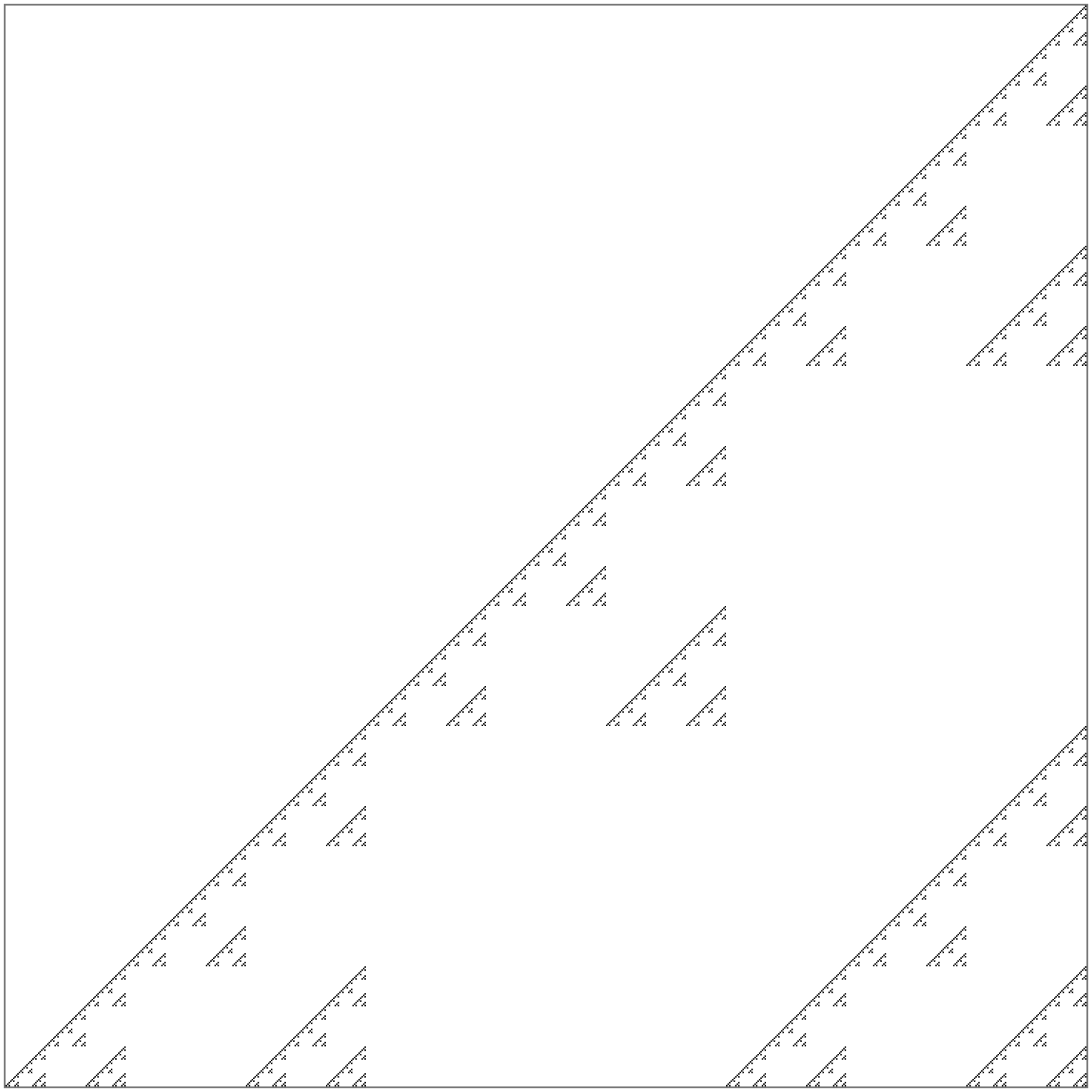} &
 \includegraphics[width=2.75cm]{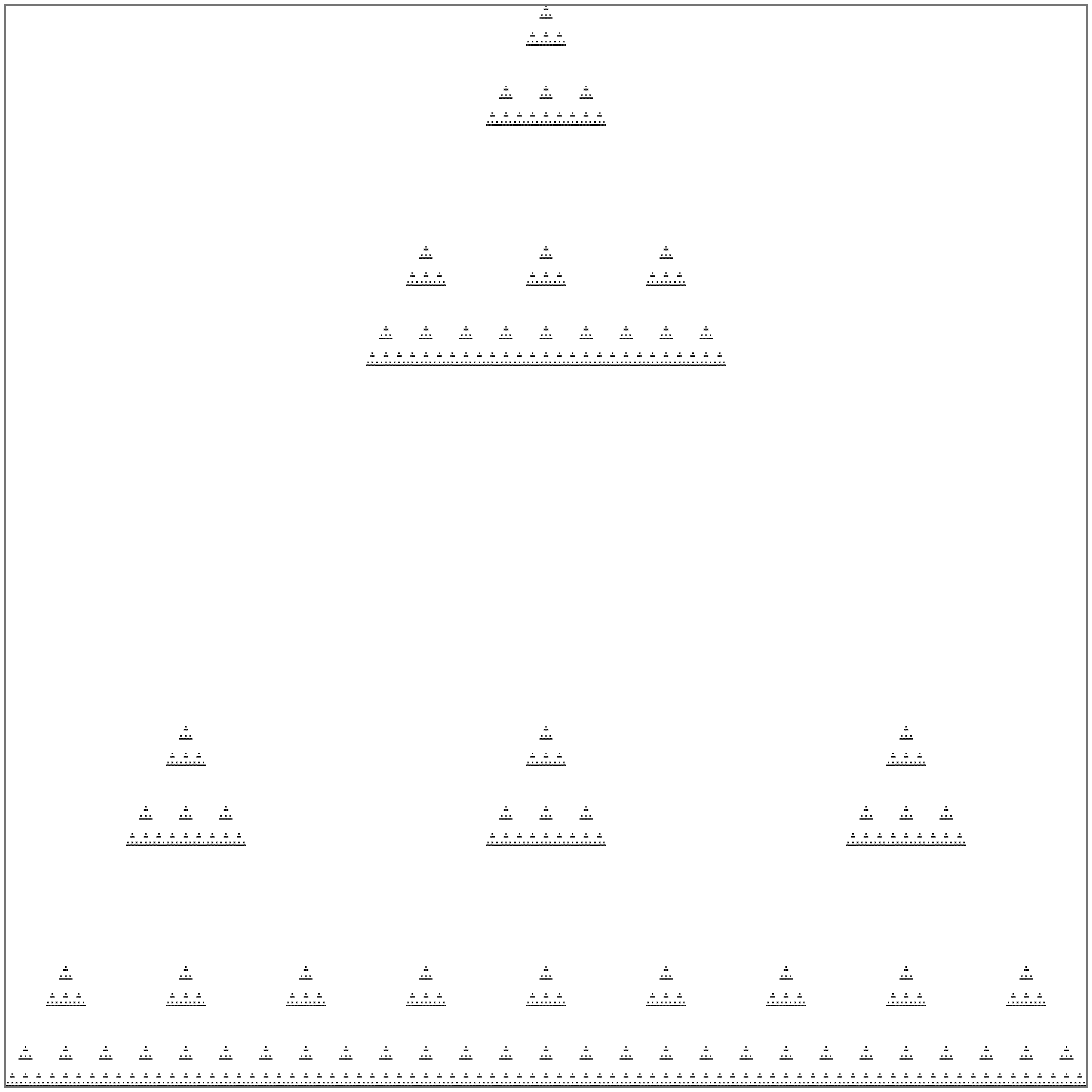}&  
 \includegraphics[width=2.75cm]{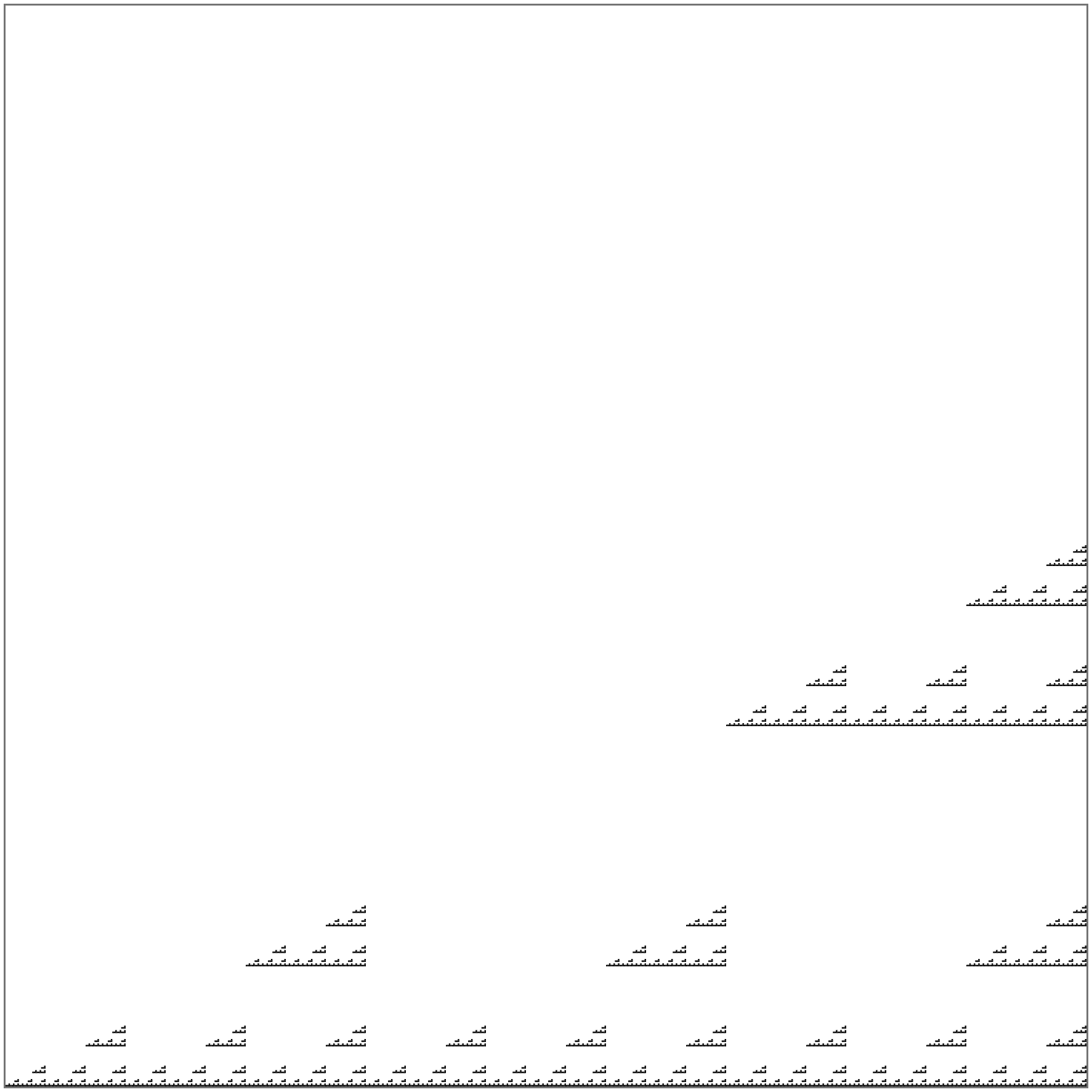}&
 \includegraphics[width=2.75cm]{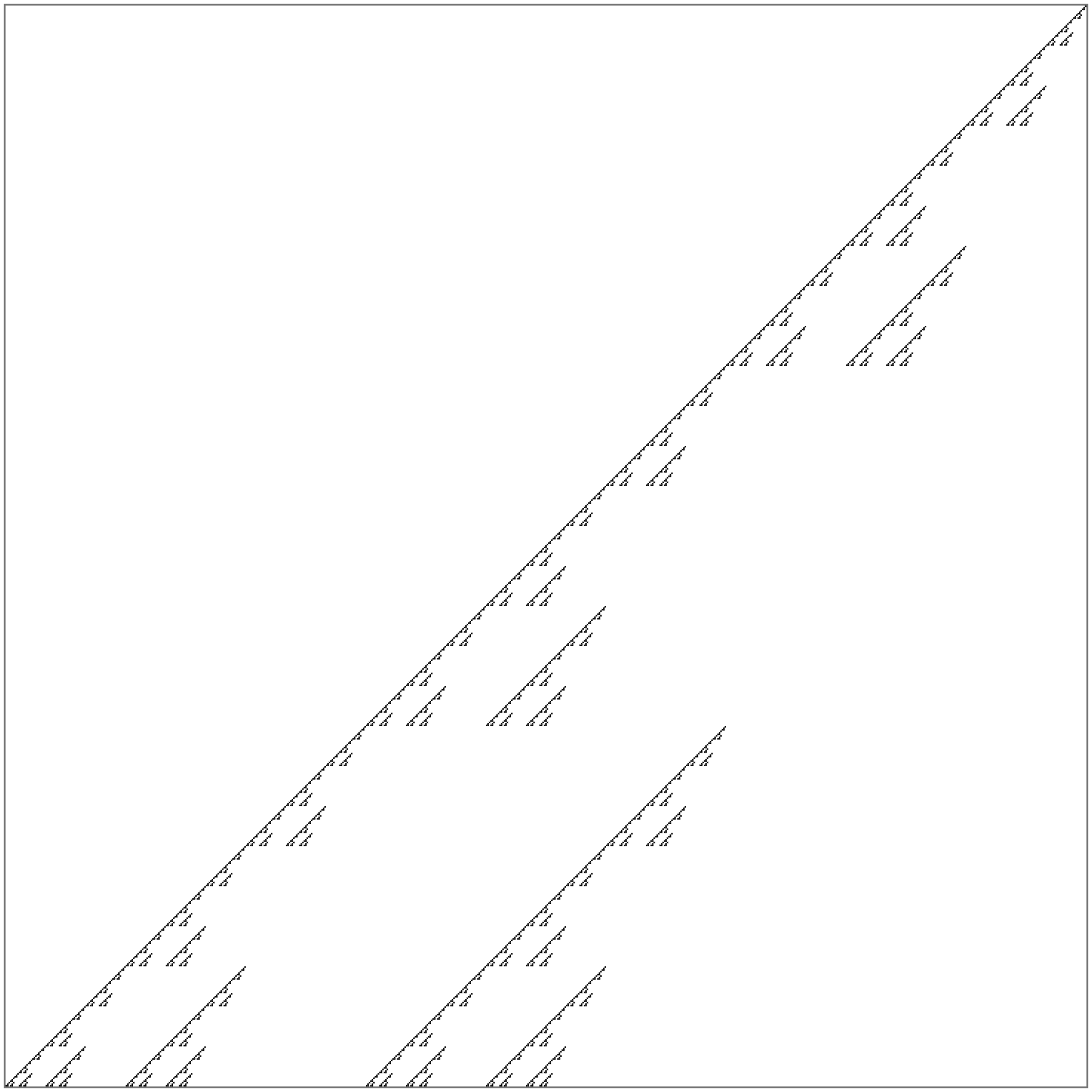}
\end{tabular}
\end{center}
\caption{Leftmost column: the first and the sixth approximations of $K_0=K\left(3,\mathcal{D}_0\right)$. Other columns: the same approximations for three other fractal squares. The four attractors are homeomorphic.}\label{fig:q=4}
\end{figure}
\end{example}

Therefore, we have.
\begin{theorem}[{\bf \cite[Theorem 3.3]{Luo-Liu16} and \cite[Theorem 4.5]{WZD12}}]\label{theo:q=3/4}
$\chi(3)=\chi(4)=2$.
\end{theorem}

The fractal squares $K(3,\mathcal{D})$ with $\#\mathcal{D}=5$ are  divided into $21$ congruent classes. See \cite{Luo-Liu16} for a complete list.

For five of them, as illustrated in \cite[Fig. A.1]{Luo-Liu16}, $K(3,\mathcal{D})$ is a Cantor set.

For six of them, as illustrated in \cite[Fig. A.2]{Luo-Liu16},  $K(3,\mathcal{D})$ is connected hence is a Peano continuum. Two of those continua are homeomorphic, with nontrivial fundamental group. The other four are dendrites, such as $K(3,\mathcal{D}_i)(i=2,3)$ with 
\begin{eqnarray}\label{eq:D_T}
\mathcal{D}_2&=&\{(i,0):0\le i\le 2\}\cup\{(1,1), (1,2)\}\\
\mathcal{D}_3&=&\{(i,1):0\le i\le 2\}\cup\{(1,0), (1,2)\}
\end{eqnarray} 
It is known that $K(3,\mathcal{D}_2)$ and $K(3,\mathcal{D}_3)$ are not homeomorphic \cite[Theorem 3.10]{Luo-Liu16}. On the other hand, the other two are  homeomorphic with $K(3,\mathcal{D}_3)$. Here, by checking the digit sets one may figure out how to choose the underlying homeomorphisms to be even an affine map.

The rest ten, as illustrated in \cite[Fig. A.3]{Luo-Liu16}, have both point components and non-degenerate ones,  which are parallel line segments. The one with digit set $\mathcal{D}_1$ 
has been addressed in Example \ref{exmp:Hata_Graph}. Those with digit sets $\mathcal{D}_j\ (4\le j\le 7)$ are illustrated in Example \ref{exmp:N=3_classes}, where
\begin{eqnarray}\label{eq:D_nondendrite}
\mathcal{D}_4&=&\{(i,0):0\le i\le 2\}\cup\{(0,2),(1,2)\}\\
\mathcal{D}_5&=&\{(i,0):0\le i\le 2\}\cup\{(0,2), (2,2)\}\\
\mathcal{D}_6&=&\{(i,i):0\le i\le 2\}\cup\{(1,0), (0,2)\}\\
\mathcal{D}_7&=&\{(i,0):0\le i\le 2\}\cup\{(0,2), (2,1)\}
\end{eqnarray} 

\begin{example}\label{exmp:N=3_classes}
Let $K_j=K\left(3,\mathcal{D}_j\right)\ (4\le j\le 7)$ be the fractal square with digit sets $\mathcal{D}_j$, that are respectively  given in Equations (7)-(10). Figure \ref{fig:order_3} illustrates for $4\le j\le 7$ the approximations  $K_j^{(n)}$ with  $n=1,6$. By  \cite[Theorems 2.1 and 5.1]{ZhuYang18} and by a family of self-similar sets that are analyzed in \cite[Example 2.1]{ZhuYang18}, we know that $K_4$ differs from $K_1=K(3,\mathcal{D}_1)$ by a bi-Lipschitz map. 
\begin{figure}[ht]    
\begin{center}
\begin{tabular}{cccc}
\includegraphics[width=2.75cm]{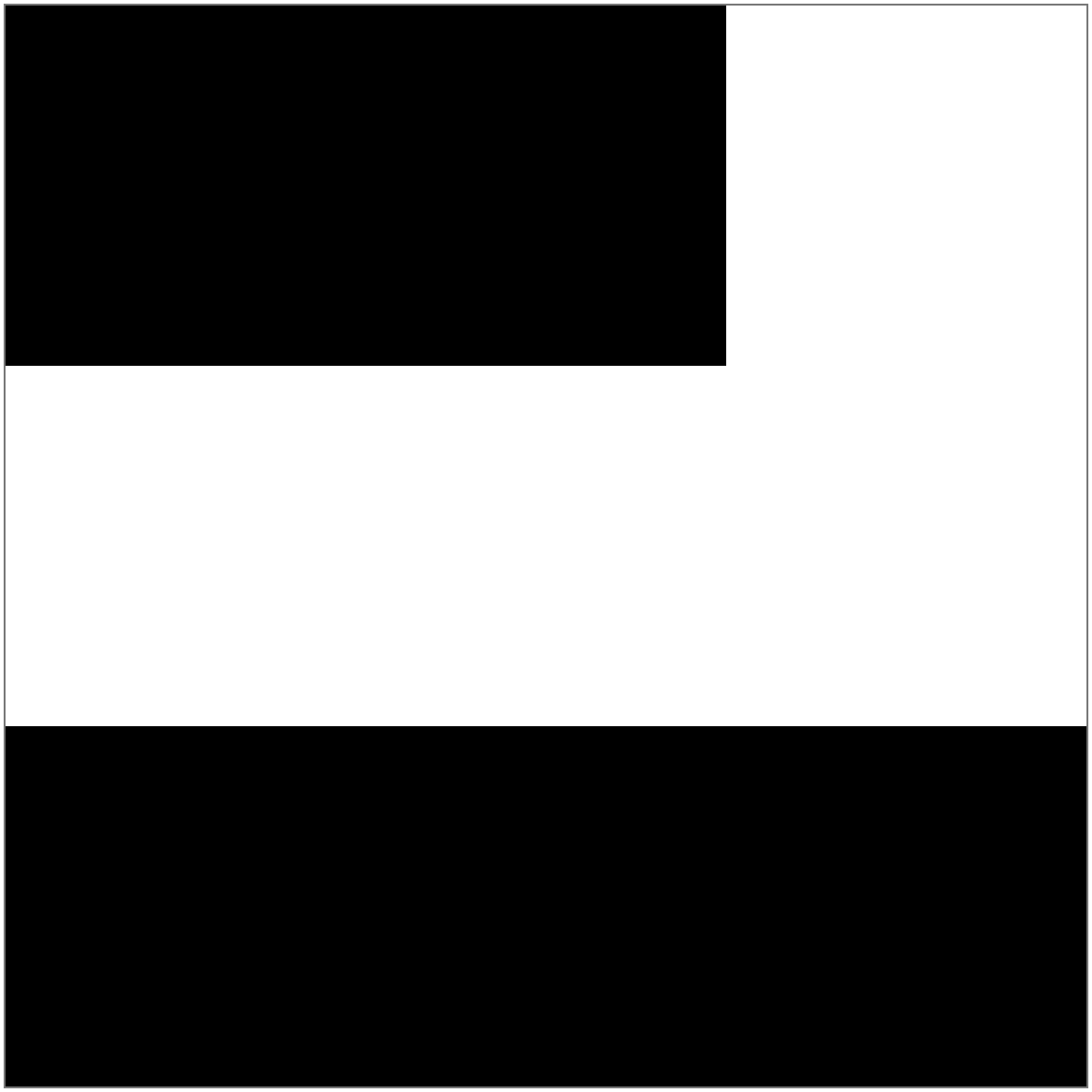}
&
\includegraphics[width=2.75cm]{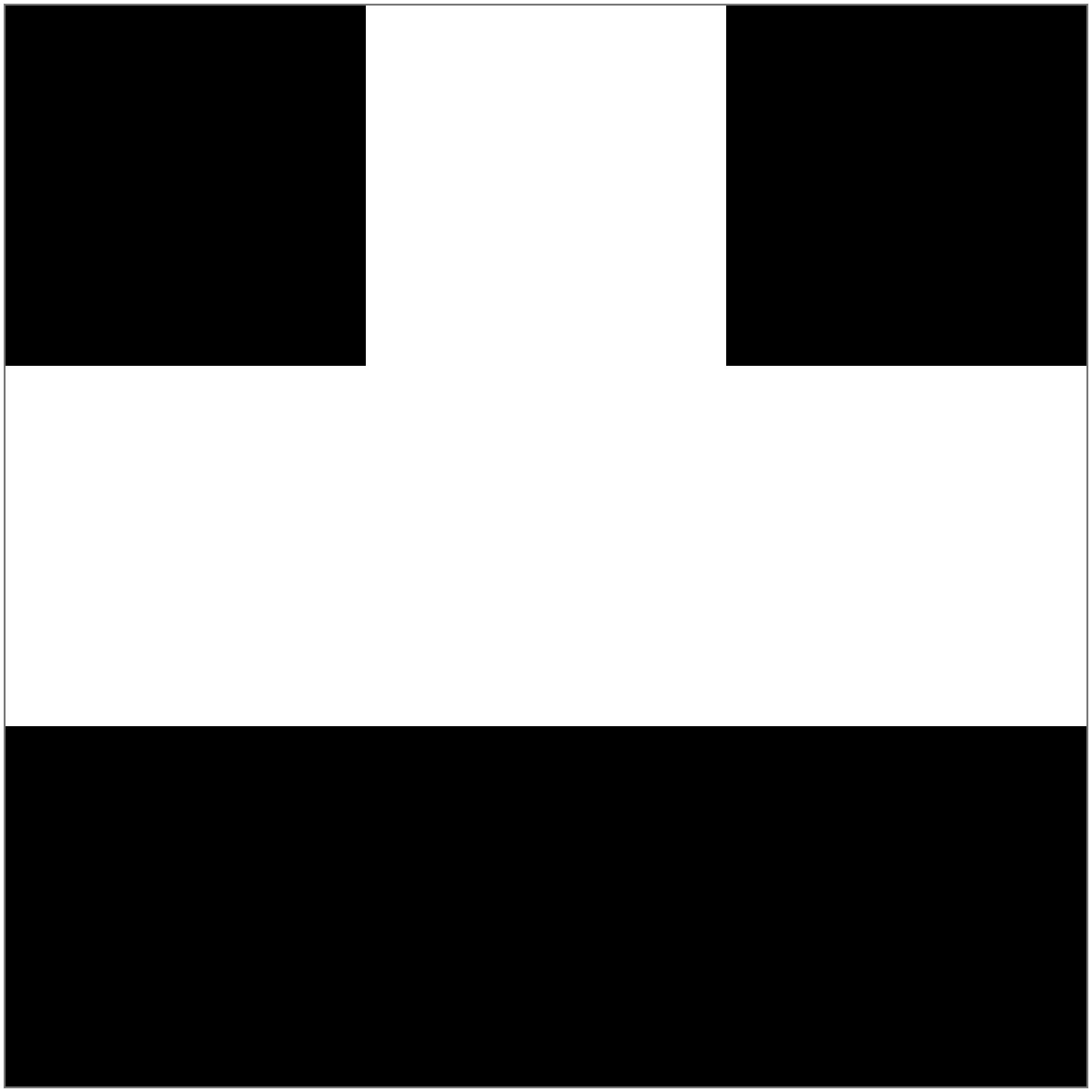}
&
\includegraphics[width=2.75cm]{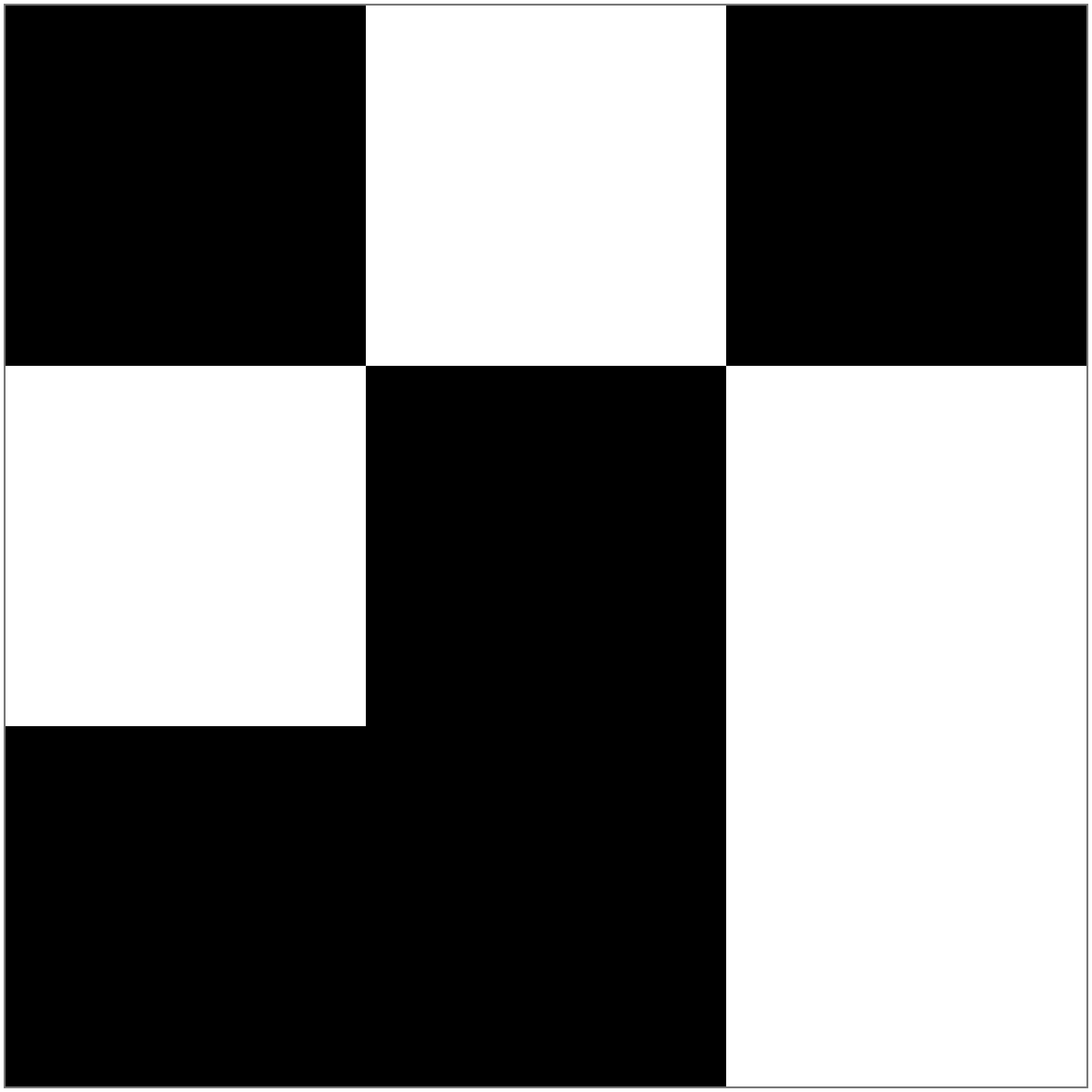}
&
\includegraphics[width=2.75cm]{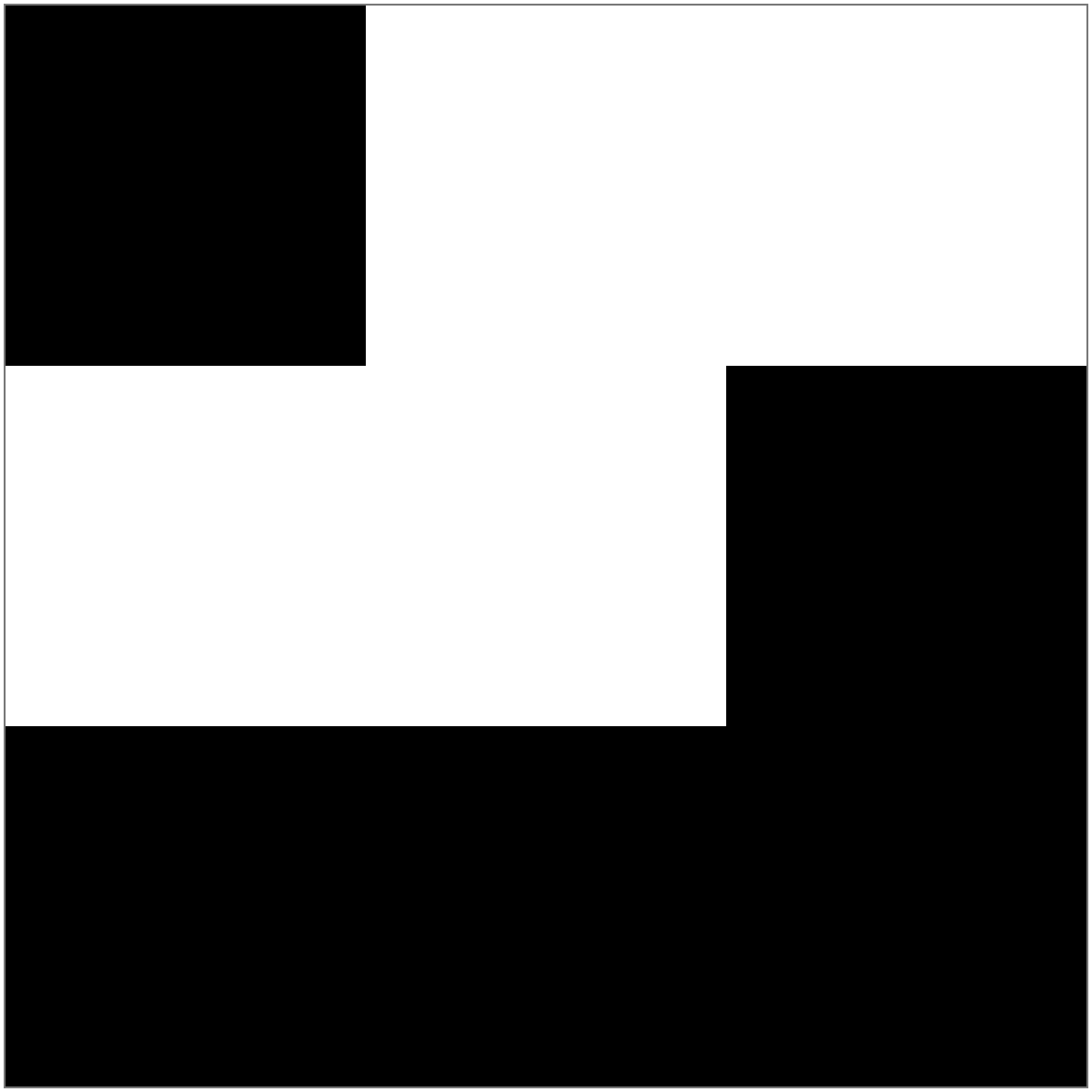}
\\
\includegraphics[width=2.75cm]{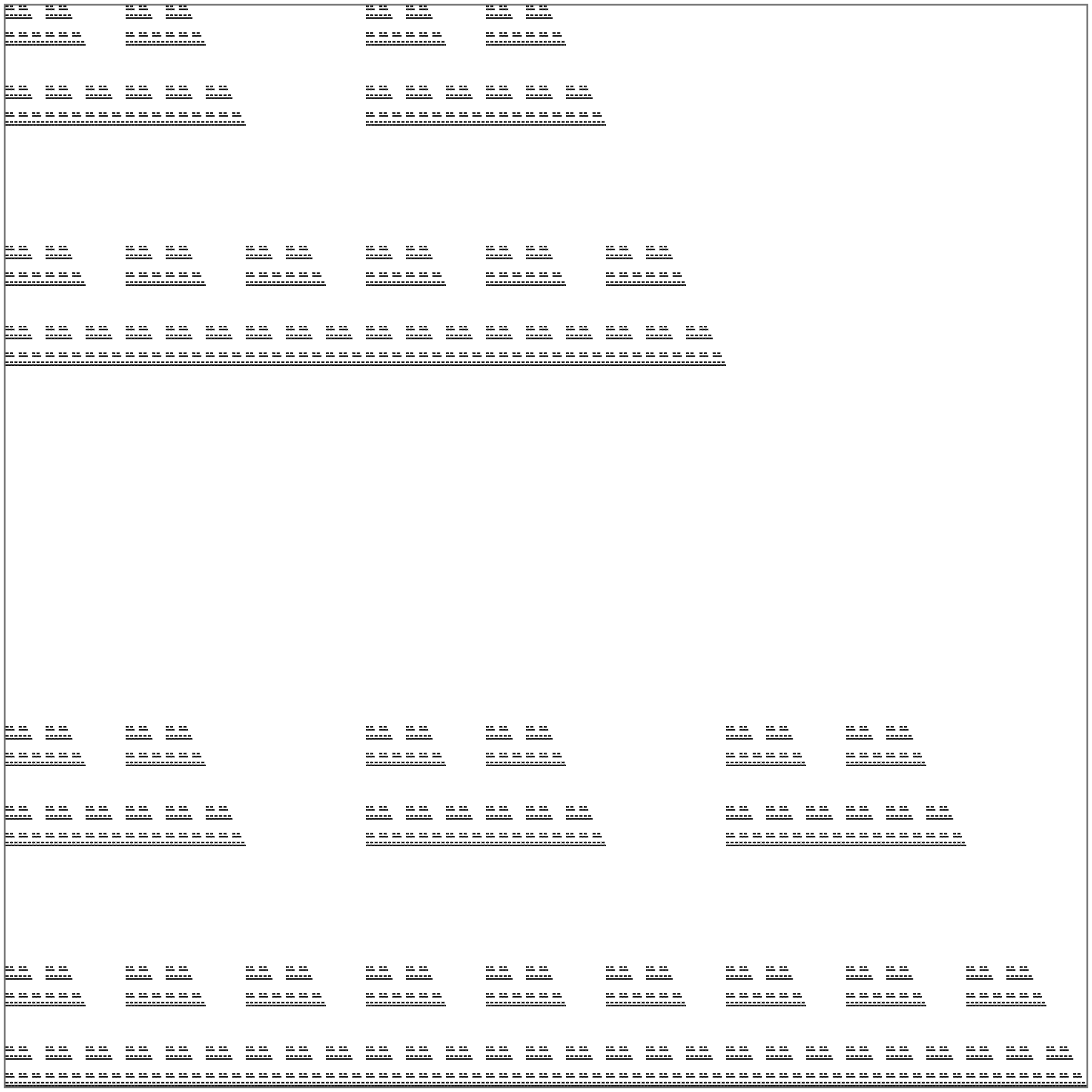}
&
\includegraphics[width=2.75cm]{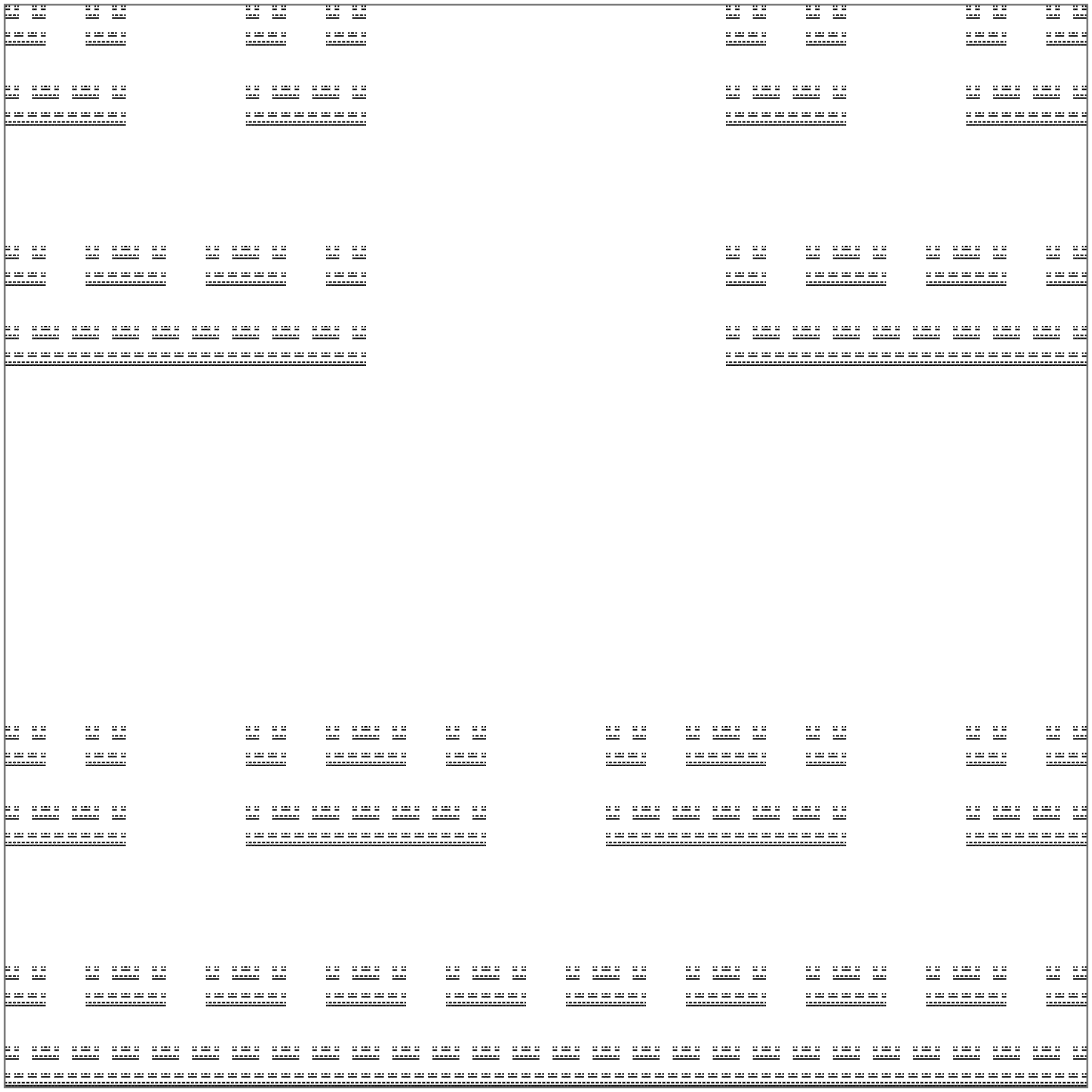}
&
\includegraphics[width=2.75cm]{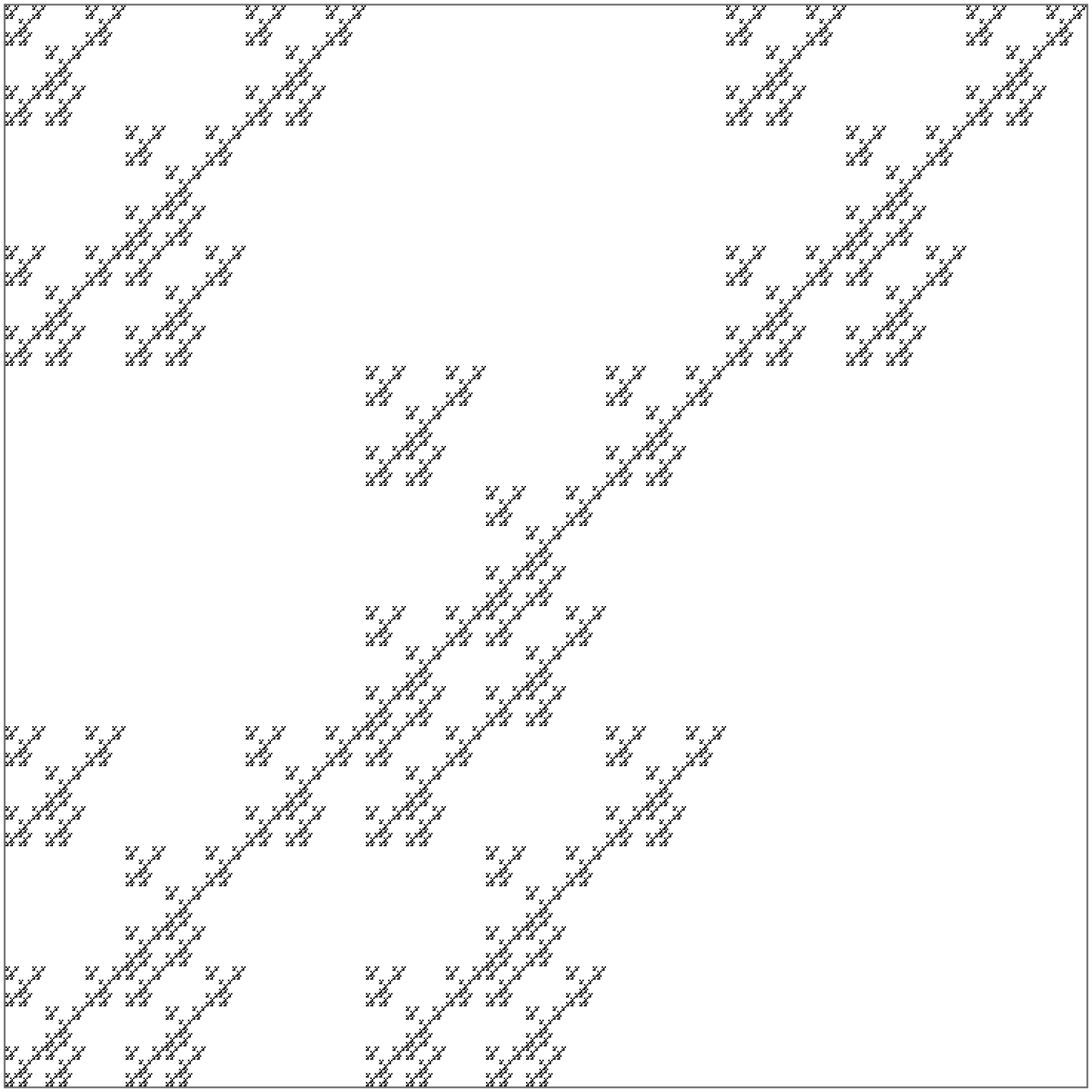}
&
\includegraphics[width=3.0cm]{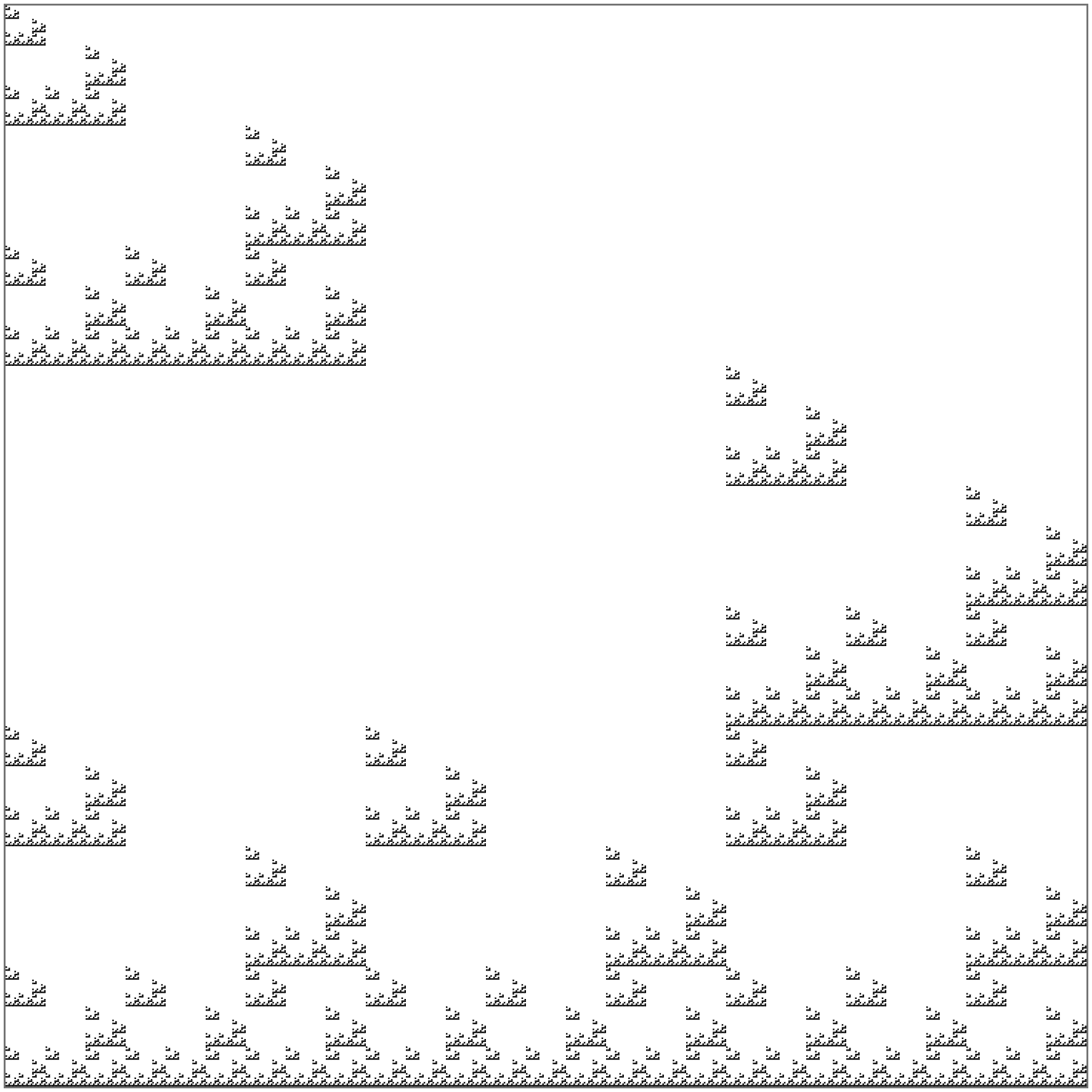}\\
$K_4^{(1)},\ K_4^{(6)}$ & $K_5^{(1)},\ K_5^{(6)}$ & $K_6^{(1)},\ K_6^{(6)}$ & $K_7^{(1)},\ K_7^{(6)}$
\end{tabular}\end{center}
\caption{The two approximations $K_j^{(n)}\ (n=1,6)$ for $K_j\ (4\le j\le 7)$.}\label{fig:order_3}
\end{figure}
\end{example}

\begin{example}\label{exmp:D1_D8}
Let $K_8=K\left(3,\mathcal{D}_8\right)$ be the fractal square  whose digit $\mathcal{D}_8$ is given by
\begin{eqnarray}\label{eq:D8}
\mathcal{D}_8&=&\{(i,0):0\le i\le 2\}\cup\{(1,2), (2,1)\}.
\end{eqnarray}
With the help of finite state automata  a bi-Lipschitz map is constructed in \cite{ZhuRao21} between $K(3,\mathcal{D}_1)$ given in Example \ref{exmp:Hata_Graph} and $K_8$.  Figure \ref{fig:topology_automaton} illustrates the approximations $K_8^{(1)}$ and $K_8^{(6)}$, together with $K_1^{(1)}$ and $K_1^{(6)}$.
\begin{figure}[ht]    
\begin{center}
\begin{tabular}{cccc}
\includegraphics[width=2.75cm]{PNG_3_D1_1}
&\includegraphics[width=2.75cm]{PNG_3_D1_6}
&\includegraphics[width=2.75cm]{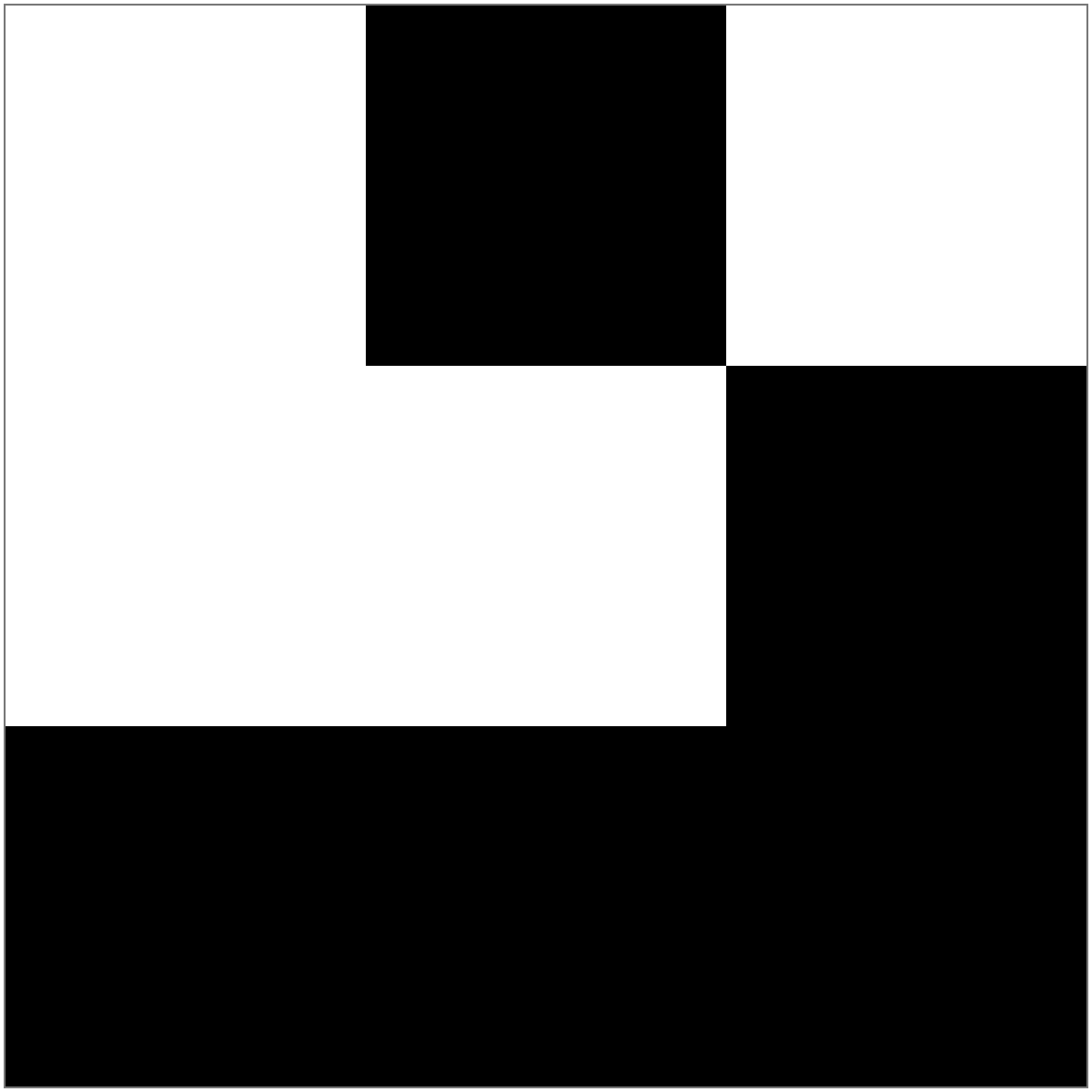}
&\includegraphics[width=2.75cm]{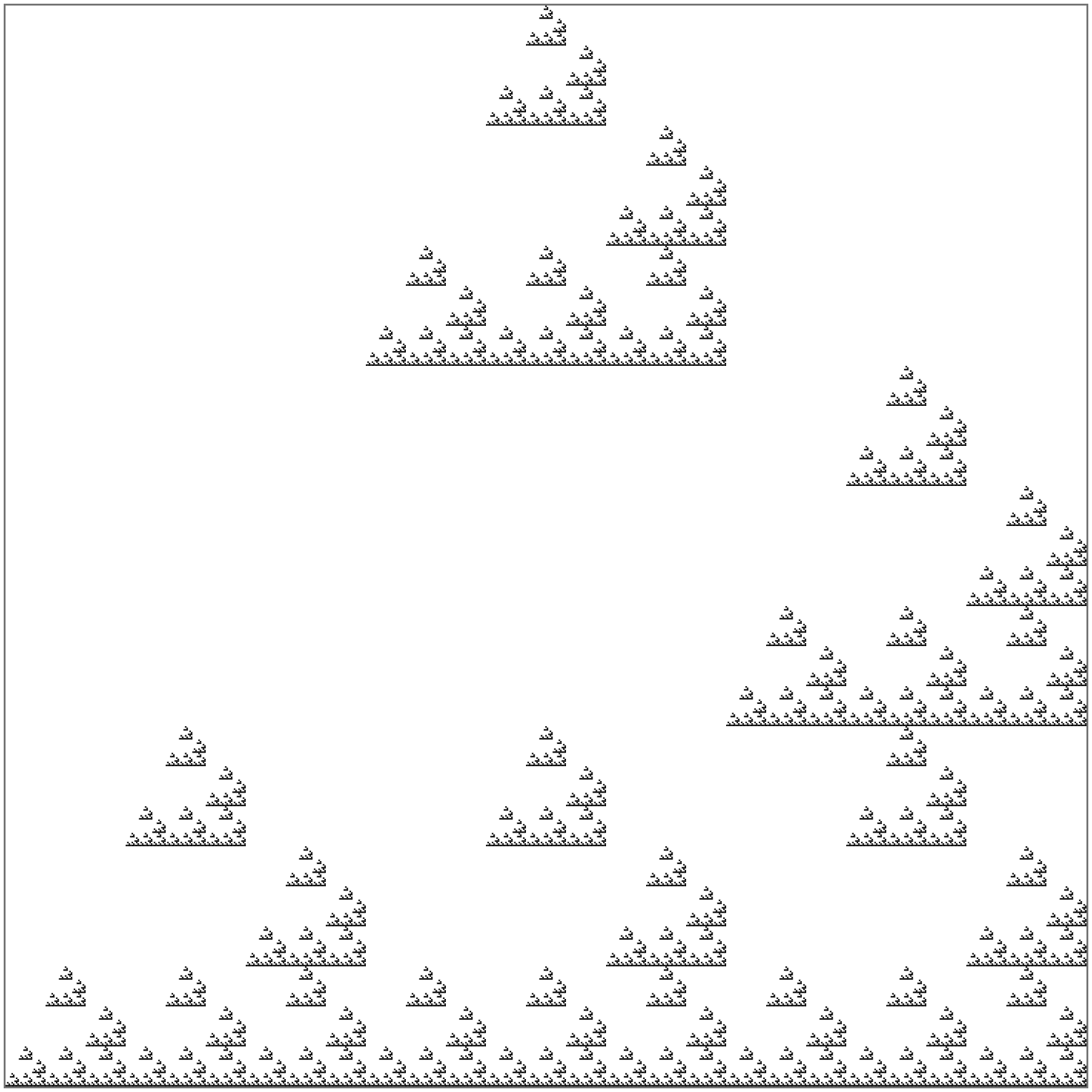}\\
$K_1^{(1)}$ & $K_1^{(6)}$ & $K_8^{(1)}$& $K_8^{(6)}$
\end{tabular}\end{center}
\caption{ The first and the sixth approximations for $K(3,\mathcal{D}_j)$ with $j=1,8$. }\label{fig:topology_automaton}
\end{figure}
\end{example}

\begin{example}\label{exmp:N=3_easy_homeo}
Each of the other four fractal squares in \cite[Fig. A.3]{Luo-Liu16} is homeomorphic to some $K_j$  with $4\le j\le7$. In deed, let $\mathcal{D}_j\ (9\le j\le12)$ be  given below. 
\begin{eqnarray}\label{eq:more_D_nondendrite}
\mathcal{D}_9 &=&\{(i,i):0\le i\le 2\}\cup\{(1,0), (0,1)\}\\
\mathcal{D}_{10}&=&\{(i,0):0\le i\le 2\}\cup\{(1,1), (2,1)\}\\
\mathcal{D}_{11}&=&\{(i,i):0\le i\le 2\}\cup\{(1,0), (2,1)\}\\
\mathcal{D}_{12}&=&\{(i,0):0\le i\le 2\}\cup\{(0,1), (2,1)\}
\end{eqnarray}
Let $K_j=K\left(3,\mathcal{D}_j\right)$. Then, one may check that $K_9$ is homeomorphic with $K_1$, both $K_{10}$ and $K_{11}$ are homeomorphic with $K_4$, while $K_{12}$ is homeomorphic with $K_5$.
\begin{figure}[ht]    
\begin{center}
\begin{tabular}{cccc}
\includegraphics[width=2.75cm]{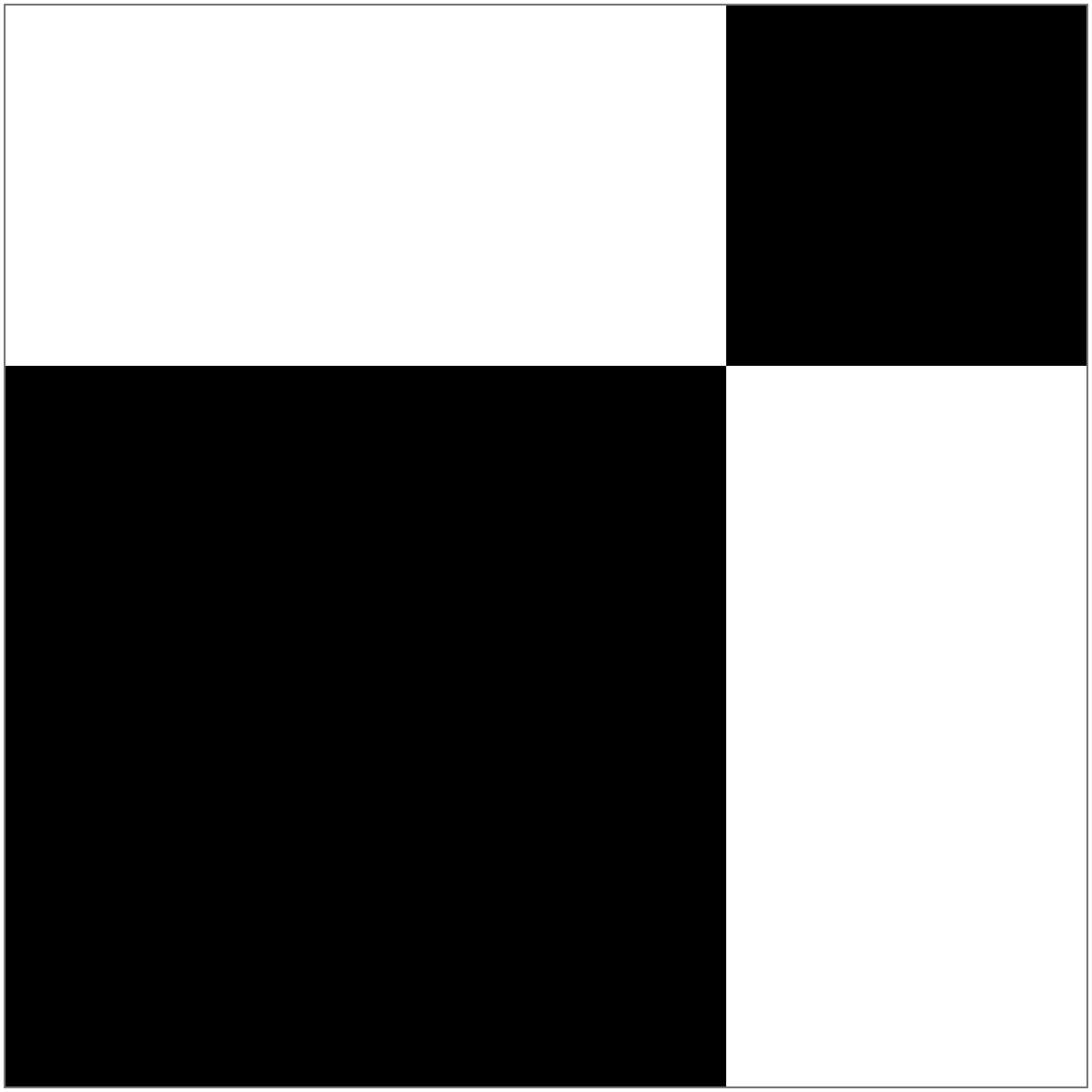}
&
\includegraphics[width=2.75cm]{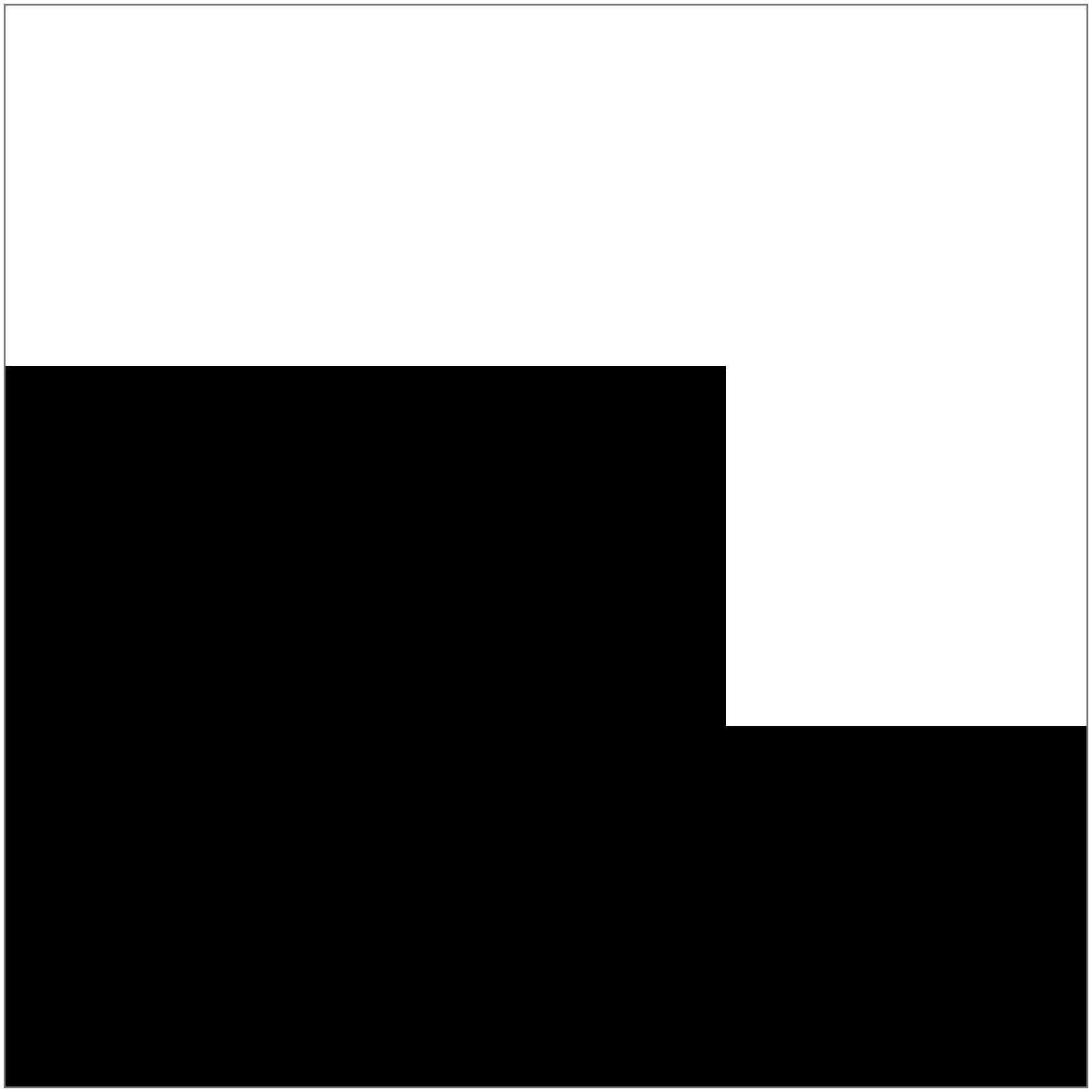}
&
\includegraphics[width=2.75cm]{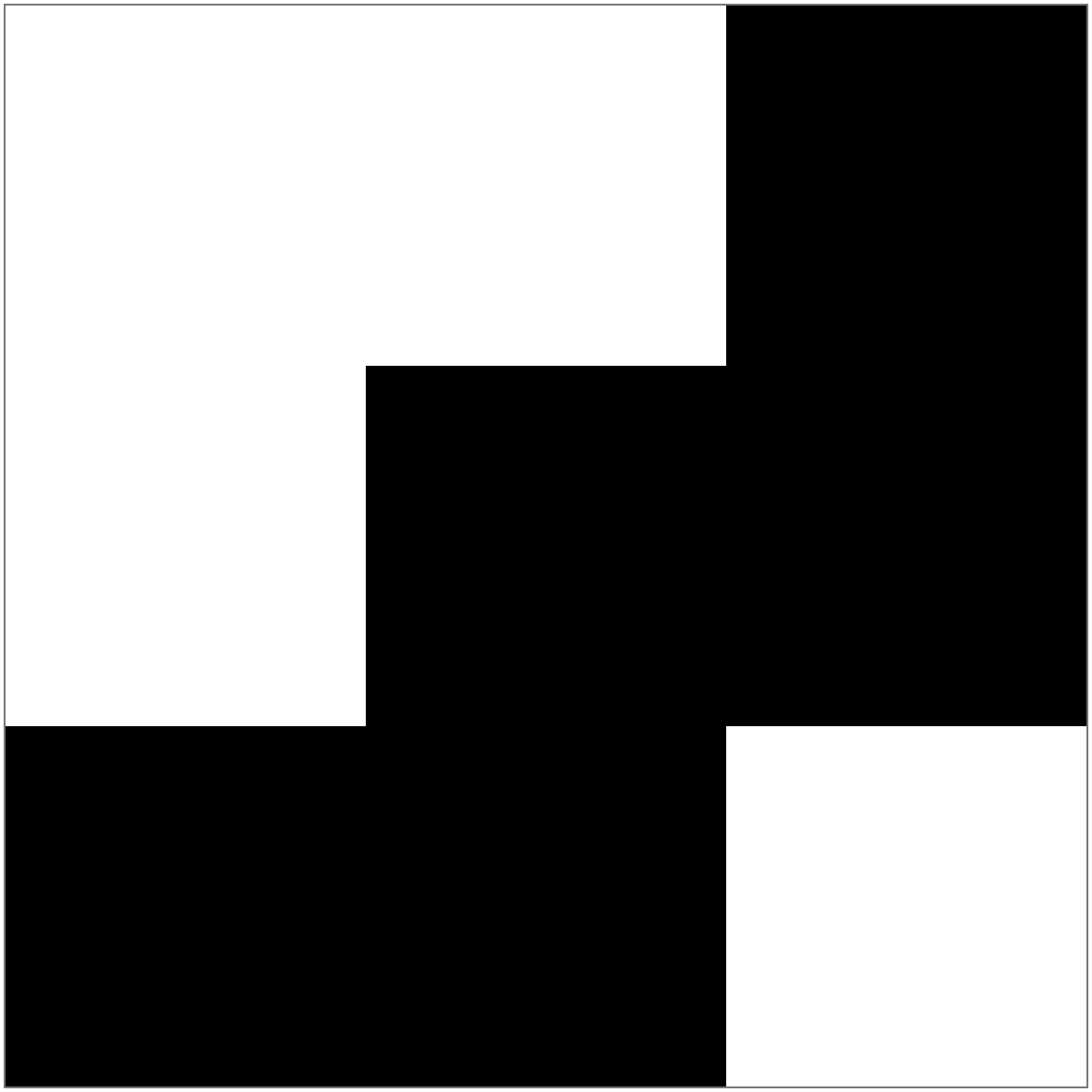}
&
\includegraphics[width=2.75cm]{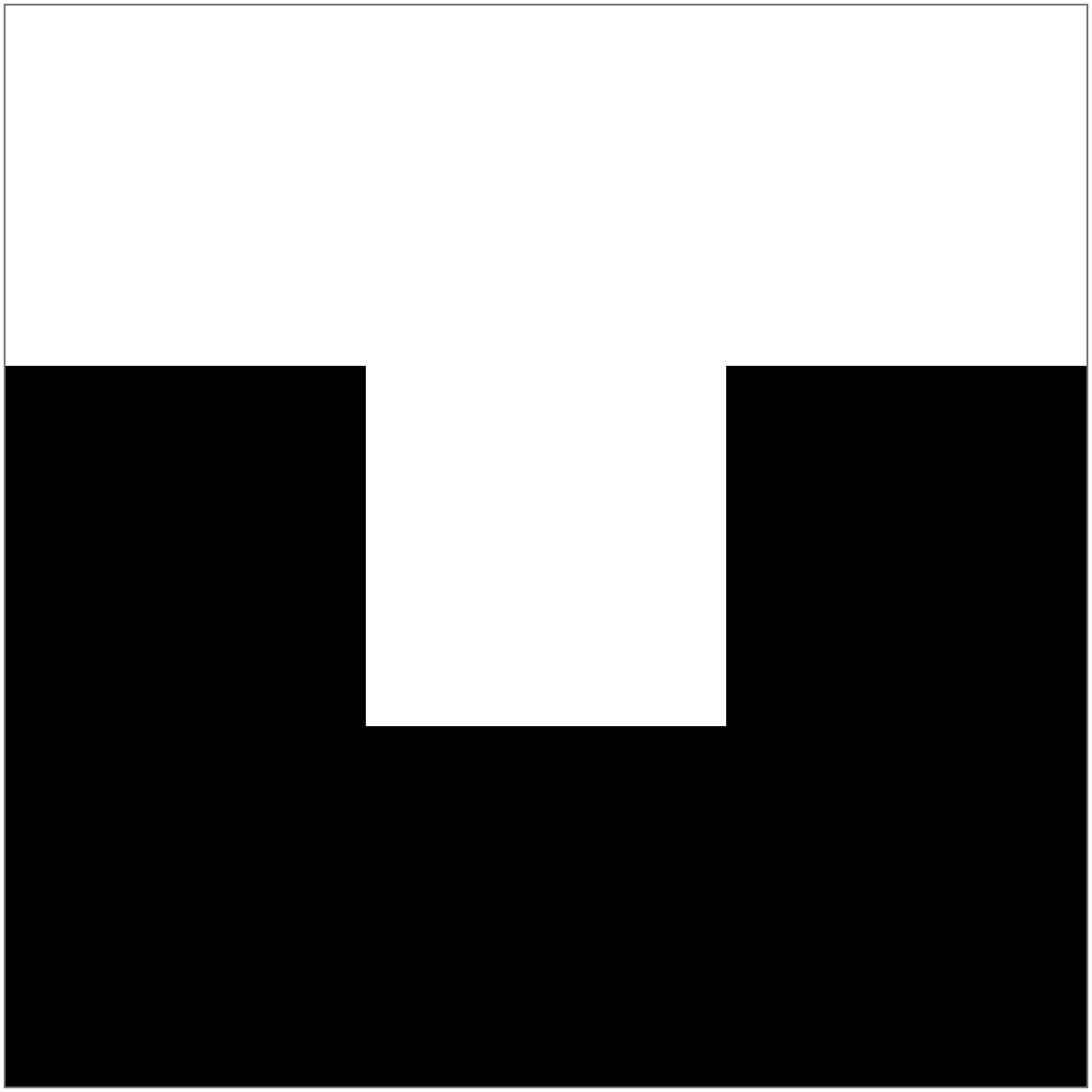}
\\
\includegraphics[width=2.75cm]{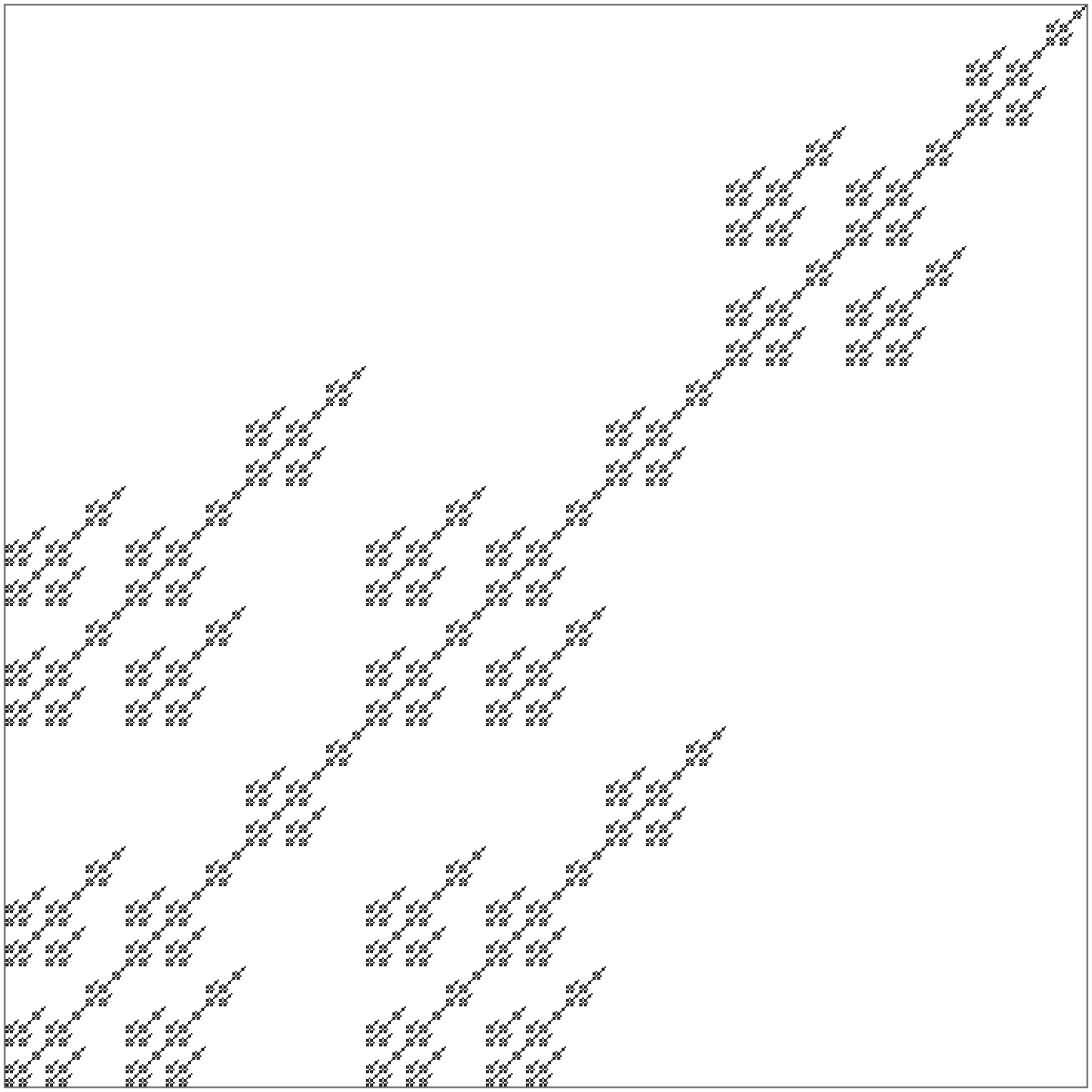}
&
\includegraphics[width=2.75cm]{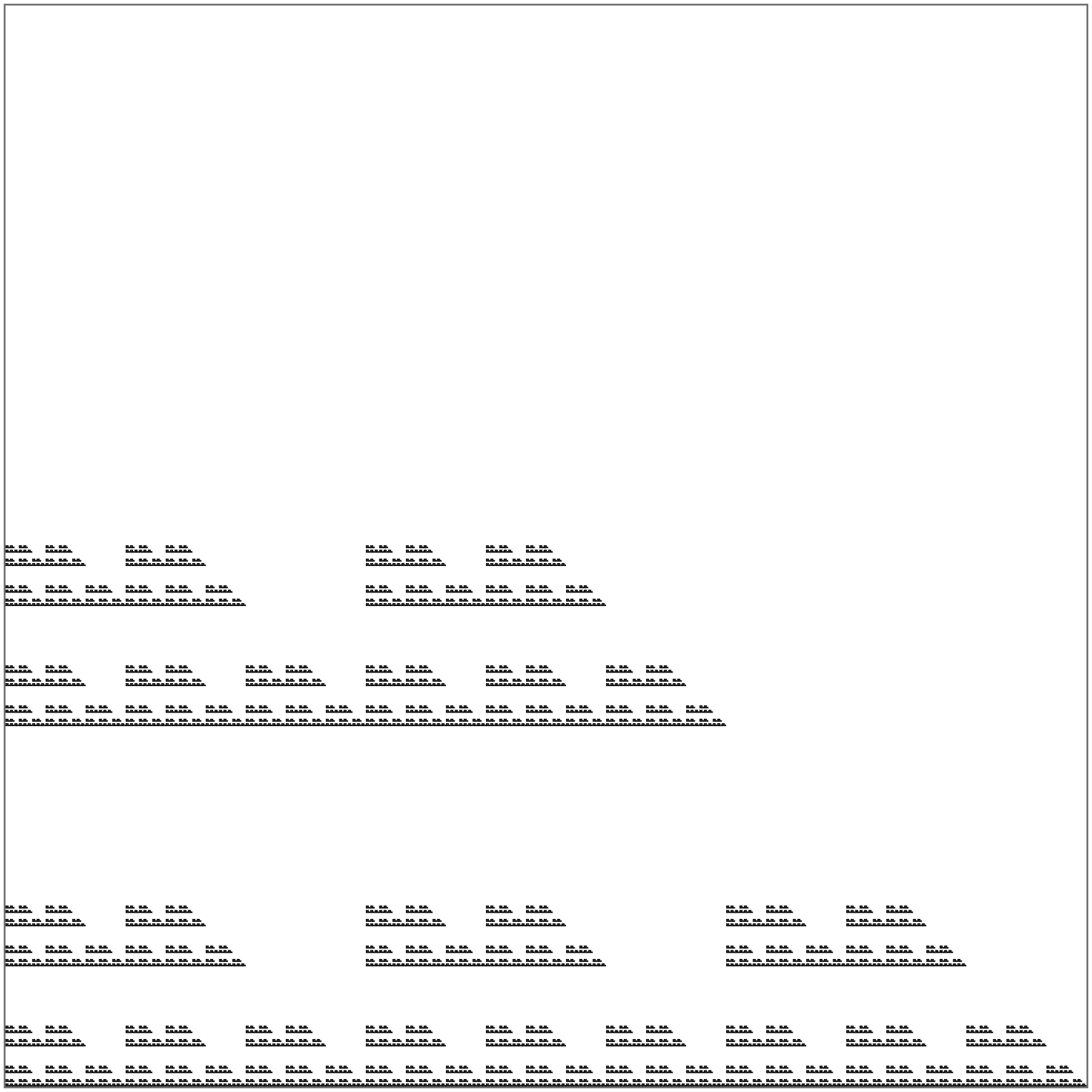}
&
\includegraphics[width=2.75cm]{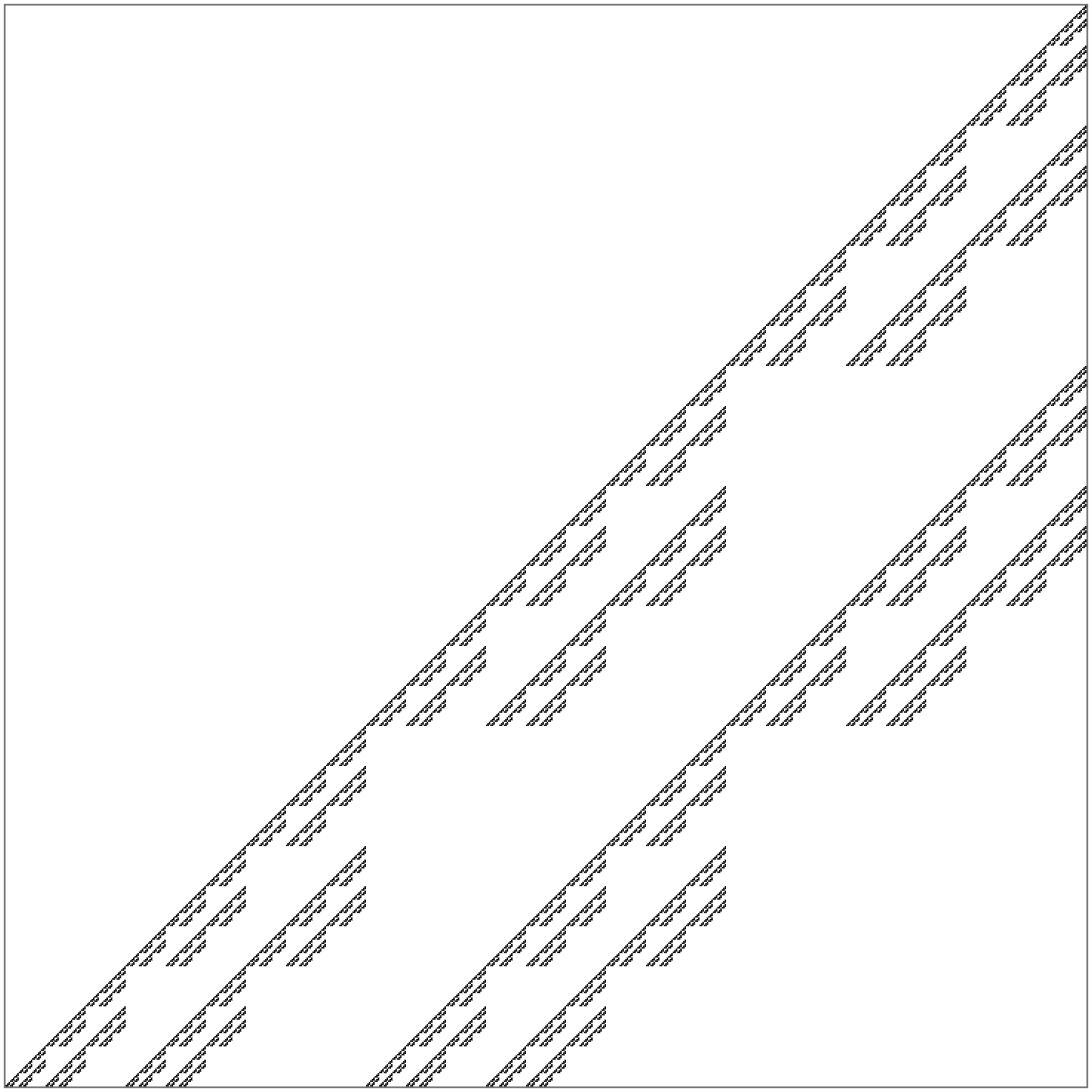}
&
\includegraphics[width=2.75cm]{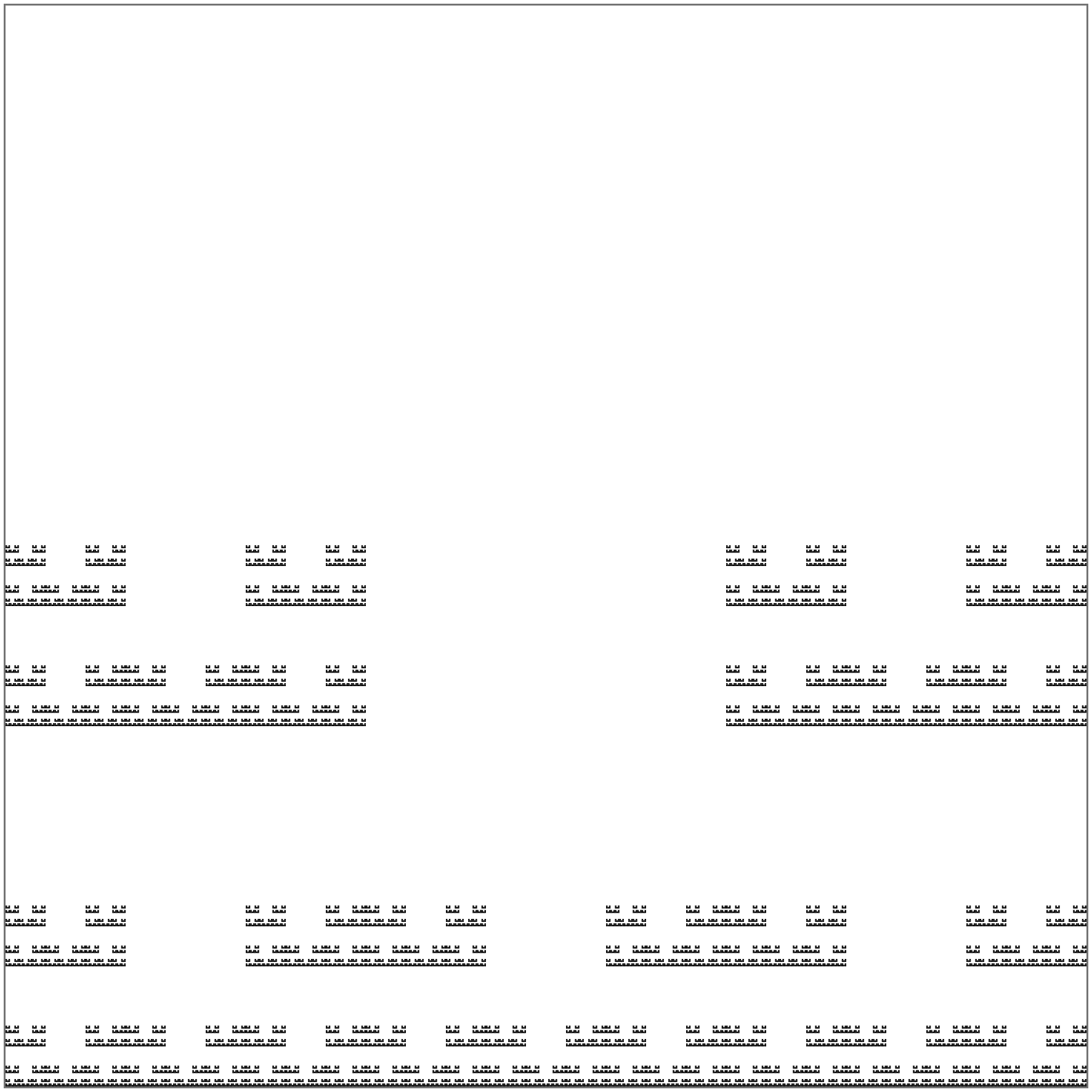}\\
$K_9^{(1)},\ K_9^{(6)}$ & $K_{10}^{(1)},\ K_{10}^{(6)}$ & $K_{11}^{(1)},\ K_{11}^{(6)}$ & $K_{12}^{(1)},\ K_{12}^{(6)}$
\end{tabular}\end{center}
\caption{The two approximations $K_j^{(n)}\ (n=1,6)$ for $K_j\ (9\le j\le 12)$.}\label{fig:order_3_more}
\end{figure}
\end{example}

Therefore, by Examples \ref{exmp:N=3_classes} to \ref{exmp:N=3_easy_homeo} and by results from \cite{ZhuYang18,ZhuRao21} we shall have.

\begin{theorem}\label{theo:topology_automaton}
Each fractal squares $K(3,\mathcal{D})\in\mathcal{K}_{3,5}$ that has both point components and line segment components is homeomorphic with some  $K_j$ with $4\le i\le7$. \end{theorem}

Consequently,  there are at most $8=1+3+4$ topological equivalence classes among all the fractal squares $K\in\mathcal{K}_{3,5}$.  In other words, we have $\chi(5)\le8$. To determine the exact value of $\chi(5)$, it suffices to resolve the following.
  
\begin{question}\label{que:q=5}
Is $K(3,\mathcal{D}_i)$ homeomorphic with $K(3,\mathcal{D}_j)$ for some $4\le i<j\le7$?
\end{question}

Fractal squares $K\in\mathcal{K}_{3,6}$ are carefully discussed in \cite{Rao-Wang-Wen17}. Every of those $K(3,\mathcal{D})$ is isometric with one of the sixteen given in \cite[Fig.1]{Rao-Wang-Wen17}. In other words, there are sixteen isometric classes. 
By \cite[Theorem 1]{Rao-Wang-Wen17}, six of those  $K(3,\mathcal{D})$ are special. Among them, two of them are the product of a Cantor set with $[0,1]$ hence are affine equivalent; two are affine equivalent but they only contain countably many line segment components; and the last two are Lipschitz equivalent. The rest ten satisfy two requirements. First, no two of them are topologically equivalent. Second, none of them is homeomorphic with any of the previous six  $K(3,\mathcal{D})$. Therefore, we have the following.

\begin{theorem}[{\cite[Theorem 1]{Rao-Wang-Wen17}}]\label{theo:q=5}
$\chi(6)=13$.
\end{theorem}

Fractal squares $K\in\mathcal{K}_{3,7}$ are thoroughly analyzed in \cite{Ruan-Wang17}, which concentrates on bi-Lipschitz equivalence between those fractals squares. By basic observations, it is easy to check that  every $K\in\mathcal{K}_{3,7}$ is isometric to one of the eight fractal squares given in \cite[Fig.4 to 11]{Ruan-Wang17}.  The complements in $\{0,1,2\}^2$ of their digit sets are respectively given below:
\begin{eqnarray}\label{eq:q=7}
\{(1,1),(1,2)\},  \quad   \{(1,0),(1,2)\},   & \{(1,2),(2,1)\},  \quad  \{(0,2),(2,2)\},  \\  
\{(0,1),(2,2)\},  \quad   \{(1,2),(2,2)\},   &  \{(0,0),(2,2)\},  \quad  \{(1,1),(2,2)\}.
\end{eqnarray}
For the sake of convenience, denote the resulting fractal squares by $K_j^*(1\le j\le 8)$.   The following five results are known, respectively  from \cite[Lemma 4.1 to 4.5]{Ruan-Wang17}.
\begin{itemize} 
\item[(1)] Each of $K^*_1$ and $K^*_2$ has a unique (global) cut point. 
\item[(2)] Each of $K^*_j(3\le j\le 5)$ has infinitely many (global) cut points. 
\item[(3)] $K^*_6$ has cuttings that contain exactly two points and there is no cut point. 
\item[(4)] $K^*_7$ has cuttings that contain exactly three points  and there is no cut point.  
\item[(5)] $K^*_8$ does not have finite cuttings.  
\end{itemize}
By \cite[Lemmas 4.7 and 4.8]{Ruan-Wang17}, we already know that no two of $K^*_j(3\le j\le 5)$ are homeomorphic. Therefore, $\chi(7)\ge7$ and the only issue is to determine whether $K^*_1$ and $K^*_2$ are homeomorphic or not. To this question, a definite answer is given in the following.

\begin{theorem}[{\cite[Theorem 2]{Luo-Yao21}}]\label{theo:Luo-Yao21}
Let $G_i(i=1,2)$ be the group of homeomorphisms of $K^*_i$ onto itself. Then $G_1$ has two elements while $G_2$ has eight.
\end{theorem}
 Therefore, we have.
 \begin{theorem}\label{theo:q=7}
 $\chi(7)=8$.
 \end{theorem}

Fractal squares $K(3,\mathcal{D})$ with $\#\mathcal{D}=8$ allow three possibilities. First, it is isometric with  the classical Sierpinski's carpet $K(3,\mathcal{D}_S)$. See Example \ref{exmp:sierpinski}. Second,  it has no local cut point hence is homeomorphic with $K(3,\mathcal{D}_S)$. See for instance the proof for \cite[Theorem 1.2]{Ruan-Wang17}. Third, it  is isometric to the fractal square given in Example \ref{exmp:N=3_local_cut}. 
  
 \begin{example}\label{exmp:N=3_local_cut}
 Let  $K=K(3,\mathcal{D}_L)$   be given with   $\mathcal{D}_L=\{0,1,2\}^2\setminus\{(2,2)\}$. See Figure \ref{fig:local_cut}
 for a depiction of $K(3,\mathcal{D}_L)$ and of its first two approximation. It is routine to verify that $K(3,\mathcal{D}_L)$ has infinitely many local cut points, one of which is $x_0=(2,2)$. 
 \begin{figure}[ht]    
\begin{center}
\begin{tabular}{cc}
\includegraphics[width=3cm]{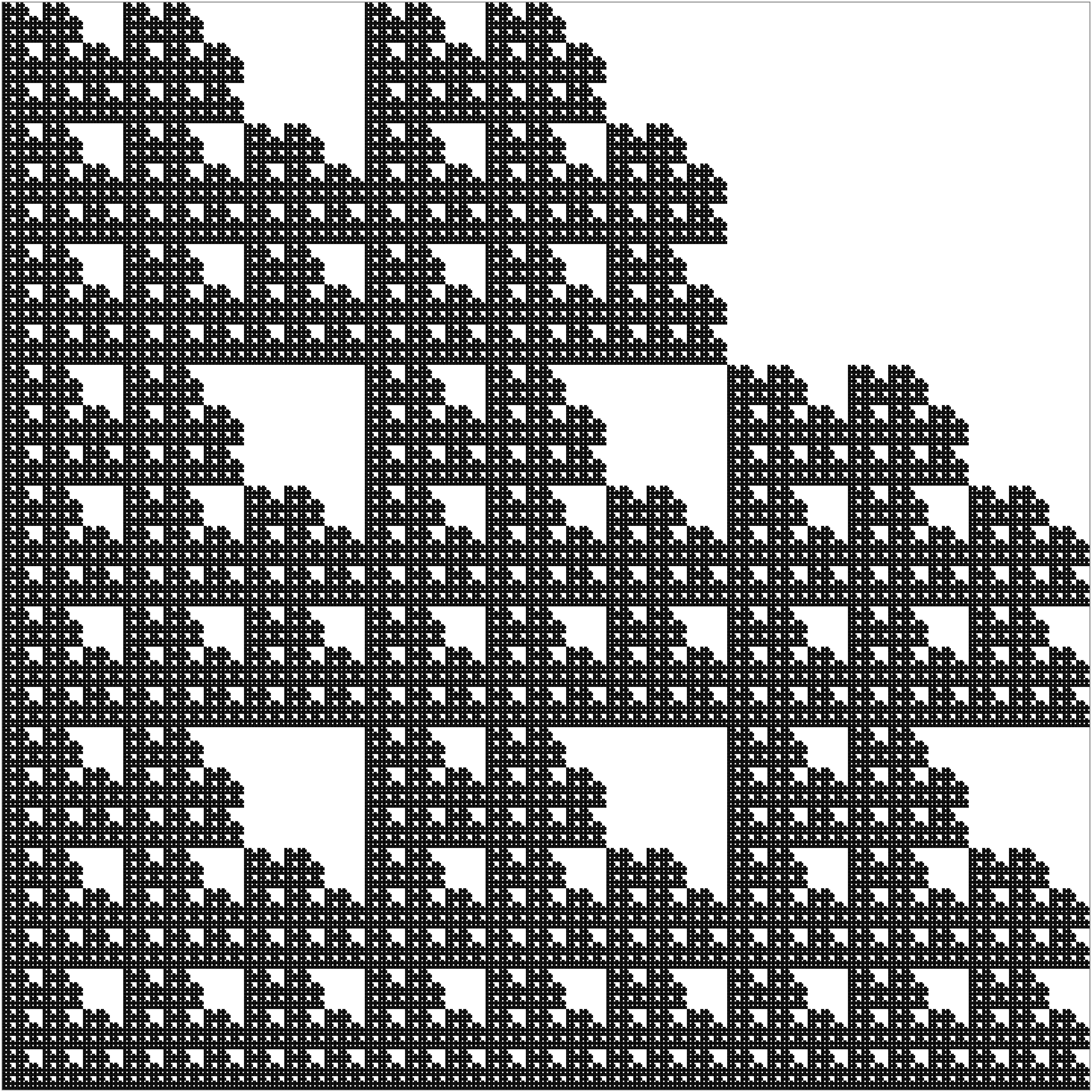}
&
\includegraphics[width=3cm]{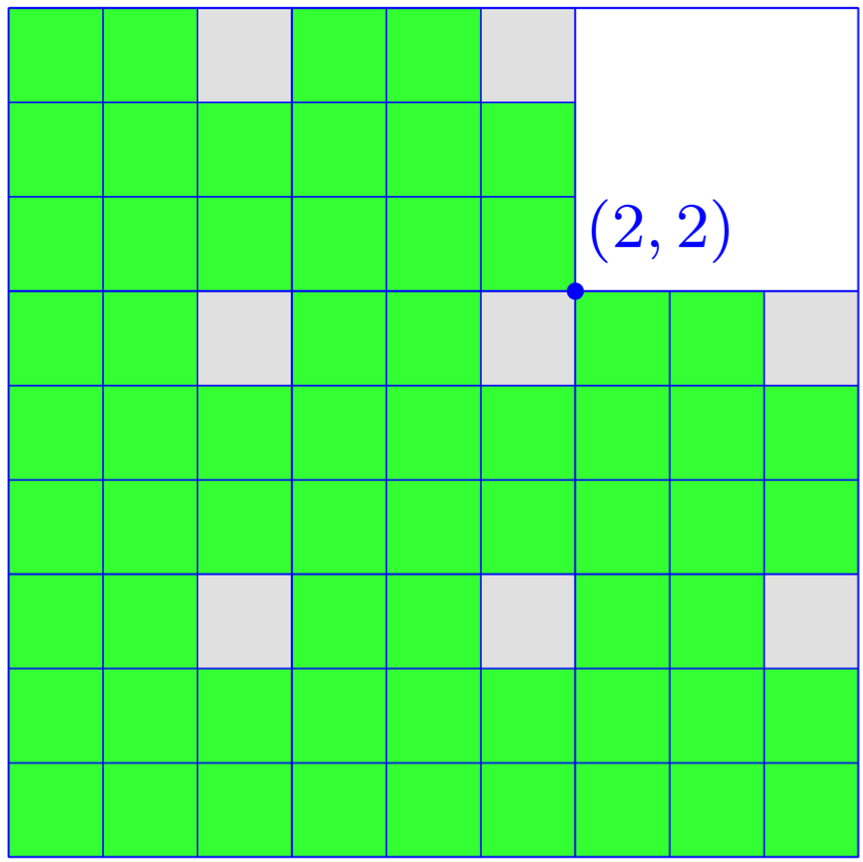}
\end{tabular}\end{center}
\vspace{-0.618cm}
\caption{The fractal square $K(3,\mathcal{D}_L)$ and its first two approximations.}\label{fig:local_cut}
\end{figure}
\end{example}
 
The following is immediate.
 \begin{theorem}\label{theo:q=8}
 $\chi(8)=2$.
 \end{theorem}

\section{Note on Fundamental Groups}\label{sec:pi_1}

Fractal squares $K(N,\mathcal{D})$, except for the case $\mathcal{D}=\{0,\ldots,N-1\}^2$, have no interior point. 
By Theorem \ref{theo:LRX}, 
every component $P$ of $K$  is a Peano continuum. Moreover, the topological dimension of $P$ is one. This section concerns the fundamental group of connected fractal squares or of the components of a disconnected fractal square. We shall need the following terminology.

\begin{definition}\label{def:semi_lsc}
A Peano continuum $X$ is said to be {\em semi-locally simply connected} at  $x_0\in X$ provided that for each neighborhood $U$ of $x_0$ there is a neighborhood $V$ of $x_0$, with $V\subset U$, such that every loop in $V$ based at $x_0$  can be contracted in the original space $X$. If $X$ is semi-locally simply connected at none of its points, we also say that it is nowhere semi-locally simply connected. \end{definition}

If $K$ is actually connected then it is locally path-connected; if further it has a nontrivial fundamental group then it is semi-locally simply connected at none of its points.  The following result sheds some light on how the topology of such a continuum is related to its fundamental group.
\begin{theorem}[{\cite[Theorem 1.3]{Eda02}}]\label{theo:Eda02}
Let $X,Y$ be one-dimensional Peano continua  that are nowhere semi-locally simply connected. Then for $X$, $Y$ to be homeomorphic it is necessary and sufficient that the fundamental groups $\pi_1(X)$, $\pi_1(Y)$ are isomorphic.
\end{theorem}

More precisely, one may recover the topology of such a Peano continuum from its fundamental group. See for instance \cite{ConnerEda05}.
Thus we are motivated  to propose the following.

\begin{question}\label{algo:pi1}
Given a connected fractal square $K$, find conditions on the initial approximations $K^{(j)}$, say for $j\le4$, for $\pi_1(K)$ to be nontrivial. The same question for the components of $K$, when $K$ itself is NOT connected.
\end{question}   

An important motivation of Question \ref{algo:pi1} comes from the following fundamental connection.
\begin{theorem}[{\bf \cite[Theorem 2.8]{CKLY-2025}}]\label{theo:CKLY2.8}
Let $K$ be a connected fractal square of order $N\ge2$ satisfying $\mathcal{D}\ne\{0,1,\ldots,N-1\}^2$.
If  the fundamental group $\pi_1\left(K^{(3)}\right)$ is trivial so is $\pi_1(K)$.
\end{theorem}

Question \ref{algo:pi1} is also motivated by simple examples of four connected fractal squares $K$ that illustrate the possibilities of the fundamental group 
$\pi_1\left(K^{(n)}\right)$ for $n=1,2,4,6$.  One of them is mentioned in \cite[Example 2.9]{CKLY-2025}, the other three are variants that have a slightly larger digit set. 

\begin{example}\label{exmp:pi1_possibilities}
Let $\mathcal{D}_j\ (1\le j\le3)$ be three subsets of $\{0,1,2,3\}^2$ and $K_j$ the resulting fractal squares. Here $K_2$ is also given in \cite[Example 2.9]{CKLY-2025}, while $K_1$ and $K_3$ are obtained by removing or adding  digits.   
\begin{figure}[ht]    
\begin{center}
\begin{tabular}{cccc}
\includegraphics[width=2.75cm]{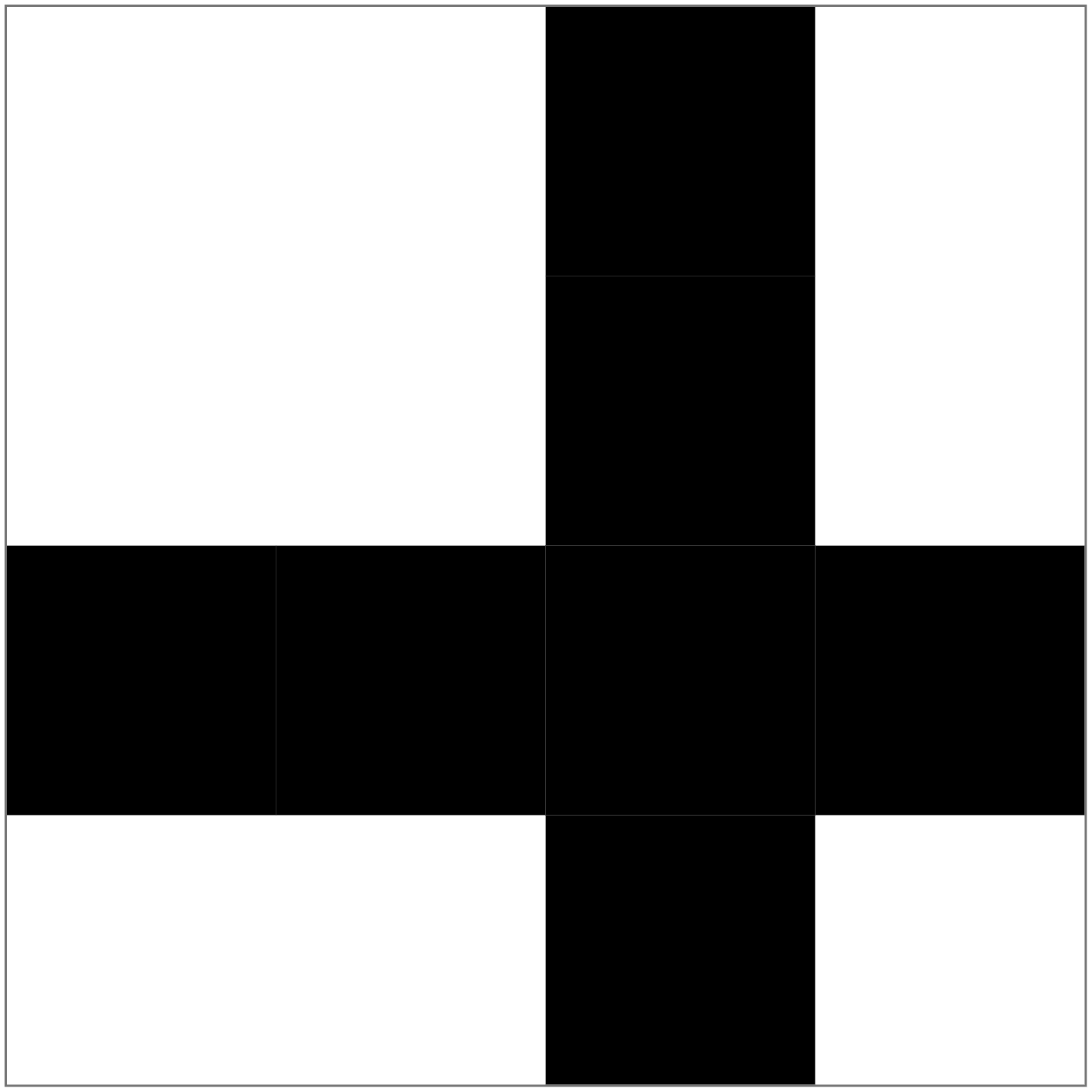}&
\includegraphics[width=2.75cm]{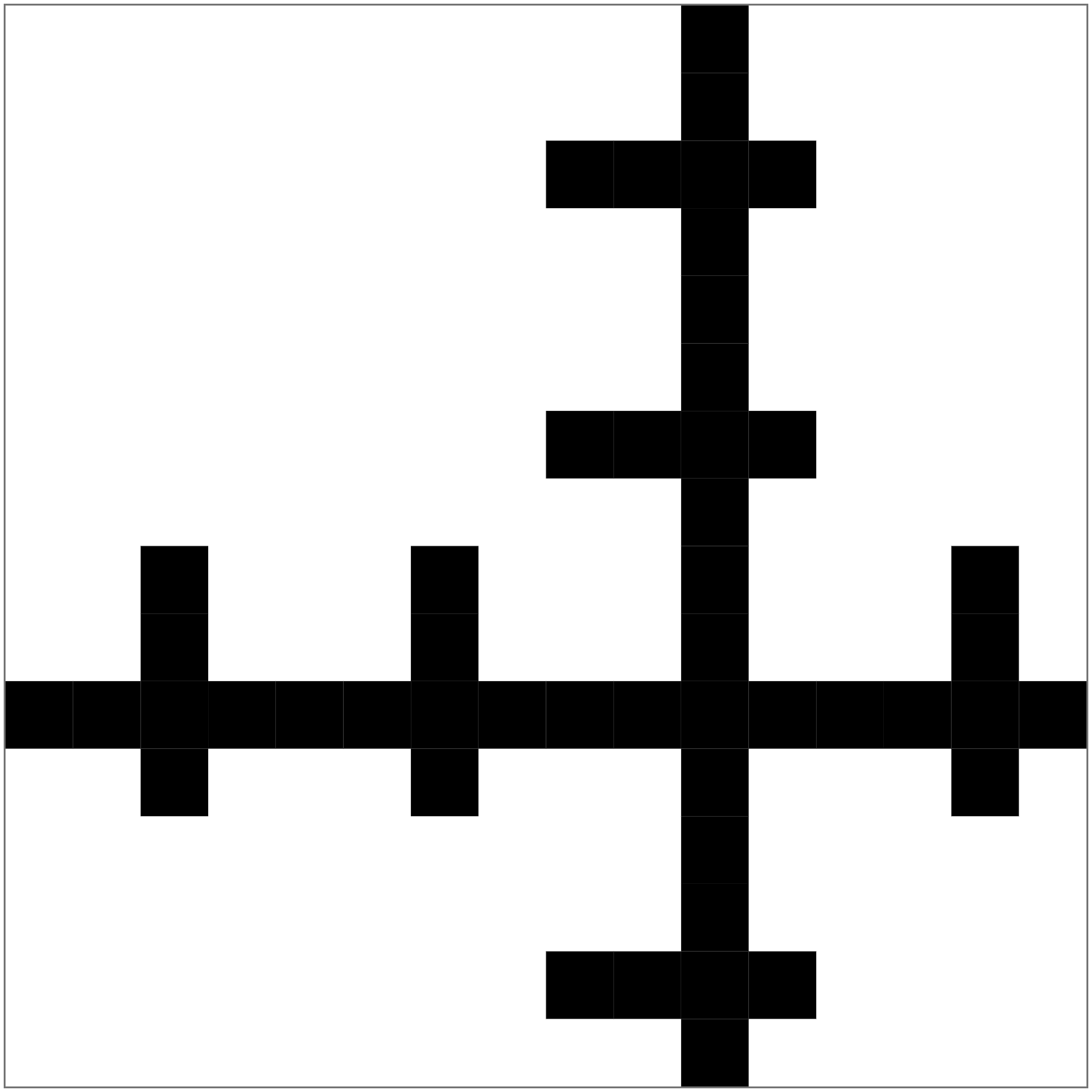}&
\includegraphics[width=2.75cm]{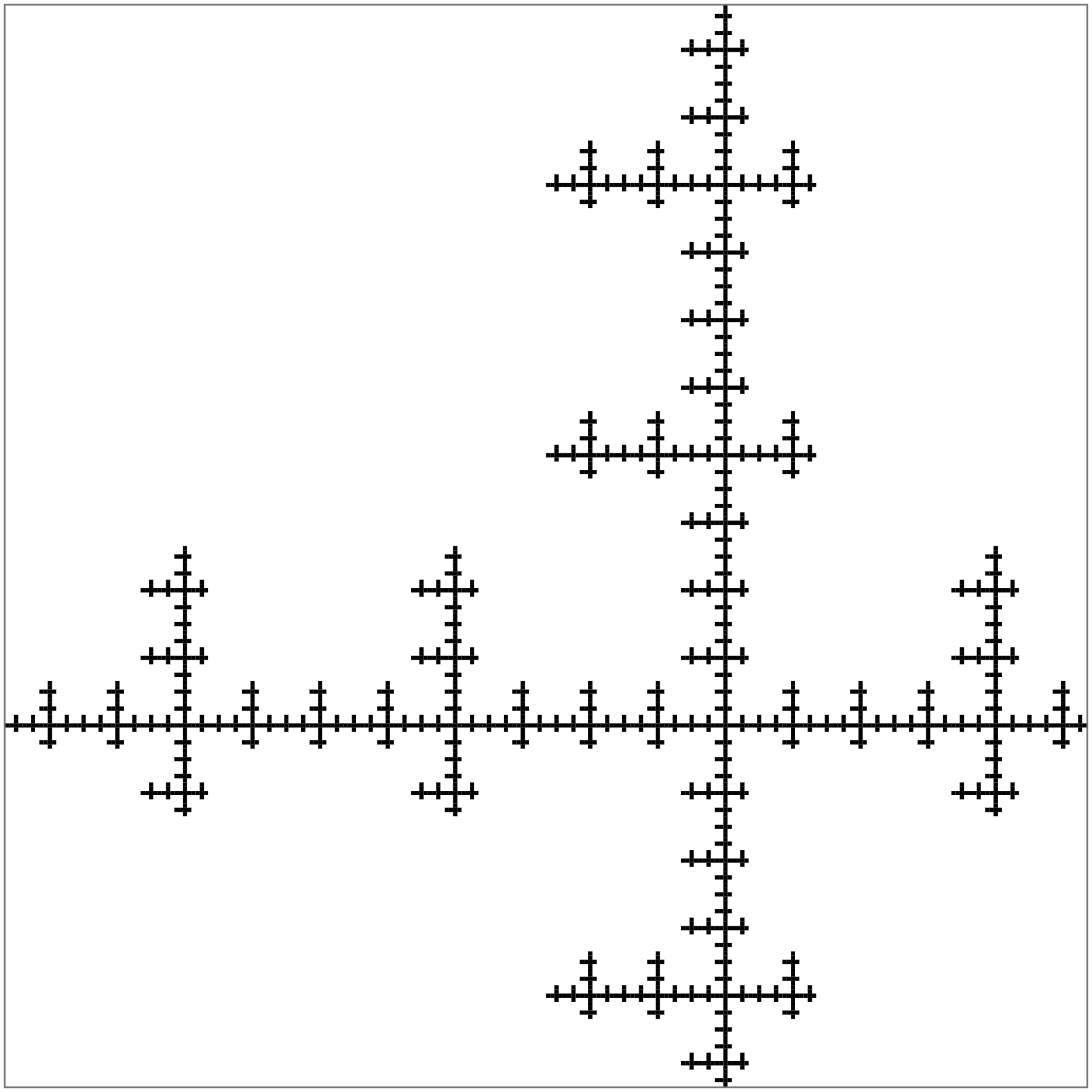}&
\includegraphics[width=2.75cm]{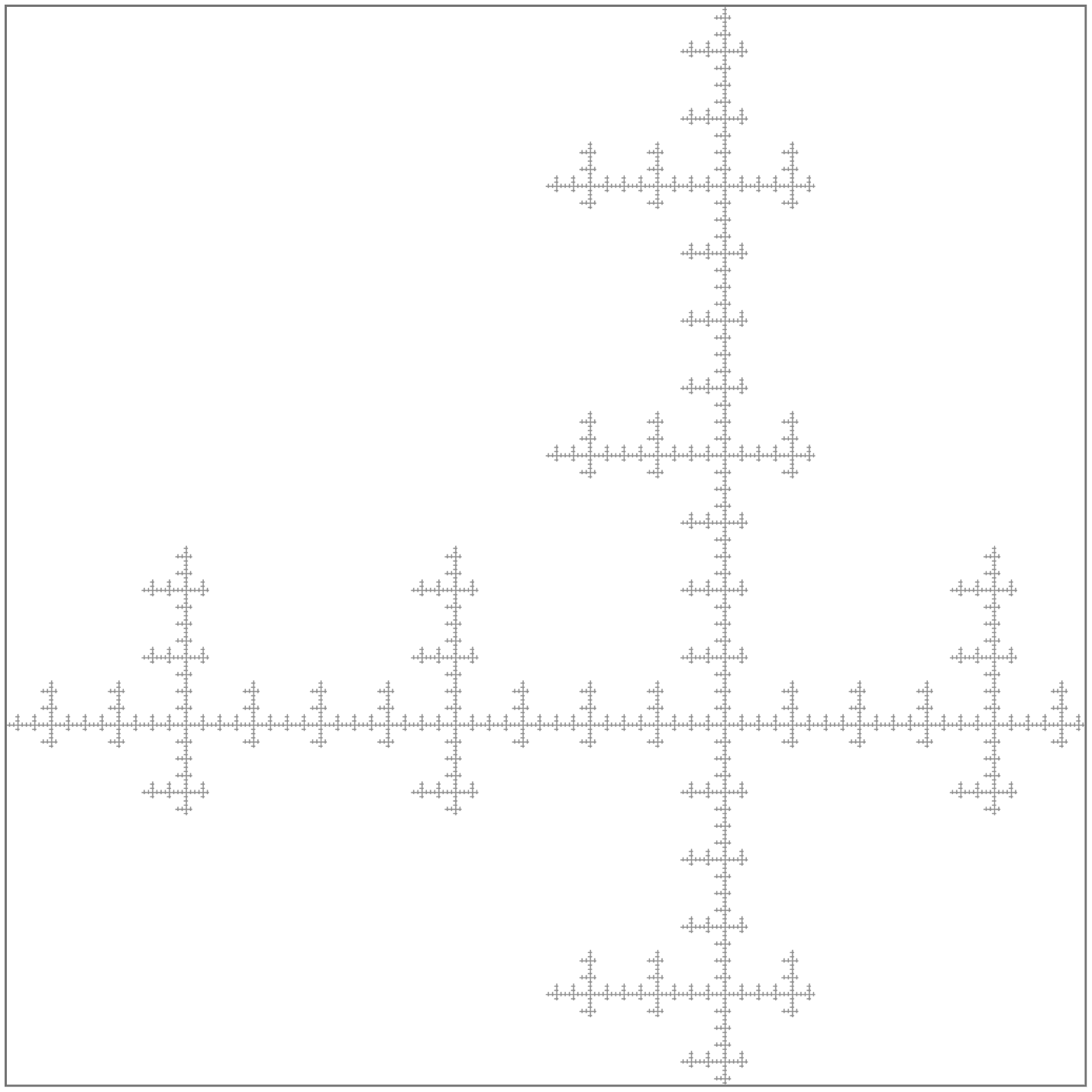}
\\
\includegraphics[width=2.75cm]{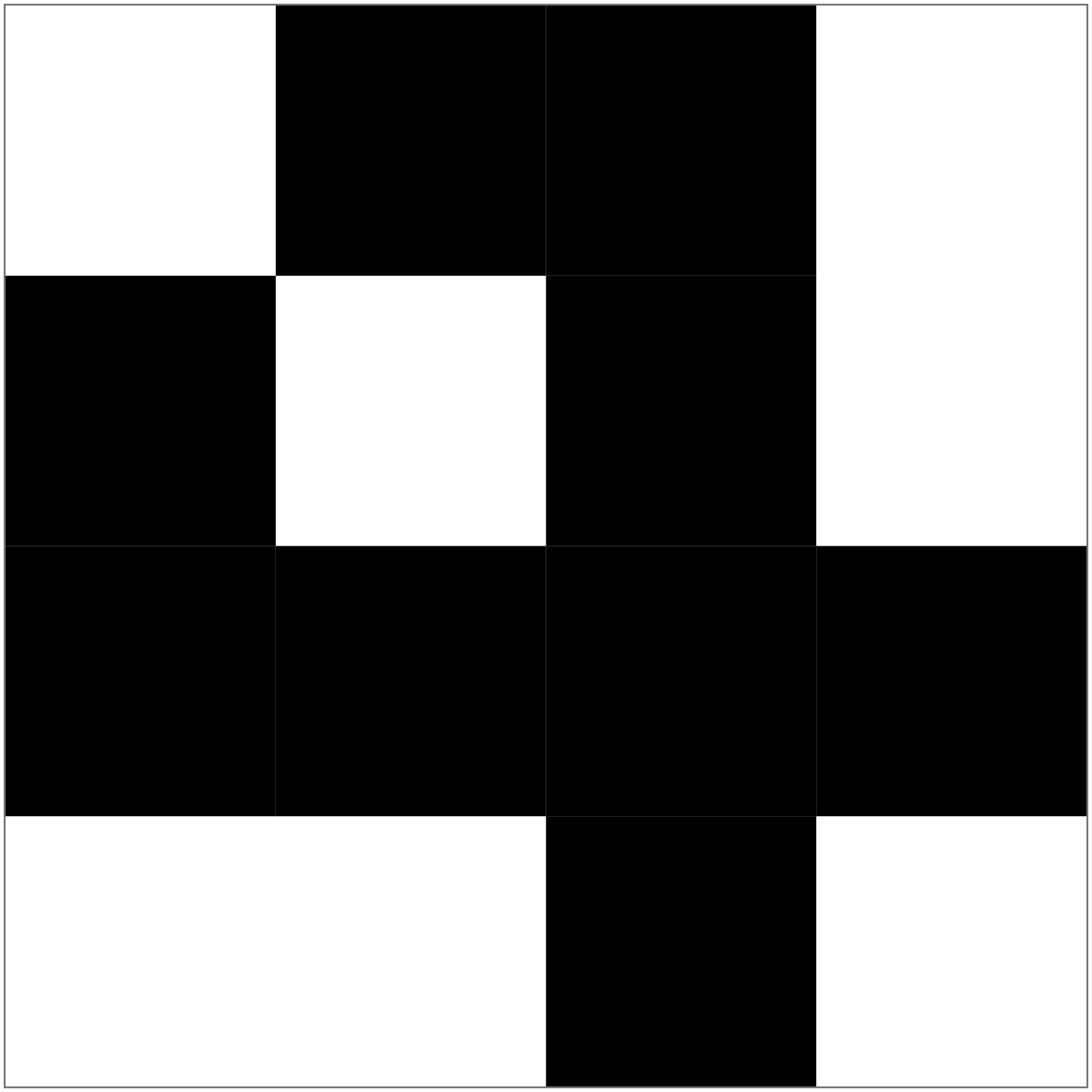}&
\includegraphics[width=2.75cm]{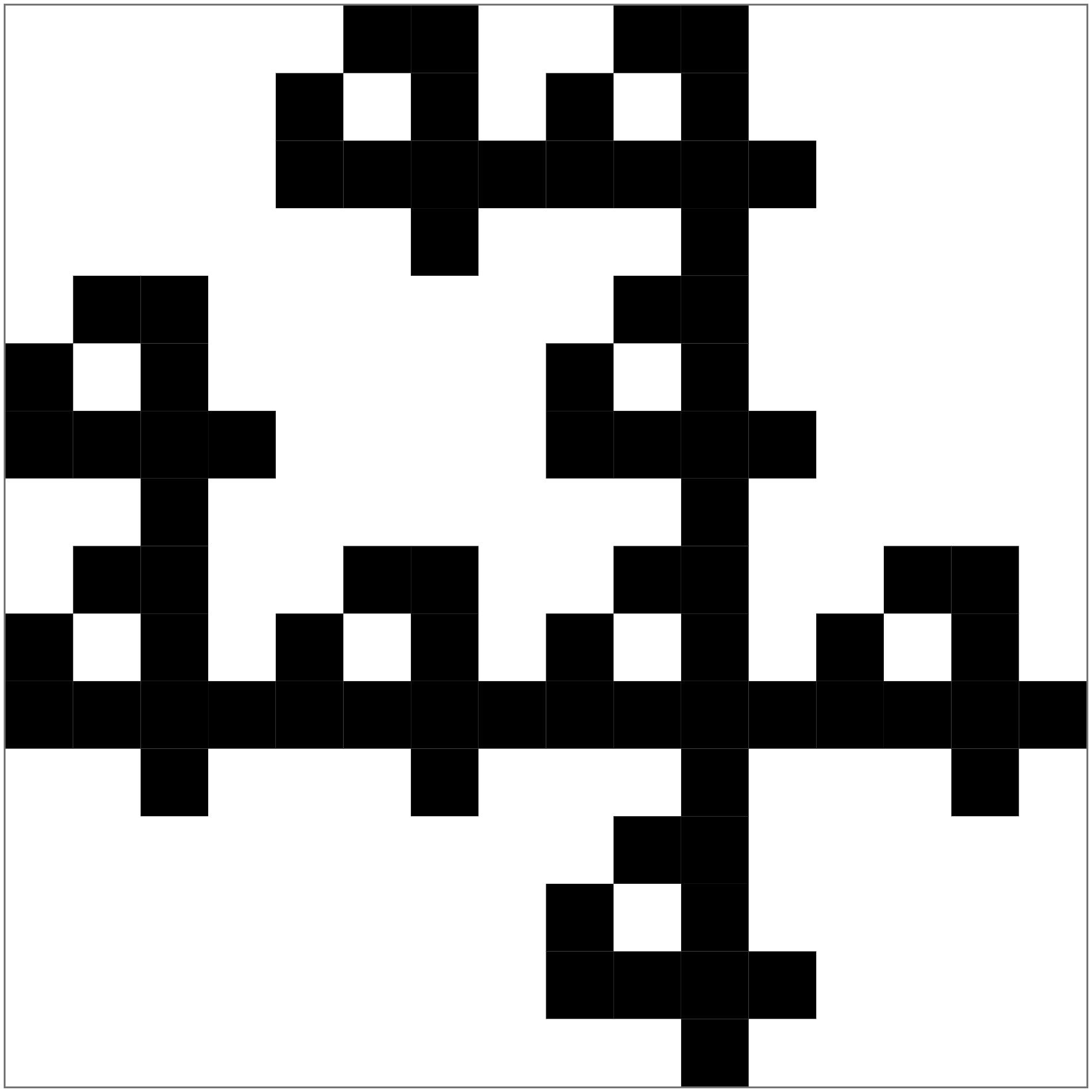}&
\includegraphics[width=2.75cm]{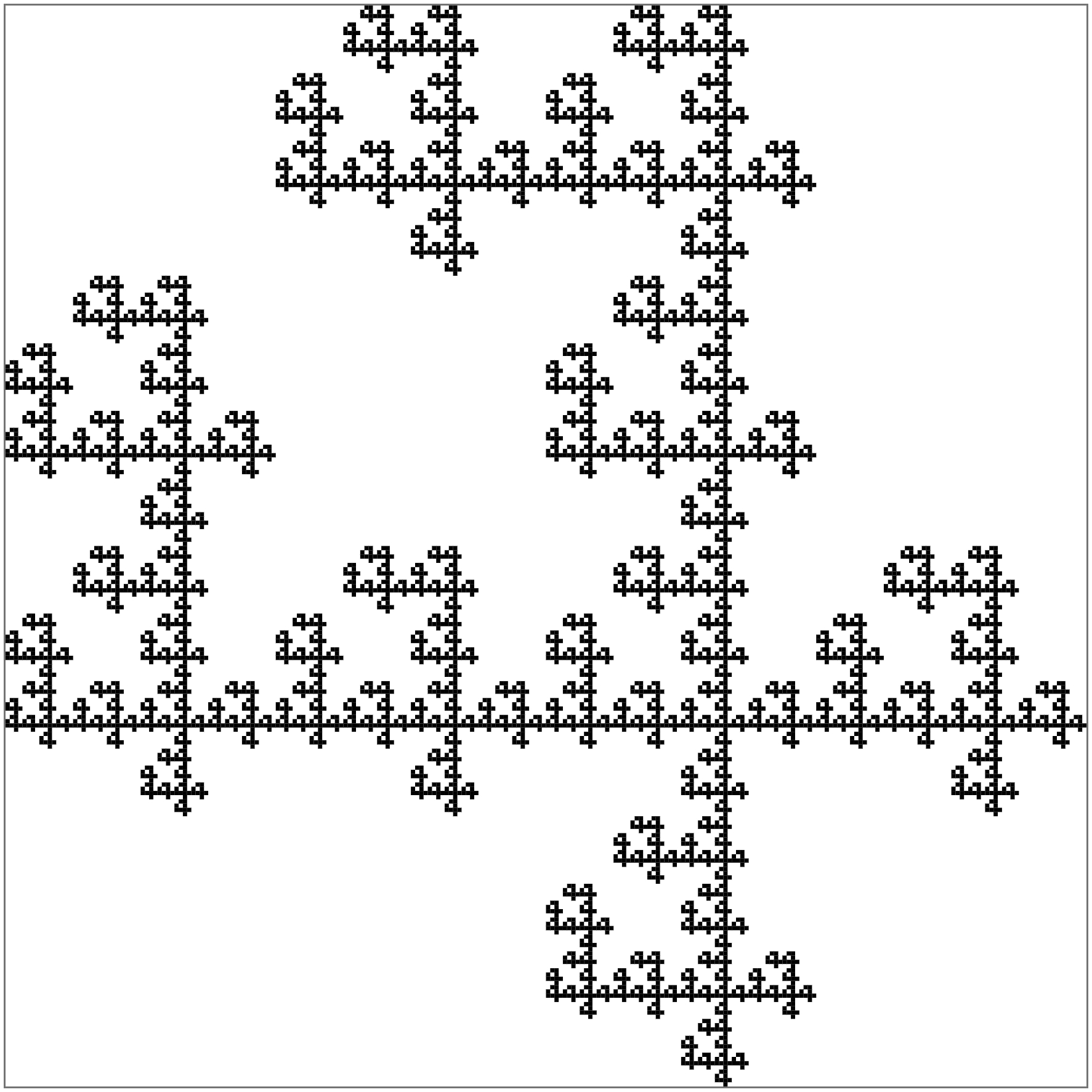}&
\includegraphics[width=2.75cm]{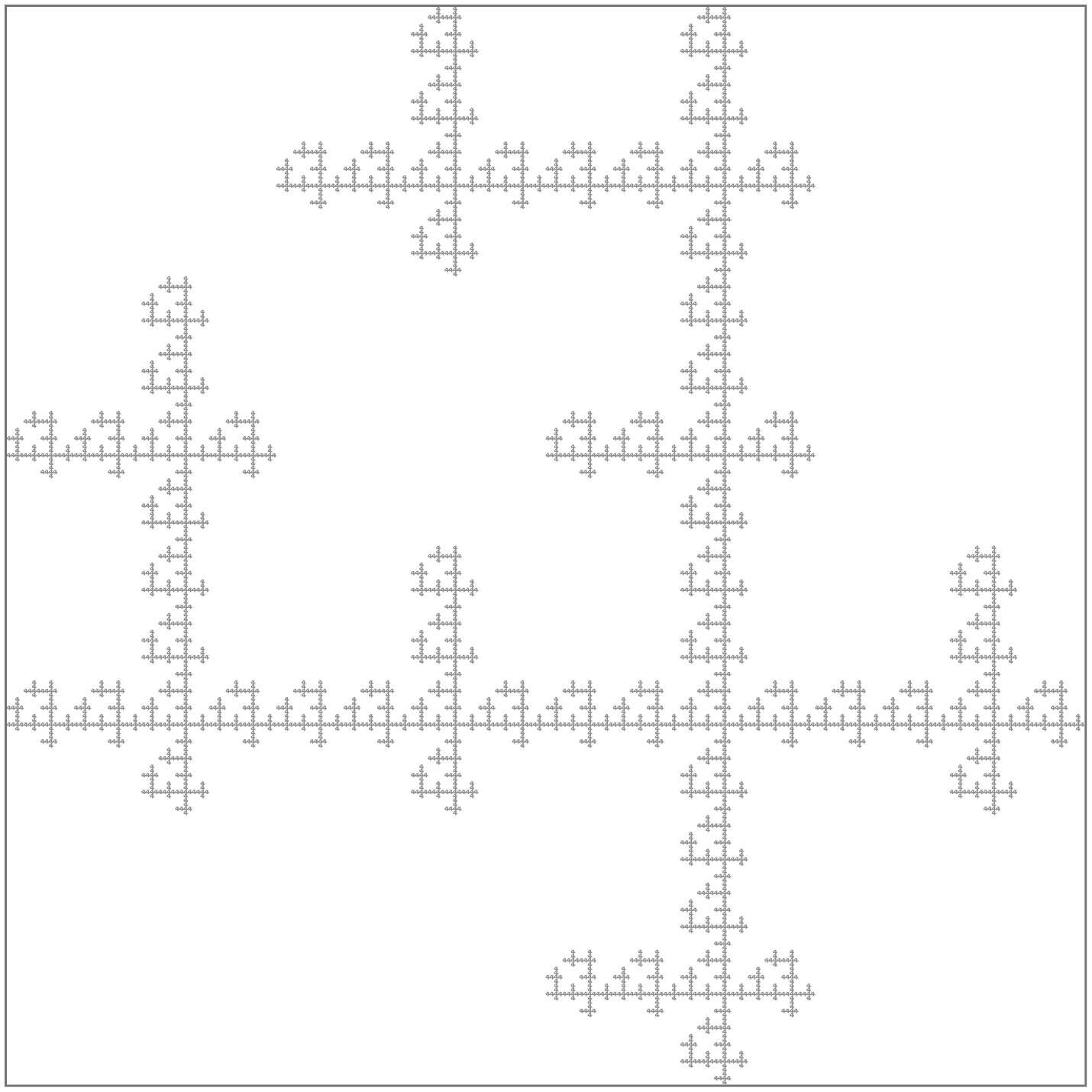}\\
\includegraphics[width=2.75cm]{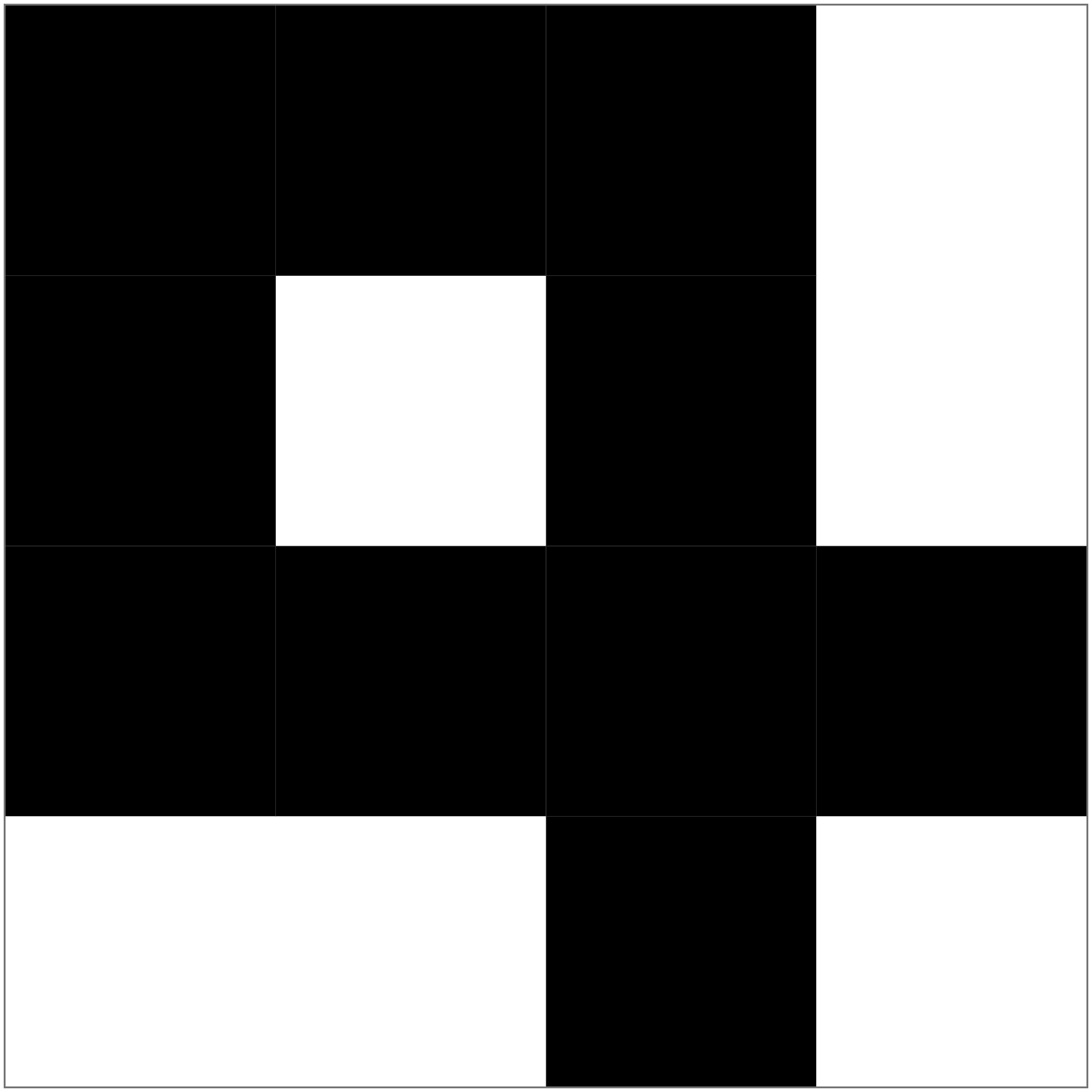}&
\includegraphics[width=2.75cm]{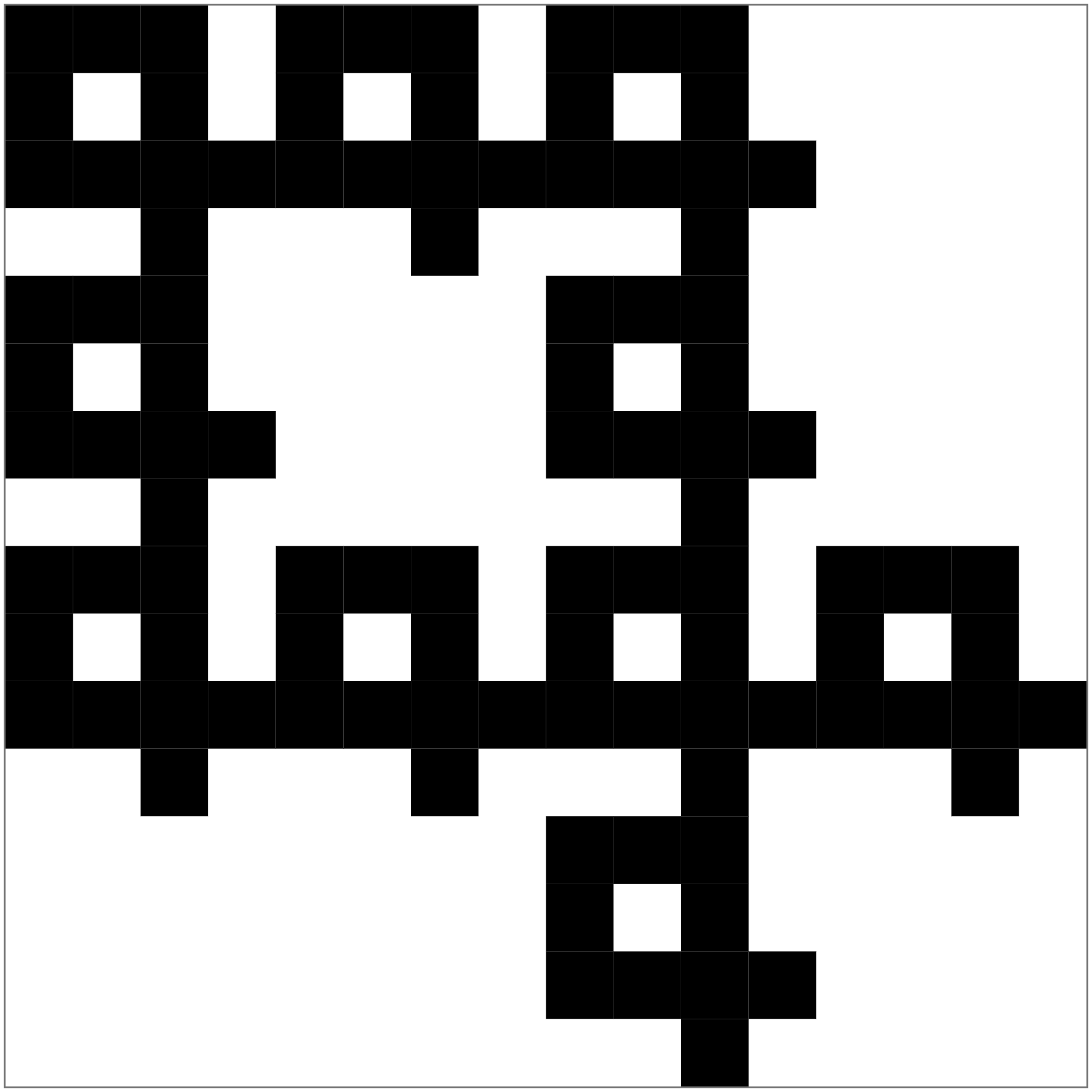}&
\includegraphics[width=2.75cm]{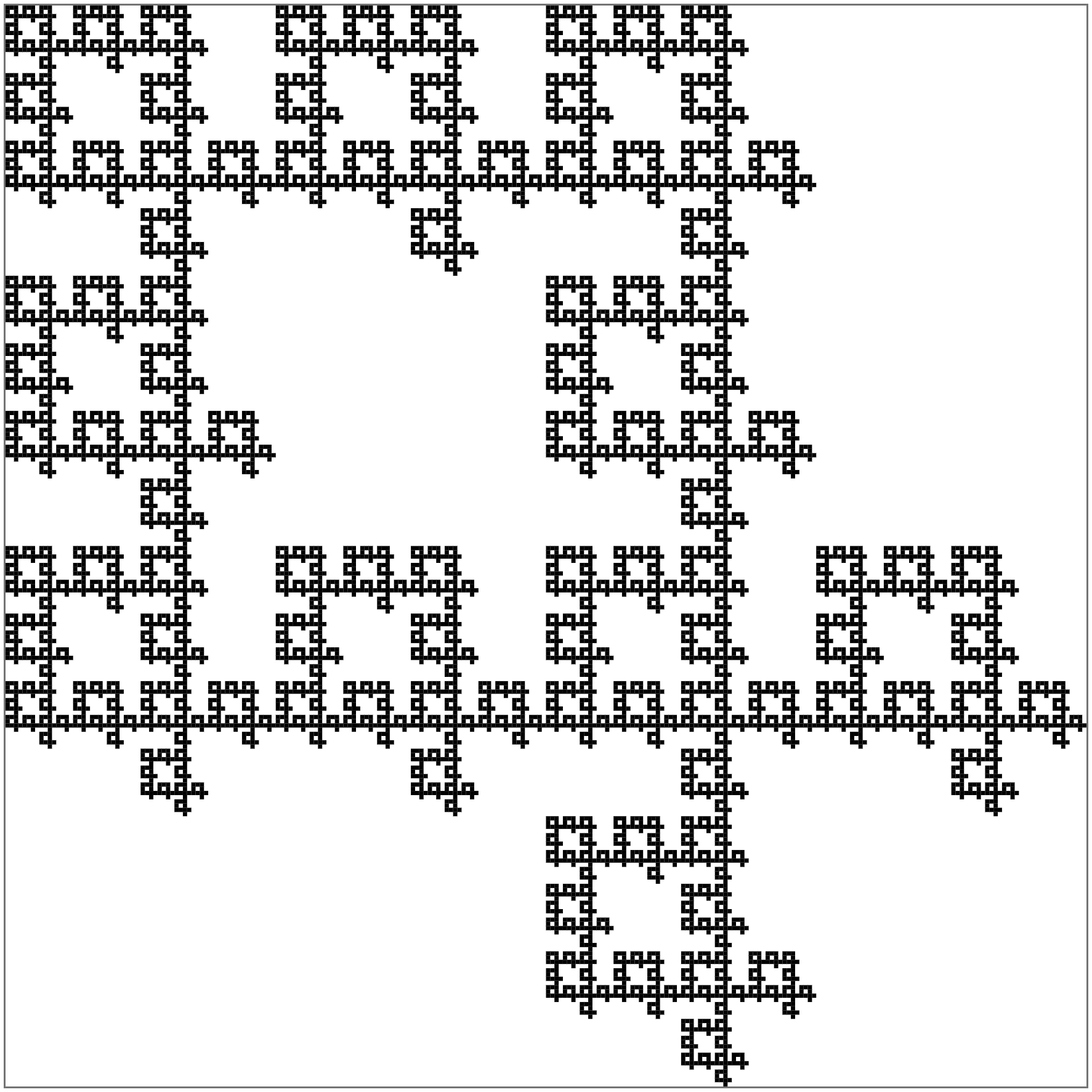}&
\includegraphics[width=2.75cm]{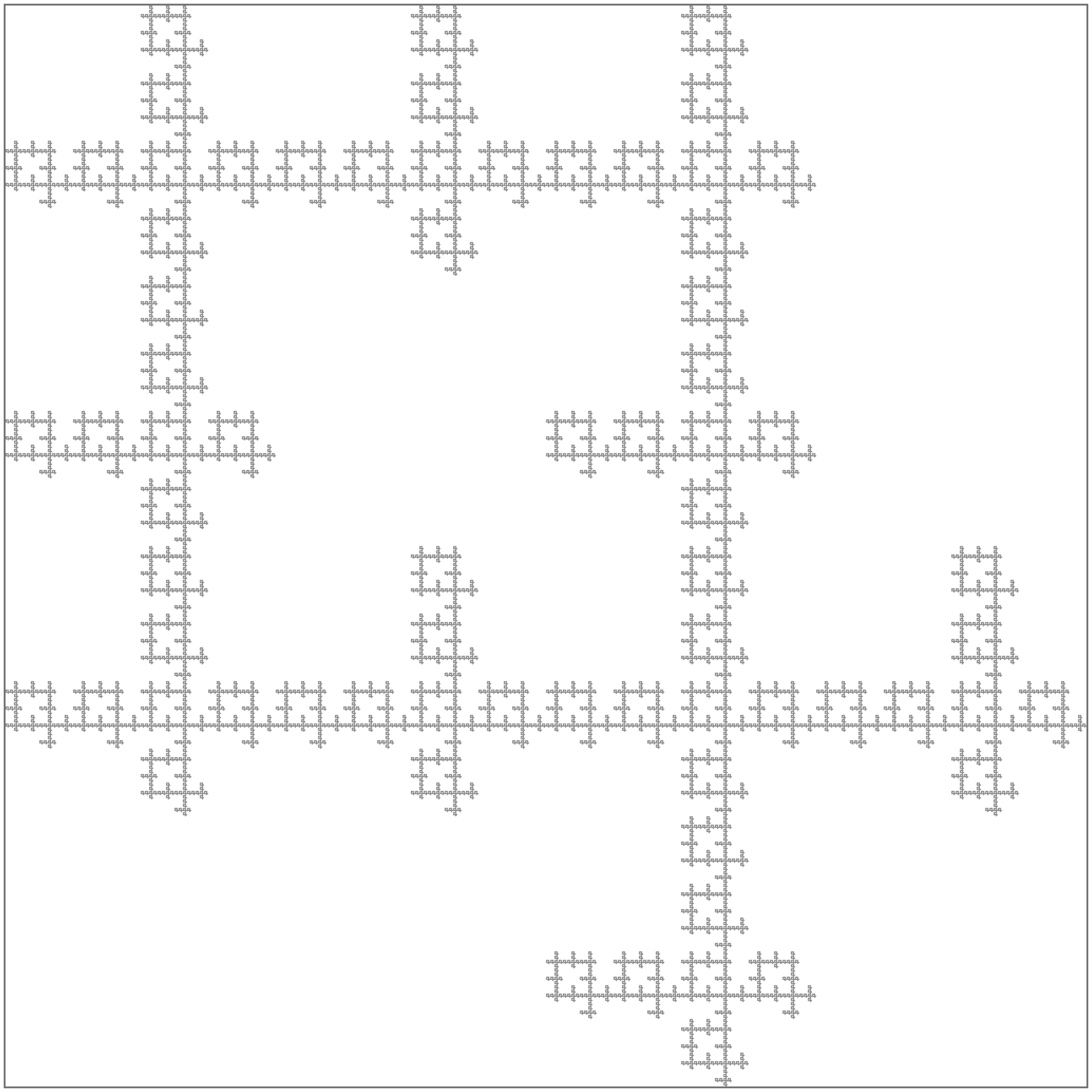}
\end{tabular}
\end{center}
\caption{The four approximations $K^{(n)}_j\ (n=1,2,4,6)$ for $K_j\ (j=1,2,3)$.}
\end{figure}
The next observations are routine. First, $\pi_1\left(K^{(n)}_1\right)=\{0\}$ for all $n\ge1$ hence by Theorem \ref{theo:CKLY2.8} we shall have $\pi_1(K_1)=\{0\}$. Second,  $\pi_2\left(K^{(n)}_1\right)\ne\{0\}$ for all $n\ge1$ while $\pi_1(K_2)=\{0\}$. This implies that the converse of Theorem \ref{theo:CKLY2.8} is not true.  Third, $\pi_1\left(K^{(n)}_3\right)\ne\{0\}$ for all $n\ge1$ and $\pi_1(K_3)\ne\{0\}$.
\end{example}

There are many other examples of fractal squares  $K$, such that some of its components  or $K$ itself are dendrites (possibly line segments). So, one may further wonder whether such a dendrite is actually a Jordan arc. Therefore, we also propose the following.

\begin{question}\label{algo:component}
Fix  a  fractal square $K$ and one of its non-degenerate components, say $P$, such that $\pi_1(P)$ is trivial. Find conditions on the initial approximations $K^{(j)}$, say for $j\le4$, for $P$ to be a Jordan arc.
\end{question}   

From the literature one can find quite a lot of efforts of characterizing  special or even general self-similar sets $K$ that are homeomorphic with $[0,1]$. See for instance \cite[Theorem 6.7]{Hata85}. See also \cite{BandtKeller91} for useful description of the topology of a self-similar set (resulted from an IFS $\mathcal{F}=\{f_1,\ldots,f_q\}$)  based on the symbolic  representation $\mathcal{S}:\{1,\ldots,q\}^\infty\rightarrow K$,  and \cite{Tetenov06} for graph-directed IFS's that generate  attractors each of which is a Jordan arc. 

\vspace{0.618cm} 

\noindent 
{\bf Acknowledgement}  This study was funded by the National Key R\&D Program of China [No. 2024YFA 1013700]. All pictures are in PNG form. The two colored ones, in Figures \ref{fig:cut} and  \ref{fig:local_cut}, are drawn in tikz packages. The other pictures are made by a new program, which is recently composed by Greg Conner at Brigham Young University in 2025. The authors are also grateful to the referee(s) for  comments and suggestions that help to improve the presentation of the paper.

\end{document}